\documentclass[12pt]{article}
\usepackage{scicite}

\usepackage{times}

\newenvironment{sciabstract}{%
\begin{quote} \bf}
{\end{quote}}

\usepackage{graphicx}
\usepackage{amssymb}
\usepackage{amsmath}
\usepackage{textcomp}
\usepackage{setspace}
\usepackage{xcolor}
\usepackage{tikz}
\usetikzlibrary{decorations.pathreplacing,decorations.markings,shapes}
\usetikzlibrary{arrows,scopes}
\usepackage{xfrac}
\usepackage{float}
\usepackage{comment}
\usepackage{lineno}
\usepackage{ dsfont }

\newcommand{\todoprivate}[1]{}
\newcommand{\todoscience}[1]{}

\newcommand{\removecite}[1]{}

\newcommand{\todocite}[1]{\textcolor{green}{}}
\newcommand{\forthesis}[1]{\textcolor{blue}{}}

\usepackage{hyperref}
\usepackage{graphicx,dblfloatfix}
\usepackage{mathrsfs}
\usepackage{eufrak}
\usepackage{yfonts}
\usepackage{changepage}

\usepackage{color}
\definecolor{womencolor}{RGB}{84,39,143}
\definecolor{mencolor}{RGB}{0,109,44}
\definecolor{unknowncolor}{RGB}{230,85,13}

\newcommand{\tikzEdgeWidth}{0.667}
\newcommand{\tikzBaseline}{0.5}
\newcommand{\tikzNodeSize}{0.15}
\newcommand{\directedArrowSize}{0.275}
\newcommand{\directedArrowSizeBig}{0.45}
\newcommand{\directedEdgeWidth}{0.667}
\newcommand{\tikzBaselineNodeLocal}{1.333} 
\newcommand{\tikzEdgeWidthBig}{1.0}
\newcommand{\tikzBaselineBig}{1.25}

\newcommand{\tikzSize}{0.2}
\newcommand{\tikzSizeBig}{0.3}
\newcommand{\groupG}{\mathbb{G}}
\newcommand{\SIreftext}{supplementary materials}
\newcommand{\SIreftextInSI}{supplementary materials}
\newcommand{\figreftext}{fig.}
\newcommand{\figreftextSI}{fig.}
\newcommand{\figreftexts}{figs.}
\newcommand{\figreftextSIs}{figs.}
\newcommand{\groupg}{g}
\newcommand{\graphG}{G}
\newcommand{\graphGrv}{\mathcal{G}}

\newcommand{\oneedge}{\begin{tikzpicture}[scale=\tikzSize,baseline=\tikzBaseline]
\draw[rounded corners=0.025em,round cap-round cap,line width=\tikzEdgeWidth] (0.0,0.0) -- (0.5,1.0);
\end{tikzpicture}}
\newcommand{\twowedge}{\begin{tikzpicture}[scale=\tikzSize,baseline=\tikzBaseline]
\draw[rounded corners=0.025em,round cap-round cap,line width=\tikzEdgeWidth] (0.0,0.0) -- (0.5,1.0) -- (1.0,0.0);
\end{tikzpicture}}
\newcommand{\twoparallel}{\begin{tikzpicture}[scale=\tikzSize,baseline=\tikzBaseline]
\draw[rounded corners=0.025em,round cap-round cap,line width=\tikzEdgeWidth] (0.0,0.0) -- (0.5,1.0);
\draw[rounded corners=0.025em,-round cap,line width=\tikzEdgeWidth] (0.333,0.0) -- (0.833,1.0);
\end{tikzpicture}}
\newcommand{\threetriangle}{\begin{tikzpicture}[scale=\tikzSize,baseline=\tikzBaseline]
\draw[rounded corners=0.025em,round cap-round cap,line width=\tikzEdgeWidth] (0.0,0.0) -- (0.5,1.0) -- (1.0,0.0) -- (0.0,0.0) -- (0.5,1.0);
\end{tikzpicture}}
\newcommand{\threeclaw}{\begin{tikzpicture}[scale=\tikzSize,baseline=\tikzBaseline]
\draw[rounded corners=0.025em,round cap-round cap,line width=\tikzEdgeWidth] (0.0,0.0) -- (0.5,0.333) -- (0.5,1.0);
\draw[rounded corners=0.025em,round cap-round cap,line width=\tikzEdgeWidth] (0.5,1.0) -- (0.5,0.333) -- (1.0,0.0);
\draw[rounded corners=0.025em,round cap-round cap,line width=\tikzEdgeWidth] (1.0,0.0) -- (0.5,0.333) -- (0.0,0.0);
\end{tikzpicture}}
\newcommand{\threeline}{\hspace{0.025em}\begin{tikzpicture}[scale=\tikzSize,baseline=\tikzBaseline]
\draw[rounded corners=0.025em,round cap-round cap,line width=\tikzEdgeWidth] (0.0,0.0) -- (0.0,1.0) -- (1.0,1.0) -- (1.0,0.0);
\end{tikzpicture}}
\newcommand{\threeedgewedge}{\begin{tikzpicture}[scale=\tikzSize,baseline=\tikzBaseline]
\draw[rounded corners=0.025em,round cap-round cap,line width=\tikzEdgeWidth] (0.0,0.0) -- (0.5,1.0);
\draw[rounded corners=0.025em,-round cap,line width=\tikzEdgeWidth] (0.333,0.0) -- (0.833,1.0) -- (1.333,0.0);
\end{tikzpicture}}
\newcommand{\threeparallel}{\begin{tikzpicture}[scale=\tikzSize,baseline=\tikzBaseline]
\draw[rounded corners=0.025em,round cap-round cap,line width=\tikzEdgeWidth] (0.0,0.0) -- (0.5,1.0);
\draw[rounded corners=0.025em,-round cap,line width=\tikzEdgeWidth] (0.333,0.0) -- (0.833,1.0);
\draw[rounded corners=0.025em,round cap-round cap,line width=\tikzEdgeWidth] (0.667,0.0) -- (1.167,1.0);
\end{tikzpicture}}

\newcommand{\fourtriangleedge}{\begin{tikzpicture}[scale=\tikzSize,baseline=\tikzBaseline]
\draw[rounded corners=0.025em,round cap-round cap,line width=\tikzEdgeWidth] (0.0,0.0) -- (0.5,1.0) -- (1.0,0.0) -- (0.0,0.0) -- (0.5,1.0) -- (1.0,0.0) -- (1.5,1.0);
\draw[rounded corners=0.025em,round cap-round cap,line width=\tikzEdgeWidth] (1.5,1.0) -- (1.0,0.0) -- (0.0,0.0);
\end{tikzpicture}}
\newcommand{\foursquare}{\begin{tikzpicture}[scale=\tikzSize,baseline=\tikzBaseline]
\draw[rounded corners=0.025em,round cap-round cap,line width=\tikzEdgeWidth] (0.0,0.0) -- (1.0,0.0) -- (1.0,1.0) -- (0.0,1.0) -- (0.0,0.0) -- (1.0,0.0);
\end{tikzpicture}}

\newcommand{\fivediagsquare}{\begin{tikzpicture}[scale=\tikzSize,baseline=\tikzBaseline]
\draw[rounded corners=0.025em,round cap-round cap,line width=\tikzEdgeWidth] (0.0,0.0) -- (1.0,0.0) -- (1.0,1.0) -- (0.0,1.0) -- (0.0,0.0)  -- (1.0,1.0) -- (0.0,1.0) -- (0.0,0.0) -- (1.0,0.0) -- (1.0,1.0) -- (0.0,0.0) -- (1.0,0.0);
\end{tikzpicture}}
\newcommand{\sixKfour}{\begin{tikzpicture}[scale=\tikzSize,baseline=\tikzBaseline]
\draw[rounded corners=0.025em,round cap-round cap,line width=\tikzEdgeWidth] (0.0,0.0) -- (1.0,0.0) -- (1.0,1.0) -- (0.0,1.0) -- (0.0,0.0)  -- (1.0,1.0) -- (0.0,1.0) -- (0.0,0.0) -- (1.0,0.0) -- (1.0,1.0) -- (0.0,0.0) -- (1.0,0.0) -- (0.0,1.0) -- (0.0,0.0) -- (1.0,0.0) -- (1.0,1.0) -- (0.0,1.0) -- (1.0,0.0) -- (1.0,1.0);
\end{tikzpicture}}

\newcommand{\tikzNodeLocalNodeSize}{\tikzNodeSize}
\newcommand{\oneedgelocalnodeself}{\begin{tikzpicture}[scale=\tikzSize,baseline=\tikzBaselineNodeLocal]
\draw[rounded corners=0.025em,round cap-round cap,line width=\tikzEdgeWidth] (0.0,0.0) -- (0.5,1.0);
\path [draw=black,fill=white,line width=\tikzEdgeWidth] (0.5,1.0) circle (\tikzNodeLocalNodeSize);
\end{tikzpicture}}
\newcommand{\oneedgelocalnodeother}{\begin{tikzpicture}[scale=\tikzSize,baseline=\tikzBaselineNodeLocal]
\draw[rounded corners=0.025em,round cap-round cap,line width=\tikzEdgeWidth] (0.0,0.0) -- (1.0,0.0);
\path [draw=black,fill=white,line width=\tikzEdgeWidth] (0.5,1.0) circle (\tikzNodeLocalNodeSize);
\end{tikzpicture}}
\newcommand{\twowedgelocalnodecenter}{\begin{tikzpicture}[scale=\tikzSize,baseline=\tikzBaselineNodeLocal]
\draw[rounded corners=0.025em,round cap-round cap,line width=\tikzEdgeWidth] (0.0,0.0) -- (0.5,1.0) -- (1.0,0.0);
\path [draw=black,fill=white,line width=\tikzEdgeWidth] (0.5,1.0) circle (\tikzNodeLocalNodeSize);
\end{tikzpicture}}
\newcommand{\twowedgelocalnodeend}{\begin{tikzpicture}[scale=\tikzSize,baseline=\tikzBaselineNodeLocal]
\draw[rounded corners=0.025em,round cap-round cap,line width=\tikzEdgeWidth] (0.5,1.0) -- (0.0,0.0) -- (1.0,0.0);
\path [draw=black,fill=white,line width=\tikzEdgeWidth] (0.5,1.0) circle (\tikzNodeLocalNodeSize);
\end{tikzpicture}}
\newcommand{\threetrianglelocalnode}{\begin{tikzpicture}[scale=\tikzSize,baseline=\tikzBaselineNodeLocal]
\draw[rounded corners=0.025em,round cap-round cap,line width=\tikzEdgeWidth] (0.0,0.0) -- (0.5,1.0) -- (1.0,0.0) -- (0.0,0.0) -- (0.5,1.0);
\path [draw=black,fill=white,line width=\tikzEdgeWidth] (0.5,1.0) circle (\tikzNodeLocalNodeSize);
\end{tikzpicture}}

\newcommand{\tikzBaselineEdgeLocal}{\tikzBaseline}
\newcommand{\tikzFourStarRatio}{3.5}
\newcommand{\tikzFourStarSize}{0.09}
\newcommand{\oneedgelocaledge}{\begin{tikzpicture}[scale=\tikzSize,baseline=\tikzBaselineEdgeLocal]
\draw[rounded corners=0.025em,round cap-round cap,line width=\tikzEdgeWidth] (0.0,0.0) -- (0.5,1.0);
\draw (0.25,0.5) node[star,fill=black,scale=\tikzFourStarSize,rotate=63.435,star points=4,star point ratio=\tikzFourStarRatio]{};
\end{tikzpicture}}
\newcommand{\twowedgelocaledge}{\begin{tikzpicture}[scale=\tikzSize,baseline=\tikzBaselineEdgeLocal]
\draw[rounded corners=0.025em,round cap-round cap,line width=\tikzEdgeWidth] (0.0,0.0) -- (0.5,1.0) -- (1.0,0.0);
\draw (0.25,0.5) node[star,fill=black,scale=\tikzFourStarSize,rotate=63.435,star points=4,star point ratio=\tikzFourStarRatio]{};
\end{tikzpicture}}
\newcommand{\threetrianglelocaledge}{\begin{tikzpicture}[scale=\tikzSize,baseline=\tikzBaselineEdgeLocal]
\draw[rounded corners=0.025em,round cap-round cap,line width=\tikzEdgeWidth] (0.0,0.0) -- (0.5,1.0) -- (1.0,0.0) -- (0.0,0.0) -- (0.5,1.0);
\draw (0.25,0.5) node[star,fill=black,scale=\tikzFourStarSize,rotate=63.435,star points=4,star point ratio=\tikzFourStarRatio]{};
\end{tikzpicture}}
\newcommand{\oneedgedetachedlocaledge}{\begin{tikzpicture}[scale=\tikzSize,baseline=\tikzBaselineEdgeLocal]
\draw[rounded corners=0.025em,round cap-round cap,line width=\tikzEdgeWidth]  (0.5,1.0) -- (1.0,0.0);
\draw (0.25,0.5) node[star,fill=black,scale=\tikzFourStarSize,rotate=63.435,star points=4,star point ratio=\tikzFourStarRatio]{};
\end{tikzpicture}}
\newcommand{\twowedgedetachedlocaledge}{\begin{tikzpicture}[scale=\tikzSize,baseline=\tikzBaselineEdgeLocal]
\draw[rounded corners=0.025em,round cap-round cap,line width=\tikzEdgeWidth]  (0.5,1.0) -- (1.0,0.0) -- (0.0,0.0);
\draw (0.25,0.5) node[star,fill=black,scale=\tikzFourStarSize,rotate=63.435,star points=4,star point ratio=\tikzFourStarRatio]{};
\end{tikzpicture}}

\newcommand{\tikzEdgeWidthDirected}{\directedEdgeWidth}
\newcommand{\tikzBaselineDirected}{\tikzBaseline}
\newcommand{\fourdirecttriangle}{\begin{tikzpicture}[scale=\tikzSize,baseline=\tikzBaselineDirected]
\draw[rounded corners=0.025em,round cap-round cap,line width=\tikzEdgeWidthDirected] (0.0,0.0) -- (0.5,1.0) -- (1.0,0.0) arc (180+123.69006752597979:180+56.309932474020215:0.9013878188659973) arc (123.69006752597979:56.309932474020215:0.9013878188659973) arc (180+123.69006752597979:180+56.309932474020215:0.9013878188659973) -- (0.5,1.0);
\end{tikzpicture}}

\newcommand{\tworeciprocal}{\begin{tikzpicture}[scale=\tikzSize,baseline=\tikzBaselineDirected]
\draw[rounded corners=0.025em,round cap-round cap,line width=\tikzEdgeWidthDirected] (0.5,1.0) arc (123.12922234117626:183.7406753046678:1.1078112962242475) arc (303.1292223411763:360+3.7406753046677887:1.1078112962242475) arc (123.12922234117626:183.7406753046678:1.1078112962242475);
\end{tikzpicture}}

\newcommand{\threereciprocal}{\begin{tikzpicture}[scale=\tikzSize,baseline=\tikzBaselineDirected]
\draw[rounded corners=0.025em,round cap-round cap,line width=\tikzEdgeWidthDirected] (0.5,1.0) arc (123.12922234117626:183.7406753046678:1.1078112962242475) arc (303.1292223411763:360+3.7406753046677887:1.1078112962242475) arc (123.12922234117626:183.7406753046678:1.1078112962242475);
\draw[rounded corners=0.025em,round cap-round cap,line width=\tikzEdgeWidthDirected] (0.0,0.0) -- (0.5,1.0);
\end{tikzpicture}}

\newcommand{\threereciprocalwedge}{\begin{tikzpicture}[scale=\tikzSize,baseline=\tikzBaselineDirected]
\draw[rounded corners=0.025em,round cap-round cap,line width=\tikzEdgeWidthDirected] (1.0,0.0) -- (0.5,1.0) arc (123.12922234117626:183.7406753046678:1.1078112962242475) arc (303.1292223411763:360+3.7406753046677887:1.1078112962242475) arc (123.12922234117626:183.7406753046678:1.1078112962242475);
\end{tikzpicture}}

\tikzset{
  on each segment/.style={
    decorate,
    decoration={
      show path construction,
      moveto code={},
      lineto code={
        \path [#1]
        (\tikzinputsegmentfirst) -- (\tikzinputsegmentlast);
      },
      curveto code={
        \path [#1] (\tikzinputsegmentfirst)
        .. controls
        (\tikzinputsegmentsupporta) and (\tikzinputsegmentsupportb)
        ..
        (\tikzinputsegmentlast);
      },
      closepath code={
        \path [#1]
        (\tikzinputsegmentfirst) -- (\tikzinputsegmentlast);
      },
    },
  },
  mid arrow/.style={postaction={decorate,decoration={
        markings,
        mark=at position #1 with {\arrow[black,scale=\directedArrowSize]{triangle 60}}
      }}},
  mid arrow big/.style={postaction={decorate,decoration={
        markings,
        mark=at position #1 with {\arrow[black,scale=\directedArrowSizeBig]{triangle 60}}
      }}},
}

\newcommand{\oneedgedir}{\begin{tikzpicture}[scale=\tikzSize,baseline=\tikzBaselineDirected]
\path[draw=white,postaction={on each segment={mid arrow=0.85}}] 
(0.0,0.0) -- (0.5,1.0);
\draw[rounded corners=0.025em,round cap-round cap,line width=\tikzEdgeWidthDirected] 
(0.0,0.0) -- (0.5,1.0);
\end{tikzpicture}}

\newcommand{\twoparalleldir}{\begin{tikzpicture}[scale=\tikzSize,baseline=\tikzBaselineDirected]
\path[draw=white,postaction={on each segment={mid arrow=0.85}}] 
(0.0,0.0) -- (0.5,1.0);
\draw[rounded corners=0.025em,round cap-round cap,line width=\tikzEdgeWidthDirected] 
(0.0,0.0) -- (0.5,1.0);
\path[draw=white,postaction={on each segment={mid arrow=0.85}}] 
(0.333,0.0) -- (0.833,1.0);
\draw[rounded corners=0.025em,round cap-round cap,line width=\tikzEdgeWidthDirected] 
(0.333,0.0) -- (0.833,1.0);
\end{tikzpicture}}

\newcommand{\tworeciprocaldir}{\begin{tikzpicture}[scale=\tikzSize,baseline=\tikzBaselineDirected]
\path[draw=white,postaction={on each segment={mid arrow=1.0}}] 
(0.0487578, 0.404908) -- (0.00246769, 0.175586);
\path[draw=white,postaction={on each segment={mid arrow=1.0}}] 
(0.294672, 0.281951) -- (0.480458, 0.708528);
\draw[rounded corners=0.025em,round cap-round cap,line width=\tikzEdgeWidthDirected] (0.5,1.0) arc (123.12922234117626:183.7406753046678:1.1078112962242475) arc (303.1292223411763:360+3.7406753046677887:1.1078112962242475) arc (123.12922234117626:183.7406753046678:1.1078112962242475);
\end{tikzpicture}}

\newcommand{\twowedgeii}{\begin{tikzpicture}[scale=\tikzSize,baseline=\tikzBaselineDirected]
\path[draw=white,postaction={on each segment={mid arrow=0.85}}] 
(0.0,0.0) -- (0.5,1.0);
\path[draw=white,postaction={on each segment={mid arrow=0.85}}] 
(1.0,0.0) -- (0.5,1.0);
\draw[rounded corners=0.025em,round cap-round cap,line width=\tikzEdgeWidthDirected] 
(0.0,0.0) -- (0.5,1.0) -- (1.0,0.0);
\end{tikzpicture}}

\newcommand{\twowedgeoo}{\begin{tikzpicture}[scale=\tikzSize,baseline=\tikzBaselineDirected]
\path[draw=white,postaction={on each segment={mid arrow=0.85}}] 
(0.5,1.0) -- (0.0,0.0);
\path[draw=white,postaction={on each segment={mid arrow=0.85}}] 
(0.5,1.0) -- (1.0,0.0);
\draw[rounded corners=0.025em,round cap-round cap,line width=\tikzEdgeWidthDirected] 
(0.0,0.0) -- (0.5,1.0) -- (1.0,0.0);
\end{tikzpicture}}

\newcommand{\twowedgeio}{\begin{tikzpicture}[scale=\tikzSize,baseline=\tikzBaselineDirected]
\path[draw=white,postaction={on each segment={mid arrow=0.85}}] 
(0.0,0.0) -- (0.5,1.0);
\path[draw=white,postaction={on each segment={mid arrow=0.85}}] 
(0.5,1.0) -- (1.0,0.0);
\draw[rounded corners=0.025em,round cap-round cap,line width=\tikzEdgeWidthDirected] 
(0.0,0.0) -- (0.5,1.0) -- (1.0,0.0);
\end{tikzpicture}}

\newcommand{\threereciprocalwedgei}{\begin{tikzpicture}[scale=\tikzSize,baseline=\tikzBaselineDirected]
\path[draw=white,postaction={on each segment={mid arrow=1.0}}] 
(0.0487578, 0.404908) -- (0.00246769, 0.175586);
\path[draw=white,postaction={on each segment={mid arrow=1.0}}] 
(0.294672, 0.281951) -- (0.480458, 0.708528);
\path[draw=white,postaction={on each segment={mid arrow=0.85}}] 
(1.0,0.0) -- (0.5,1.0);
\draw[rounded corners=0.025em,round cap-round cap,line width=\tikzEdgeWidthDirected] (1.0,0.0) -- (0.5,1.0) arc (123.12922234117626:183.7406753046678:1.1078112962242475) arc (303.1292223411763:360+3.7406753046677887:1.1078112962242475) arc (123.12922234117626:183.7406753046678:1.1078112962242475);
\end{tikzpicture}}

\newcommand{\threereciprocalwedgeo}{\begin{tikzpicture}[scale=\tikzSize,baseline=\tikzBaselineDirected]
\path[draw=white,postaction={on each segment={mid arrow=1.0}}] 
(0.0487578, 0.404908) -- (0.00246769, 0.175586);
\path[draw=white,postaction={on each segment={mid arrow=1.0}}] 
(0.294672, 0.281951) -- (0.480458, 0.708528);
\path[draw=white,postaction={on each segment={mid arrow=0.85}}] 
(0.5,1.0) -- (1.0,0.0);
\draw[rounded corners=0.025em,round cap-round cap,line width=\tikzEdgeWidthDirected] (1.0,0.0) -- (0.5,1.0) arc (123.12922234117626:183.7406753046678:1.1078112962242475) arc (303.1292223411763:360+3.7406753046677887:1.1078112962242475) arc (123.12922234117626:183.7406753046678:1.1078112962242475);
\end{tikzpicture}}

\newcommand{\threetrianglelinear}{\begin{tikzpicture}[scale=\tikzSize,baseline=\tikzBaselineDirected]
\path[draw=white,postaction={on each segment={mid arrow=0.85}}] 
(0.0,0.0) -- (0.5,1.0);
\path[draw=white,postaction={on each segment={mid arrow=0.85}}] 
(0.5,1.0) -- (1.0,0.0);
\path[draw=white,postaction={on each segment={mid arrow=0.85}}] 
(0.0,0.0) -- (1.0,0.0);
\draw[rounded corners=0.025em,round cap-round cap,line width=\tikzEdgeWidthDirected] (0.0,0.0) -- (0.5,1.0) -- (1.0,0.0) -- (0.0,0.0) -- (0.5,1.0);
\end{tikzpicture}}

\newcommand{\threetrianglecyclic}{\begin{tikzpicture}[scale=\tikzSize,baseline=\tikzBaselineDirected]
\path[draw=white,postaction={on each segment={mid arrow=0.85}}] 
(0.0,0.0) -- (0.5,1.0);
\path[draw=white,postaction={on each segment={mid arrow=0.85}}] 
(0.5,1.0) -- (1.0,0.0);
\path[draw=white,postaction={on each segment={mid arrow=0.85}}] 
(1.0,0.0) -- (0.0,0.0);
\draw[rounded corners=0.025em,round cap-round cap,line width=\tikzEdgeWidthDirected] (0.0,0.0) -- (0.5,1.0) -- (1.0,0.0) -- (0.0,0.0) -- (0.5,1.0);
\end{tikzpicture}}

\newcommand{\fourdirectedtriangleii}{\begin{tikzpicture}[scale=\tikzSize,baseline=\tikzBaselineDirected]
\path[draw=white,postaction={on each segment={mid arrow=0.85}}] 
(0.0,0.0) -- (0.5,1.0);
\path[draw=white,postaction={on each segment={mid arrow=0.85}}] 
(1.0,0.0) -- (0.5,1.0);
\path[draw=white,postaction={on each segment={mid arrow=1.0}}] 
(0.385497, -0.186057) -- (0.887606, -0.0657329);
\path[draw=white,postaction={on each segment={mid arrow=1.0}}] 
(0.744482, 0.178405) -- (0.231723, 0.117825);
\draw[rounded corners=0.025em,round cap-round cap,line width=\tikzEdgeWidthDirected] (0.0,0.0) -- (0.5,1.0) -- (1.0,0.0) arc (180+123.69006752597979:180+56.309932474020215:0.9013878188659973) arc (123.69006752597979:56.309932474020215:0.9013878188659973) arc (180+123.69006752597979:180+56.309932474020215:0.9013878188659973) -- (0.5,1.0);
\end{tikzpicture}}

\newcommand{\fourdirectedtriangleoo}{\begin{tikzpicture}[scale=\tikzSize,baseline=\tikzBaselineDirected]
\path[draw=white,postaction={on each segment={mid arrow=0.85}}] 
(0.5,1.0) -- (0.0,0.0);
\path[draw=white,postaction={on each segment={mid arrow=0.85}}] 
(0.5,1.0) -- (1.0,0.0);
\path[draw=white,postaction={on each segment={mid arrow=1.0}}] 
(0.385497, -0.186057) -- (0.887606, -0.0657329);
\path[draw=white,postaction={on each segment={mid arrow=1.0}}] 
(0.744482, 0.178405) -- (0.231723, 0.117825);
\draw[rounded corners=0.025em,round cap-round cap,line width=\tikzEdgeWidthDirected] (0.0,0.0) -- (0.5,1.0) -- (1.0,0.0) arc (180+123.69006752597979:180+56.309932474020215:0.9013878188659973) arc (123.69006752597979:56.309932474020215:0.9013878188659973) arc (180+123.69006752597979:180+56.309932474020215:0.9013878188659973) -- (0.5,1.0);
\end{tikzpicture}}

\newcommand{\fourdirectedtriangleio}{\begin{tikzpicture}[scale=\tikzSize,baseline=\tikzBaselineDirected]
\path[draw=white,postaction={on each segment={mid arrow=0.85}}] 
(0.0,0.0) -- (0.5,1.0);
\path[draw=white,postaction={on each segment={mid arrow=0.85}}] 
(0.5,1.0) -- (1.0,0.0);
\path[draw=white,postaction={on each segment={mid arrow=1.0}}] 
(0.385497, -0.186057) -- (0.887606, -0.0657329);
\path[draw=white,postaction={on each segment={mid arrow=1.0}}] 
(0.744482, 0.178405) -- (0.231723, 0.117825);
\draw[rounded corners=0.025em,round cap-round cap,line width=\tikzEdgeWidthDirected] (0.0,0.0) -- (0.5,1.0) -- (1.0,0.0) arc (180+123.69006752597979:180+56.309932474020215:0.9013878188659973) arc (123.69006752597979:56.309932474020215:0.9013878188659973) arc (180+123.69006752597979:180+56.309932474020215:0.9013878188659973) -- (0.5,1.0);
\end{tikzpicture}}

\newcommand{\fourreciprocalreciprocaldir}{\begin{tikzpicture}[scale=\tikzSize,baseline=\tikzBaselineDirected]
\path[draw=white,postaction={on each segment={mid arrow=1.0}}] 
(0.0487578, 0.404908) -- (0.00246769, 0.175586);
\path[draw=white,postaction={on each segment={mid arrow=1.0}}] 
(0.294672, 0.281951) -- (0.480458, 0.708528);
\path[draw=white,postaction={on each segment={mid arrow=1.0}}] 
(0.794672, 0.718049) -- (0.638988, 0.892674);
\path[draw=white,postaction={on each segment={mid arrow=1.0}}] 
(0.548758, 0.595092) -- (0.778547, 0.190517);
\draw[rounded corners=0.025em,round cap-round cap,line width=\tikzEdgeWidthDirected] (1.0,0.0) arc (356.2593246953322:360+56.870777658823755:1.1078112962242475) arc (176.25932469533223:236.87077765882376:1.1078112962242475) arc (356.2593246953322:360+56.870777658823755:1.1078112962242475) arc (123.12922234117626:183.7406753046678:1.1078112962242475) arc (303.1292223411763:360+3.7406753046677887:1.1078112962242475) arc (123.12922234117626:183.7406753046678:1.1078112962242475);
\end{tikzpicture}}

\newcommand{\fivedirectedtriangledir}{\begin{tikzpicture}[scale=\tikzSize,baseline=\tikzBaselineDirected]
\path[draw=white,postaction={on each segment={mid arrow=1.0}}] 
(0.0487578, 0.404908) -- (0.00246769, 0.175586);
\path[draw=white,postaction={on each segment={mid arrow=1.0}}] 
(0.294672, 0.281951) -- (0.480458, 0.708528);
\path[draw=white,postaction={on each segment={mid arrow=1.0}}] 
(0.794672, 0.718049) -- (0.638988, 0.892674);
\path[draw=white,postaction={on each segment={mid arrow=1.0}}] 
(0.548758, 0.595092) -- (0.778547, 0.190517);
\path[draw=white,postaction={on each segment={mid arrow=0.85}}] 
(0.0,0.0) -- (1.0,0.0);
\draw[rounded corners=0.025em,round cap-round cap,line width=\tikzEdgeWidthDirected] (0.5,1.0) arc (123.12922234117626:183.7406753046678:1.1078112962242475) arc (303.1292223411763:360+3.7406753046677887:1.1078112962242475) arc (123.12922234117626:183.7406753046678:1.1078112962242475) -- (1.0,0.0) arc (356.2593246953322:360+56.870777658823755:1.1078112962242475) arc (176.25932469533223:236.87077765882376:1.1078112962242475) arc (356.2593246953322:360+56.870777658823755:1.1078112962242475) arc (123.12922234117626:183.7406753046678:1.1078112962242475);
\end{tikzpicture}}

\newcommand{\sixdirectedtriangledir}{\begin{tikzpicture}[scale=\tikzSize,baseline=\tikzBaselineDirected]
\path[draw=white,postaction={on each segment={mid arrow=1.0}}] 
(0.0487578, 0.404908) -- (0.00246769, 0.175586);
\path[draw=white,postaction={on each segment={mid arrow=1.0}}] 
(0.294672, 0.281951) -- (0.480458, 0.708528);
\path[draw=white,postaction={on each segment={mid arrow=1.0}}] 
(0.794672, 0.718049) -- (0.638988, 0.892674);
\path[draw=white,postaction={on each segment={mid arrow=1.0}}] 
(0.548758, 0.595092) -- (0.778547, 0.190517);
\path[draw=white,postaction={on each segment={mid arrow=1.0}}] 
(0.385497, -0.186057) -- (0.887606, -0.0657329);
\path[draw=white,postaction={on each segment={mid arrow=1.0}}] 
(0.744482, 0.178405) -- (0.231723, 0.117825);
\draw[rounded corners=0.025em,round cap-round cap,line width=\tikzEdgeWidthDirected] (0.5,1.0) arc (123.12922234117626:183.7406753046678:1.1078112962242475) arc (303.1292223411763:360+3.7406753046677887:1.1078112962242475) arc (123.12922234117626:183.7406753046678:1.1078112962242475) arc (180+56.309932474020215:180+123.69006752597979:0.9013878188659973) arc (56.309932474020215:123.69006752597979:0.9013878188659973) arc (180+56.309932474020215:180+123.69006752597979:0.9013878188659973) arc (356.2593246953322:360+56.870777658823755:1.1078112962242475) arc (176.25932469533223:236.87077765882376:1.1078112962242475) arc (356.2593246953322:360+56.870777658823755:1.1078112962242475) arc (123.12922234117626:183.7406753046678:1.1078112962242475);
\end{tikzpicture}}

\newcommand{\nodeAttributeScale}{0.2}
\newcommand{\nodeAttributeBaseline}{\tikzBaselineNodeLocal}
\newcommand{\nodeAttributeVertexSize}{0.2}
\newcommand{\nodeAttributeEdgeWidth}{\tikzEdgeWidth}
\newcommand{\colorPurple}{womencolor}
\newcommand{\colorGreen}{mencolor}

\newcommand{\PurpleDot}{\begin{tikzpicture}[baseline=\nodeAttributeBaseline,scale=\nodeAttributeScale]
\filldraw[ultra thin,color=\colorPurple] (0,0.5) circle (\nodeAttributeVertexSize);
\end{tikzpicture}}

\newcommand{\GreenDot}{\begin{tikzpicture}[baseline=\nodeAttributeBaseline,scale=\nodeAttributeScale]
\filldraw[ultra thin,color=\colorGreen] (0,0.5) circle (\nodeAttributeVertexSize);
\end{tikzpicture}}

\newcommand{\OnePurpleToGreenEdge}{\begin{tikzpicture}[baseline=\nodeAttributeBaseline,scale=\nodeAttributeScale]
\draw[rounded corners=0.025em,color=black, line width=\nodeAttributeEdgeWidth, round cap-round cap] (0.0,0.0) -- (0.5,1.0); 
\filldraw[ultra thin,color=\colorPurple] (0,0) circle (\nodeAttributeVertexSize);
\filldraw[ultra thin,color=\colorGreen] (0.5,1.0) circle (\nodeAttributeVertexSize);
\end{tikzpicture}}

\newcommand{\OnePurpleToPurpleEdge}{\begin{tikzpicture}[baseline=\nodeAttributeBaseline,scale=\nodeAttributeScale]
\draw[rounded corners=0.025em,color=black, line width=\nodeAttributeEdgeWidth, round cap-round cap] (0.0,0.0) -- (0.5,1.0);
\filldraw[ultra thin,color=\colorPurple] (0,0) circle (\nodeAttributeVertexSize);
\filldraw[ultra thin,color=\colorPurple] (0.5,1.0) circle (\nodeAttributeVertexSize);
\end{tikzpicture}}

\newcommand{\OneGreenToGreenEdge}{\begin{tikzpicture}[baseline=\nodeAttributeBaseline,scale=\nodeAttributeScale]
\draw[rounded corners=0.025em,color=black, line width=\nodeAttributeEdgeWidth, round cap-round cap] (0.0,0.0) -- (0.5,1.0);
\filldraw[ultra thin,color=\colorGreen] (0,0) circle (\nodeAttributeVertexSize);
\filldraw[ultra thin,color=\colorGreen] (0.5,1.0) circle (\nodeAttributeVertexSize);
\end{tikzpicture}}

\newcommand{\TwoPurplePurplePurple}{\begin{tikzpicture}[baseline=\nodeAttributeBaseline,scale=\nodeAttributeScale]
\draw[rounded corners=0.025em,color=black, line width=\nodeAttributeEdgeWidth, round cap-round cap] (0.0,0.0) -- (0.5,1.0);
\draw[rounded corners=0.025em,color=black, line width=\nodeAttributeEdgeWidth, round cap-round cap] (0.5,1.0) -- (1.0,0.0);
\filldraw[ultra thin,color=\colorPurple] (0,0) circle (\nodeAttributeVertexSize);
\filldraw[ultra thin,color=\colorPurple] (0.5,1.0) circle (\nodeAttributeVertexSize);
\filldraw[ultra thin,color=\colorPurple] (1.0,0.0) circle (\nodeAttributeVertexSize);
\end{tikzpicture}}

\newcommand{\TwoPurplePurpleGreen}{\begin{tikzpicture}[baseline=\nodeAttributeBaseline,scale=\nodeAttributeScale]
\draw[rounded corners=0.025em,color=black, line width=\nodeAttributeEdgeWidth, round cap-round cap] (0.0,0.0) -- (0.5,1.0);
\draw[rounded corners=0.025em,color=black, line width=\nodeAttributeEdgeWidth, round cap-round cap] (0.5,1.0) -- (1.0,0.0);
\filldraw[ultra thin,color=\colorPurple] (0,0) circle (\nodeAttributeVertexSize);
\filldraw[ultra thin,color=\colorPurple] (0.5,1.0) circle (\nodeAttributeVertexSize);
\filldraw[ultra thin,color=\colorGreen] (1.0,0.0) circle (\nodeAttributeVertexSize);
\end{tikzpicture}}

\newcommand{\TwoGreenPurpleGreen}{\begin{tikzpicture}[baseline=\nodeAttributeBaseline,scale=\nodeAttributeScale]
\draw[rounded corners=0.025em,color=black, line width=\nodeAttributeEdgeWidth, round cap-round cap] (0.0,0.0) -- (0.5,1.0);
\draw[rounded corners=0.025em,color=black, line width=\nodeAttributeEdgeWidth, round cap-round cap] (0.5,1.0) -- (1.0,0.0);
\filldraw[ultra thin,color=\colorGreen] (0,0) circle (\nodeAttributeVertexSize);
\filldraw[ultra thin,color=\colorPurple] (0.5,1.0) circle (\nodeAttributeVertexSize);
\filldraw[ultra thin,color=\colorGreen] (1.0,0.0) circle (\nodeAttributeVertexSize);
\end{tikzpicture}}

\newcommand{\TwoPurpleGreenPurple}{\begin{tikzpicture}[baseline=\nodeAttributeBaseline,scale=\nodeAttributeScale]
\draw[rounded corners=0.025em,color=black, line width=\nodeAttributeEdgeWidth, round cap-round cap] (0.0,0.0) -- (0.5,1.0);
\draw[rounded corners=0.025em,color=black, line width=\nodeAttributeEdgeWidth, round cap-round cap] (0.5,1.0) -- (1.0,0.0);
\filldraw[ultra thin,color=\colorPurple] (0,0) circle (\nodeAttributeVertexSize);
\filldraw[ultra thin,color=\colorGreen] (0.5,1.0) circle (\nodeAttributeVertexSize);
\filldraw[ultra thin,color=\colorPurple] (1.0,0.0) circle (\nodeAttributeVertexSize);
\end{tikzpicture}}

\newcommand{\ThreePurplePurplePurple}{\begin{tikzpicture}[baseline=\nodeAttributeBaseline,scale=\nodeAttributeScale]
\draw[rounded corners=0.025em,color=black, line width=\nodeAttributeEdgeWidth, round cap-round cap] (0.0,0.0) -- (0.5,1.0);
\draw[rounded corners=0.025em,color=black, line width=\nodeAttributeEdgeWidth, round cap-round cap] (0.5,1.0) -- (1.0,0.0);
\draw[rounded corners=0.025em,color=black, line width=\nodeAttributeEdgeWidth, round cap-round cap] (0.0,0.0) -- (1.0,0.0);
\filldraw[ultra thin,color=\colorPurple] (0,0) circle (\nodeAttributeVertexSize);
\filldraw[ultra thin,color=\colorPurple] (0.5,1.0) circle (\nodeAttributeVertexSize);
\filldraw[ultra thin,color=\colorPurple] (1.0,0.0) circle (\nodeAttributeVertexSize);
\end{tikzpicture}}

\newcommand{\ThreePurplePurpleGreen}{\begin{tikzpicture}[baseline=\nodeAttributeBaseline,scale=\nodeAttributeScale]
\draw[rounded corners=0.025em,color=black, line width=\nodeAttributeEdgeWidth, round cap-round cap] (0.0,0.0) -- (0.5,1.0);
\draw[rounded corners=0.025em,color=black, line width=\nodeAttributeEdgeWidth, round cap-round cap] (0.5,1.0) -- (1.0,0.0);
\draw[rounded corners=0.025em,color=black, line width=\nodeAttributeEdgeWidth, round cap-round cap] (0.0,0.0) -- (1.0,0.0);
\filldraw[ultra thin,color=\colorPurple] (0,0) circle (\nodeAttributeVertexSize);
\filldraw[ultra thin,color=\colorPurple] (0.5,1.0) circle (\nodeAttributeVertexSize);
\filldraw[ultra thin,color=\colorGreen] (1.0,0.0) circle (\nodeAttributeVertexSize);
\end{tikzpicture}}

\newcommand{\ThreeClawPurplePurplePurple}{\begin{tikzpicture}[baseline=\nodeAttributeBaseline,scale=\nodeAttributeScale]
\draw[rounded corners=0.025em,color=black, line width=\nodeAttributeEdgeWidth, round cap-round cap] (0.0,0.0) -- (0.5,0.333);
\draw[rounded corners=0.025em,color=black, line width=\nodeAttributeEdgeWidth, round cap-round cap] (0.5,1.0) -- (0.5,0.333);
\draw[rounded corners=0.025em,color=black, line width=\nodeAttributeEdgeWidth, round cap-round cap] (1.0,0.0) -- (0.5,0.333);
\filldraw[ultra thin,color=\colorPurple] (0.5,0.333) circle (\nodeAttributeVertexSize);
\filldraw[ultra thin,color=\colorPurple] (0,0) circle (\nodeAttributeVertexSize);
\filldraw[ultra thin,color=\colorPurple] (0.5,1.0) circle (\nodeAttributeVertexSize);
\filldraw[ultra thin,color=\colorPurple] (1.0,0.0) circle (\nodeAttributeVertexSize);
\end{tikzpicture}}

\newcommand{\ThreeClawPurplePurpleGreen}{\begin{tikzpicture}[baseline=\nodeAttributeBaseline,scale=\nodeAttributeScale]
\draw[rounded corners=0.025em,color=black, line width=\nodeAttributeEdgeWidth, round cap-round cap] (0.0,0.0) -- (0.5,0.333);
\draw[rounded corners=0.025em,color=black, line width=\nodeAttributeEdgeWidth, round cap-round cap] (0.5,1.0) -- (0.5,0.333);
\draw[rounded corners=0.025em,color=black, line width=\nodeAttributeEdgeWidth, round cap-round cap] (1.0,0.0) -- (0.5,0.333);
\filldraw[ultra thin,color=\colorPurple] (0.5,0.333) circle (\nodeAttributeVertexSize);
\filldraw[ultra thin,color=\colorPurple] (0,0) circle (\nodeAttributeVertexSize);
\filldraw[ultra thin,color=\colorPurple] (0.5,1.0) circle (\nodeAttributeVertexSize);
\filldraw[ultra thin,color=\colorGreen] (1.0,0.0) circle (\nodeAttributeVertexSize);
\end{tikzpicture}}

\newcommand{\ThreeClawPurpleGreenGreen}{\begin{tikzpicture}[baseline=\nodeAttributeBaseline,scale=\nodeAttributeScale]
\draw[rounded corners=0.025em,color=black, line width=\nodeAttributeEdgeWidth, round cap-round cap] (0.0,0.0) -- (0.5,0.333);
\draw[rounded corners=0.025em,color=black, line width=\nodeAttributeEdgeWidth, round cap-round cap] (0.5,1.0) -- (0.5,0.333);
\draw[rounded corners=0.025em,color=black, line width=\nodeAttributeEdgeWidth, round cap-round cap] (1.0,0.0) -- (0.5,0.333);
\filldraw[ultra thin,color=\colorPurple] (0.5,0.333) circle (\nodeAttributeVertexSize);
\filldraw[ultra thin,color=\colorPurple] (0,0) circle (\nodeAttributeVertexSize);
\filldraw[ultra thin,color=\colorGreen] (0.5,1.0) circle (\nodeAttributeVertexSize);
\filldraw[ultra thin,color=\colorGreen] (1.0,0.0) circle (\nodeAttributeVertexSize);
\end{tikzpicture}}

\newcommand{\ThreeClawGreenGreenGreen}{\begin{tikzpicture}[baseline=\nodeAttributeBaseline,scale=\nodeAttributeScale]
\draw[rounded corners=0.025em,color=black, line width=\nodeAttributeEdgeWidth, round cap-round cap] (0.0,0.0) -- (0.5,0.333);
\draw[rounded corners=0.025em,color=black, line width=\nodeAttributeEdgeWidth, round cap-round cap] (0.5,1.0) -- (0.5,0.333);
\draw[rounded corners=0.025em,color=black, line width=\nodeAttributeEdgeWidth, round cap-round cap] (1.0,0.0) -- (0.5,0.333);
\filldraw[ultra thin,color=\colorPurple] (0.5,0.333) circle (\nodeAttributeVertexSize);
\filldraw[ultra thin,color=\colorGreen] (0,0) circle (\nodeAttributeVertexSize);
\filldraw[ultra thin,color=\colorGreen] (0.5,1.0) circle (\nodeAttributeVertexSize);
\filldraw[ultra thin,color=\colorGreen] (1.0,0.0) circle (\nodeAttributeVertexSize);
\end{tikzpicture}}

\newcommand{\ThreeTriangleLinearBig}{\begin{tikzpicture}[baseline=\tikzBaselineBig,scale=\tikzSizeBig]
\path[draw=white,postaction={on each segment={mid arrow big=0.85}}] 
(0.0,0.0) -- (0.5,1.0);
\path[draw=white,postaction={on each segment={mid arrow big=0.85}}] 
(0.5,1.0) -- (1.0,0.0);
\path[draw=white,postaction={on each segment={mid arrow big=0.85}}] 
(0.0,0.0) -- (1.0,0.0);
\draw[rounded corners=0.025em,round cap-round cap,line width=\tikzEdgeWidthBig] (0.0,0.0) -- (0.5,1.0) -- (1.0,0.0) -- (0.0,0.0) -- (0.5,1.0);
\end{tikzpicture}}

\newcommand{\FourBifanBig}{\begin{tikzpicture}[baseline=\tikzBaselineBig,scale=\tikzSizeBig]
\path[draw=white,postaction={on each segment={mid arrow big=0.85}}] 
(0.0,0.0) -- (1.0,0.0);
\path[draw=white,postaction={on each segment={mid arrow big=0.85}}] 
(0.0,0.0) -- (0.0,1.0);
\path[draw=white,postaction={on each segment={mid arrow big=0.85}}] 
(1.0,1.0) -- (1.0,0.0);
\path[draw=white,postaction={on each segment={mid arrow big=0.85}}] 
(1.0,1.0) -- (0.0,1.0);
\draw[rounded corners=0.025em,round cap-round cap,line width=\tikzEdgeWidthBig] (0.0,0.0) -- (1.0,0.0) -- (1.0,1.0) -- (0.0,1.0) -- (0.0,0.0) -- (1.0,0.0);
\end{tikzpicture}}

\newcommand{\FourReciprocalReciprocalDirBig}{\begin{tikzpicture}[baseline=\tikzBaselineBig,scale=\tikzSizeBig]
\path[draw=white,postaction={on each segment={mid arrow big=1.0}}] 
(0.0487578, 0.404908) -- (0.00246769, 0.175586);
\path[draw=white,postaction={on each segment={mid arrow big=1.0}}] 
(0.294672, 0.281951) -- (0.480458, 0.708528);
\path[draw=white,postaction={on each segment={mid arrow big=1.0}}] 
(0.794672, 0.718049) -- (0.638988, 0.892674);
\path[draw=white,postaction={on each segment={mid arrow big=1.0}}] 
(0.548758, 0.595092) -- (0.778547, 0.190517);
\draw[rounded corners=0.025em,round cap-round cap,line width=\tikzEdgeWidthBig] (1.0,0.0) arc (356.2593246953322:360+56.870777658823755:1.1078112962242475) arc (176.25932469533223:236.87077765882376:1.1078112962242475) arc (356.2593246953322:360+56.870777658823755:1.1078112962242475) arc (123.12922234117626:183.7406753046678:1.1078112962242475) arc (303.1292223411763:360+3.7406753046677887:1.1078112962242475) arc (123.12922234117626:183.7406753046678:1.1078112962242475);
\end{tikzpicture}}

\usepackage[figurename=Fig.,labelfont=bf]{caption}

\newcommand{\numsymbol}{\raisebox{-1.25pt}{\scalebox{1}{\#}}}

\usepackage[numbers,square]{natbib}
\usepackage{hyperref,url}  

\title{\vspace{-24pt}Introducing Graph Cumulants:\\What is the Variance of Your Social Network?} 

\author
{Lee M.~Gunderson$^{\ast,1,\dagger}$ \& Gecia Bravo-Hermsdorff$^{\;\!\:\!\ast,2,\dagger}$ \\
\normalsize{$^\ast$\textit{Both authors contributed equally to this work}}\\
\\
\normalsize{$^{1}$Department of Astrophysical Sciences, Princeton University, Princeton, NJ, 08544, USA}\\
\normalsize{$^{2}$Princeton Neuroscience Institute, Princeton University, Princeton, NJ, 08544, USA}\\
\\
\normalsize{$^\dagger$Correspondence to: geciah@princeton.edu or leeg@princeton.edu}
}

\date{}

\begin{document} 


\baselineskip24pt


\maketitle


\begin{sciabstract} 
In an increasingly interconnected world, understanding and summarizing the structure of these networks becomes increasingly relevant. 
However, this task is nontrivial; proposed summary statistics are as diverse as the networks they describe, and a standardized hierarchy has not yet been established.  
In contrast, \mbox{vector-valued} random variables admit such a description in terms of their cumulants (e.g., mean, (co)variance, skew, kurtosis).  
Here, we introduce the natural analogue of cumulants for networks, building a hierarchical description based on correlations between an increasing number of connections, seamlessly incorporating additional information, such as directed edges, node attributes, and edge weights.  
These graph cumulants provide a principled and unifying framework for quantifying the propensity of a network to display any substructure of interest (such as cliques to measure clustering).\todoprivate{This unifying framework provides principled summary statistics for...} 
Moreover, they give rise to a natural hierarchical family of maximum entropy models for networks (i.e., ERGMs) that do not suffer from the ``degeneracy problem'', a common practical pitfall of other ERGMs.\todoprivate{<3_G wants to say more.}  
\todoscience{still ask for smaller abstract... maybe remove mention of degeneracy and just say natural family of null models or something. or notion is general part. or maybe just say node attributes or node and edge attributes. remove the building a hierarchical description of a network based on local correlations between an increasing number of connections. Moreover, it gives rise to a natural hierarchical family of maximum entropy models for networks. This family enjoys all the theoretical advantages associated with exponential random graph models (ERGMs), without suffering from their common practical pitfall known as the ``degeneracy problem''.} 
\end{sciabstract}

\todoscience{when submitting to science needs a one sentence summary. ``One Sentence Summary:  A brief teaser statement highlighting main result of the paper, understandable by a scientist not in your field, without jargon or abbreviations. This will appear online adjacent to the title and should not repeat phrases already present there. Please keep to under 125 characters.''} 

\todoscience{For science check if size of caption is good (should be at most 200 words).}


\todoscience{Maybe put section introduction when submitting to science}

The power of natural science relies on the ability to find abstract representations of complex systems and describe them with meaningful summary statistics.  
For example, using the pluripotent language of graph theory, the field of network science distills a variety of systems into the entities comprising them (the nodes) and their pairwise connections (the edges). 
\todoprivate{Remarkably many systems, which at the surface appear to be completely unrelated, naturally admit such a description: 
physical (e.g., electrical circuits \cite{klein1993resistance}, the cosmic web \cite{coutinho2016network}, Feymann diagrams \cite{feynman1949space}), 
biological (e.g., brains \cite{lynn2019physics}, food webs \cite{landi2018complexity}, protein interactions \cite{athanasios2017protein}), 
social (e.g., friendships \cite{wasserman1994social}, affiliations \cite{wang2009exponential}, collaborations \cite{bravohermsdorff2019gender}), 
and technological (e.g., transportation networks \cite{bell1997transportation}, financial transactions \cite{elliott2014financial}, the internet \cite{hall2012web}).} 
Remarkably many systems, from a wide variety of fields, naturally benefit from such a description, such as electrical circuits \cite{klein1993resistance}, brains \cite{lynn2019physics}, food webs \cite{landi2018complexity}, friendships \cite{wasserman1994social}, transportation networks \cite{bell1997transportation}, and the internet \cite{hall2012web}. 

As the field developed, real networks were noticed to display several recurring themes \cite{cimini2019statistical,barabasi2016network,newman2018networks}, such as: 
sparsity (a small fraction of all potential connections actually exist), 
heterogeneous degree distributions (some nodes have many more connections than average), 
and high clustering (nodes form tightly connected groups). 
As a consequence, measures are often tailored to capture such properties \cite{emmons2016analysis,goswami2018sparsity,demaine2019structural} and network models are often chosen to mimic them \cite{watts1998collective,barabasi1999emergence,colman2013complex,cai2016edge,kartun2018sparse,holme2002growing}. 
However, network statistics are often intertwined \cite{del2011all,orsini2015quantifying}, and models with notably different mechanisms \cite{cai2016edge,kartun2018sparse,holme2002growing,caron2017sparse,sampaio2015mandala,zuev2015emergence} can reproduce the same commonly observed properties. 
It is not currently clear how to compare different methods within a single principled framework; a standardized hierarchy of descriptive network statistics and associated models is needed.  

Real networks often contain a hierarchy of interconnected scales \cite{lovasz2012large}.  
This paper formalizes a novel \mbox{bottom-up} approach, where the elementary units are the edges (as compared to the nodes \cite{orsini2015quantifying}), and the hierarchy is built by considering correlations between an increasing number of them.  
Relationships between several edges can be expressed as subgraphs.  
Such substructures are often referred to as network motifs (or \mbox{anti-motifs}) when their appearance (or absence) in a network is deemed to be statistically significant \cite{milo2002network}.  
For example, triangular substructures appear in a wide range of real networks, indicating a common tendency for clustering.\todoprivate{cite}       
However, networks in different domains appear to display different motif profiles, suggesting that quantifying propensity for substructures could play a crucial role in understanding network function.\todoprivate{cite}   
For example, the \mbox{feed-forward} ($\ThreeTriangleLinearBig$) and the \mbox{bi-fan} ($\FourBifanBig$)\todoprivate{SHIFT DOWNER} substructures are prevalent in protein interaction networks\todoprivate{cite}, while the \mbox{bi-directional} \mbox{two-hop} ($\FourReciprocalReciprocalDirBig$) substructure is more prominent in transportation networks.
Current methods for quantifying a network's propensity for various substructures are generally characterized by the following two choices. 
First, one needs a network statistic to measure the propensity for the substructures of interest.  
%
Second, to determine statistical significance, one needs a null model for the observed network: an ensemble of networks that matches some set of properties deemed to be important.

For the first choice, a common approach is to simply count the number of instances of a substructure \cite{milo2002network,alon2007network,xia2019survey}.\todoprivate{As these motifs are overlapping structures made of the same material (the edges), the counts of these substructures are invariably inter-connected.} 
However, substructure counts can be misleading, and are generally not faithful measures of their propensities \cite{fretter2012subgraph,ginoza2010network}.     
For example, the number of triangles (or of any other substructure) naturally correlates with the number of larger substructures that contain them as subgraphs.  
Unfortunately, there is currently no standard measure of the propensity for an arbitrary substructure \cite{xia2019survey}; 
other choices tend to be rather \mbox{context-dependent} or specifically tailored to the substructure of interest (e.g., clustering coefficients), often leading to a proliferation of essentially incomparable statistics.\todoprivate{cite for clustering coefficients}   
This issue is exacerbated when one considers incorporating additional information, such as directed edges, node attributes, and edge weights.\todoprivate{cite}

Regardless of the measure used, in order to assess whether its value in the observed network is statistically significant, one must compare against its distribution in some appropriate null model. 
This distribution is often obtained by measuring its value in networks explicitly sampled from the chosen null model.\todoprivate{cite by computing this measure in networks sampled from the chosen null model.}    
Replicating the observed degree sequence (i.e., the configuration model) is a popular choice of null model.\todoprivate{cite}  
%
In some cases, this \mbox{node-centric} model may be appropriate, such as when the nodes in the observed network have unique identities that should indeed be preserved by the null model (e.g., generating randomized networks of transactions between some specific set of countries).\todoprivate{maybe make a general explanation, like network in another set of nodes doesnt make sense. maybe footnote about what we mean by unique identity}   
However, 
the defining symmetry of graphs is that nodes are \textit{a priori} indistinguishable (i.e., node exchangeability\todoprivate{cite}), and often one desires a null model that is similar to the observed network, but not necessarily with exactly the ``same'' nodes 
(e.g., generating typical networks of interactions between students to model the spread of infections).  
In addition, the configuration model does not treat all substructures equally. 
For example, hubs and clustering, both hallmark features of real networks, are associated with a propensity for \mbox{($k$-)stars} (i.e., a central node connected to $k$ others) and for \mbox{($k$-)cliques} (i.e., fully connected groups of $k$ nodes), respectively.\todoprivate{citations for hubs and cliques} 
However, the configuration model fully constrains the counts of all \mbox{$k$-stars}, and is therefore useless to assess their relative propensity in the observed network.   
Moreover, by design, it does not promote clustering at any level\todoprivate{cite}, and is thus inappropriate to study the \mbox{higher-order} organization of clustering.  
Even if one were to employ modifications to promote triangles\todoprivate{cite}, the issue is simply postponed to slightly larger substructures \cite{ritchie2017generation}. 

For a \mbox{real-valued} random variable, cumulants provide a hierarchical sequence of summary statistics that efficiently encode its distribution \cite{thiele1903theory,rota2000combinatorics,mccullagh2018tensor}.
Aside from their unique mathematical properties\footnote{In particular, their additive nature when applied to sums of independent random variables (see \SIreftext~\ref{SI:graphcumulantsadd}).}, \mbox{low-order} cumulants have intuitive interpretations. 
The first two orders, the mean and variance, correspond to the center of mass and spread of a distribution, and are taught in nearly every introductory statistics course \cite{gunderson2005interactive}. 
The next two orders have likewise been given unique names (skew and kurtosis), and are useful to describe data that deviate from normality, appearing in a variety of applications, such as finance \cite{wu2019forecasting}, economics \cite{vahamaa2003skewness}, and psychology \cite{blanca2013skewness}.
Generalizations of these notions have proven similarly useful: for example, joint cumulants (e.g., covariance) are the natural choice for measuring correlations in multivariate random variables. 
Unsurprisingly, cumulants are essentially universally used by the statistics community.  
\todoprivate{Additionally, deep connections with cumulants have appeared in a wide variety of contexts, such as statistical and particle physics \cite{kardar2007statistical,kubo1962generalized}.} 

By generalizing the combinatorial definition of cumulants \cite{kardar2007statistical,speed1983cumulants,speed1986cumulants}, we obtain their analogue for networks, which we refer to as \textit{graph cumulants}.  
%
After introducing the framework, we demonstrate their usefulness as principled and intuitive measures of the propensity (or aversiveness) of a network to display any substructure of interest, allowing for systematic and meaningful comparisons between networks. 
We then describe how this framework seamlessly incorporates additional network features, providing examples using real datasets containing directed edges, node attributes, and edge weights. 
Finally, we show how graph cumulants give rise to a natural hierarchical family of maximum entropy null models from a single observed network. 
\todoscience{We also provide a publicly available code that computes these graph cumulants.}   


For ease of exposition, we first introduce our framework for the case of simple graphs, i.e., undirected, unweighted networks with no \mbox{self-loops} or multiple edges.

\section*{Graph moments}

The cumulants of a \mbox{real-valued} random variable are frequently given in terms of its moments, i.e., the expected value of its powers.
A real network can be thought of as a realization of some \mbox{graph-valued} random variable, and its graph cumulants are likewise given in terms of its graph moments. 

%
To motivate our definition of graph moments, consider the measurement that provides the smallest nonzero amount of information about a network $\graphG$ with $n$ nodes. 
This measurement is a binary query, which yields $1$ if an edge exists between a random pair of nodes in $\graphG$, and $0$ otherwise\footnote{One might consider measuring a random node instead.  
However, the presence of a node \textit{per se} does not give any new information.  
On the other hand, querying the degree of a node yields an integer, containing more information than the binary result from querying the existence of a single edge.}. 
At first order, we consider repeated observations of a \textit{single} such measurement.  
We define the \mbox{first-order} graph moment $\mu_{1\oneedge}^{ }$ (the ``mean'') as the expected value of this quantity: the counts of edges in $\graphG$ normalized by the counts of edges in the complete graph with $n$ nodes. 
Hence, the \mbox{first-order} graph moment of a network is simply its edge density.\todoscience{probability of a connection. give a feel for density in terms of probability that a given substructure exists, say for substructure prob. of existence in terms of edges.} 

For the \mbox{second-order} graph moments, we again consider a binary query, but now using \textit{two simultaneous} such measurements.
This query yields $1$ if $\graphG$ contains edges between both pairs of nodes, and $0$ otherwise. 
We now must distinguish between two new cases: when the two edges share a node, thereby forming a wedge ($\mu_{2\twowedge}^{ }$); and when they do not share any node ($\mu_{2\twoparallel}^{ }$).  
Each case is associated with a different \mbox{second-order} graph moment, which we analogously define as the counts of the associated subgraph in $\graphG$, normalized by the counts of this subgraph in the complete graph with $n$ nodes. 
Hence, the \mbox{second-order} graph moments of a network are the densities of substructures with two edges.\todoscience{I think maybe here (or in the cumulants section), we could give an interpretation of the second order moments like wedge more scale free.more comments about interpretations of different cumulants that people might be interested in.}  

\begin{figure}[H]
\begin{center}
\centerline{\includegraphics[width=0.7\columnwidth]{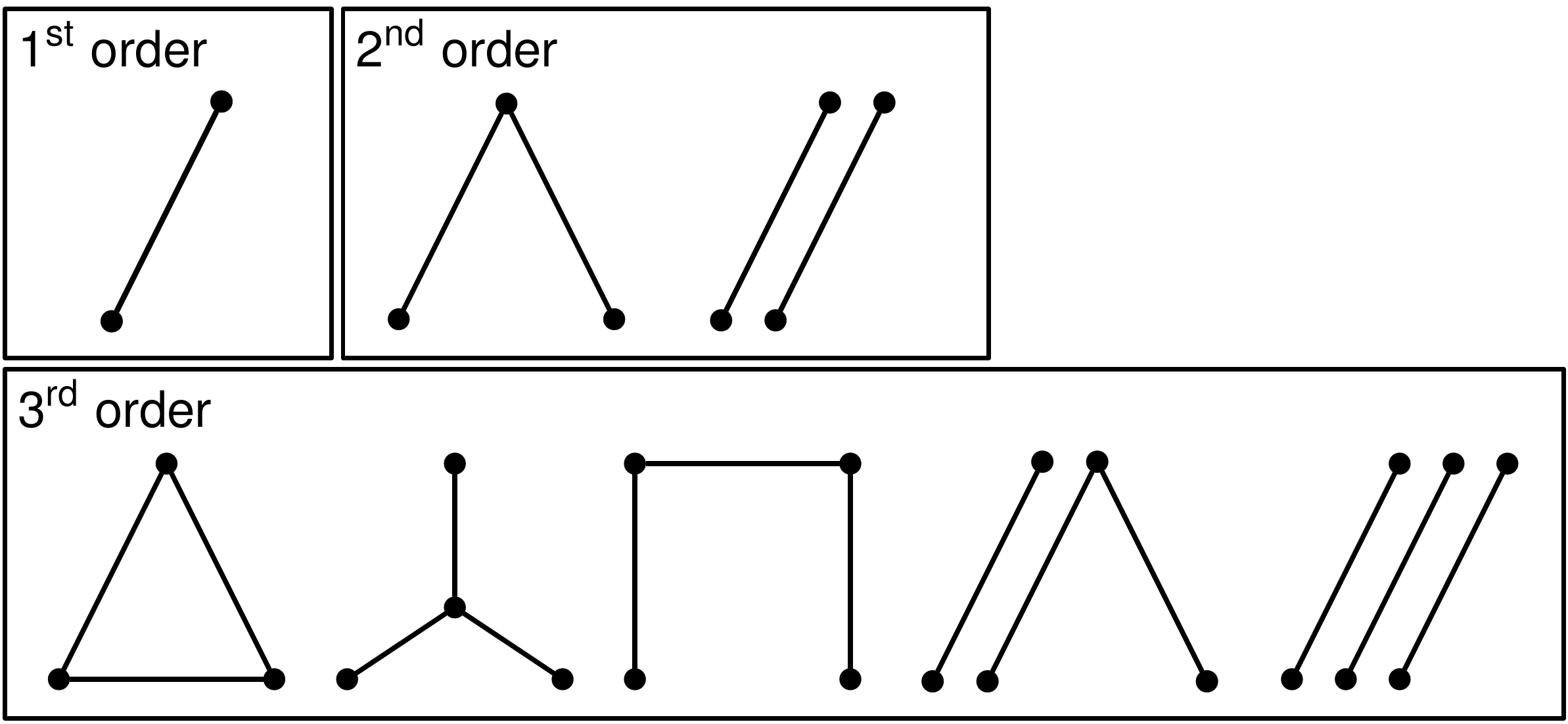}}
\caption{\textbf{Subgraphs associated with graph moments.}  
The graph moment $\mu_{rg}^{ }$ of a network $\graphG$ is defined as the counts of the subgraph $g$ (with $r$ edges), normalized by the counts of this subgraph if connections were present between all pairs of nodes in $\graphG$.\todoscience{maybe include a comment about density, and interpret it as a probability stuff for the edge. also use the name substructure.}    
Displayed here are the subgraphs $g$ associated with the graph moments up to third order (for simple graphs). 
The full set of \mbox{$r^{\text{th}}$-order} moments contains all substructures with exactly $r$ edges (including disconnected subgraphs), corresponding to all the ways that $r$ edges can relate to each other.
}
\label{fig:MomentsUndirectedUpto3rdorder} 
\end{center}
\end{figure}

Likewise, we define an \mbox{$r^\text{th}$-order} graph moment for each of the ways that $r$ edges can relate to each other (see \figreftext~\ref{fig:MomentsUndirectedUpto3rdorder} for the subgraphs associated with graph moments up to third order).\todoprivate{this is not my favorite...other options: Likewise, we define an \mbox{$r^\text{th}$-order} graph moment for each of the different configurations that $r$ distinct edges can form. Likewise, we define an \mbox{$r^\text{th}$-order} graph moment for each substructure that can be formed by $r$ edges.}     
Again, $\mu_{rg}^{ }$ is the density of subgraph $g$, defined as its counts in $\graphG$ normalized by its counts in the complete graph with the same number of nodes as $\graphG$.\todoscience{Again, $\mu_{rg}^{ }$ is the density of the substructure $g$ in $\graphG$. maybe remove from ``for the subgraphs associated with graph moments up to third order'' in parenthesis, or maybe phrase.  defined as the counts of the subgraph $g$ in $\graphG$, normalized by the counts of this subgraph in the complete graph with the same number of nodes as $\graphG$}  

The scalability of computing graph moments is determined by the computational complexity associated with counting connected subgraphs, as these imply the counts of disconnected subgraphs (see \SIreftext~\ref{SI:ConnectedToDisconnected} for details).  
This is an important fact, as the counts of disconnected subgraphs are generally orders of magnitude larger, 
and we can leverage numerous methods for efficiently counting connected subgraphs \cite{curticapean2017homomorphisms,fomin2012faster, amini2009counting,pinar2017escape}.

\section*{Graph cumulants}

\todoscience{I think here we could give an interpretation of the second order moments like wedge more scale free. more comments about interpretations of different cumulants that people might be interested in.} 
\todoscience{Figure: We can do networks and illustrate to show how the cumulants quantify certain structures. and say the properties scale-free, sparse, clustering. ideally, also think about a good one that people dont think about it.  Sketches (or actual samples?)} 

\todoprivate{Graph moments and graph cumulants are statistics of the entire network, corresponding to the density and propensity (or aversiveness) for a substructure, respectively.}

As the density of smaller substructures increases, the appearance of larger substructures that contain them as subgraphs will clearly also tend to increase. 
Hence, we would like to measure the difference between the observed value of a given graph moment and that which would be expected due to graph moments of lower order, so as to quantify the propensity (or aversiveness) for that specific network substructure. 
Cumulants are the natural statistics with this desired property. 
For example, the variance quantifies the intuitive notion of the ``spread'' of a distribution (regardless of its mean), while the skew and the kurtosis reflect, respectively, the asymmetry and contribution of large deviations to this spread.  

\begin{figure}[H]
\begin{center}
\centerline{\includegraphics[width=0.8\columnwidth]{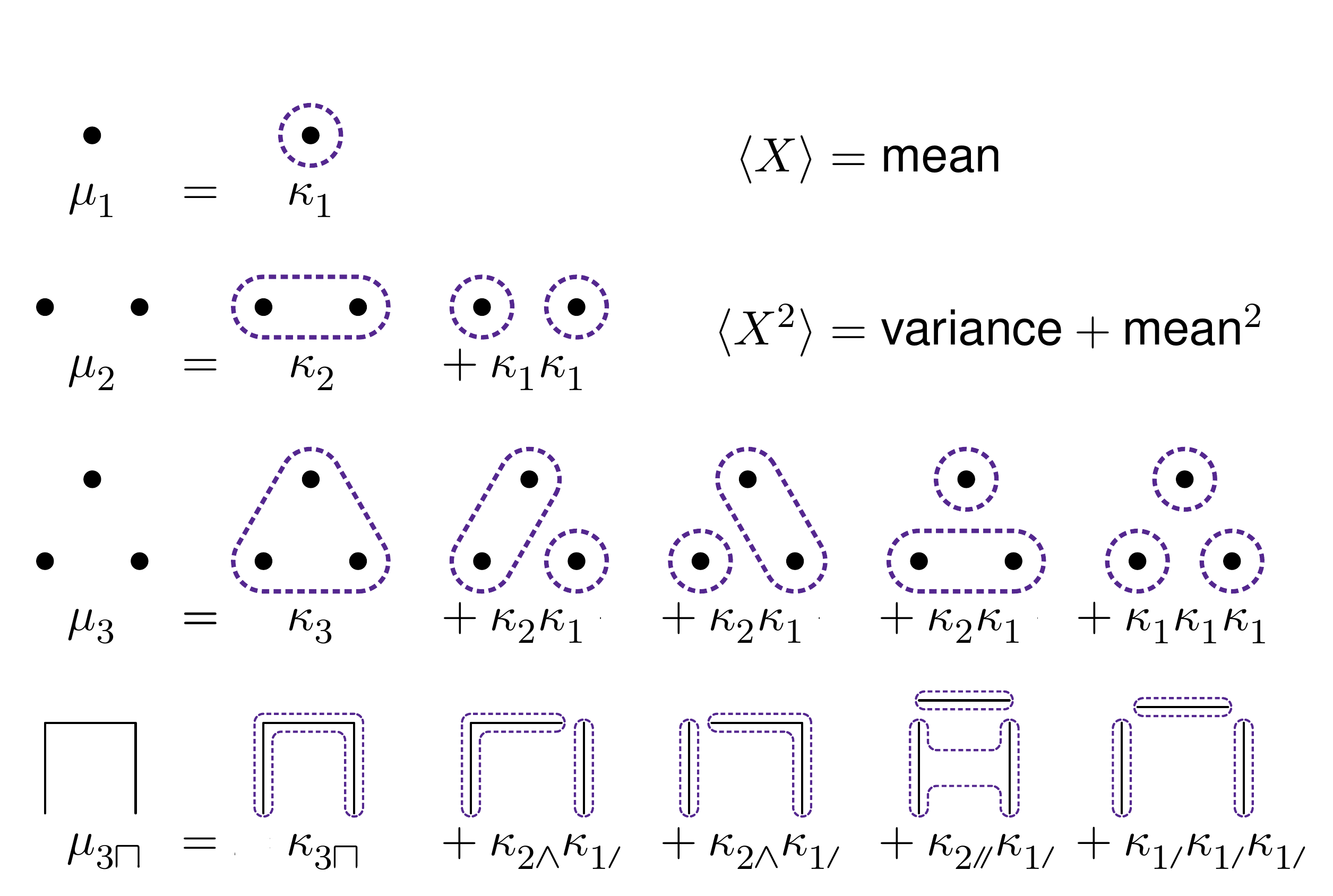}} 
\caption{
\textbf{To expand an \mbox{$r^{\text{th}}$-order} graph moment in terms of graph cumulants, enumerate all partitions of the $r$ edges comprising its subgraph.} 
The top three rows illustrate the combinatorial expansion of the classical moments in terms of cumulants.  
Analogously, the bottom row shows how to expand $\mu_{3\protect\threeline}^{ }$ in terms of graph cumulants.  
The last term ($\kappa_{1\protect\oneedge}^{3}$) corresponds to partitioning this subgraph into three subsets, each with a single edge. 
The first term ($\kappa_{3\protect\threeline}^{ }$) corresponds to ``partitioning'' this subgraph into a single subset containing all three edges, thus inheriting the connectivity of the entire subgraph.  
The remaining terms ($\kappa_{2}^{ }\kappa_{1}^{ }$) correspond to partitioning this subgraph into a subset with one edge and a subset with two edges.  
This can be done in three different ways: 
in two cases (the two $\kappa_{2\protect\twowedge}^{ }\kappa_{1\protect\oneedge}^{ }$ terms), the subset with two edges has those edges sharing a node; 
and in one case (the $\kappa_{2\protect\twoparallel}^{ }\kappa_{1\protect\oneedge}^{ }$ term), the subset with two edges has those edges not sharing any node.  
\todoscience{maybe make bottom row edges and line that group thicker.}   
}
\label{fig:Partition3Line}
\end{center}
\end{figure}

While often defined via the cumulant generating function \cite{hald2000early,gnedenko1949limit}, cumulants have an equivalent combinatorial definition \cite{kardar2007statistical,speed1983cumulants,speed1986cumulants} (see \figreftext~\ref{fig:Partition3Line}).  
At order $r$, it involves the partitions of a set of $r$ elements: 
\begin{equation}
\mu_r^{ } = \sum_{\pi \in P_r^{ }} \prod_{p \in \pi} \kappa_{|p|}^{ }, \label{eq:Partition}
\end{equation}
where $\mu_r^{ }$ is the $r^{\text{th}}$ moment, $\kappa_r^{ }$ is the $r^{\text{th}}$ cumulant, $P_r^{ }$ is the set of all partitions of a set of $r$ unique elements, $\pi$ is one such partition, $p$ is a subset of a partition $\pi$, and $|p|$ is the number of elements in a subset $p$.  

When generalizing this definition to graph moments, the partitioning $P_r^{ }$ of the edges must respect their connectivity (see \figreftext~\ref{fig:Partition3Line}, bottom row), i.e.
\begin{align}
\mu_{rg}^{ } &= \sum_{\pi \in P_{E(g)}^{ }} \prod_{p \in \pi} \kappa_{|p|\,g_p^{ }}^{ }, \label{eq:PartitionSI}
\end{align}
where $E(g)$ is the set of the $r$ edges forming subgraph $g$, $P_{E(g)}^{ }$ is the set of partitions of these edges, and $g_p^{ }$ is the subgraph formed by the edges in a subset $p$. 
These expressions can then be inverted to yield the graph cumulants in terms of graph moments (summarized in \SIreftext~\ref{SI:networkcumulantslist} and provided in our code).  

Essentially, the defining feature of cumulants is their additive nature when summing independent random variables \cite{thiele1903theory,rota2000combinatorics,gnedenko1949limit,hald2000early}.  
In \SIreftext~\ref{SI:graphcumulantsadd}, we show that the graph cumulants of independent \mbox{graph-valued} random variables also have this additive property for a natural notion of summing graphs.  
Moreover, we remark that the \mbox{Erd\H{o}s--R\'enyi} distribution has graph cumulants of zero for all orders \mbox{$r\geq2$}, similar to the classical cumulants of the normal distribution, which are zero for all orders \mbox{$r\geq3$}.\todoscience{check for n=4, then imagine more. maybe derive some realizability conditions. G wants to say stuff about ER having ``no graphical structure''......}  

\section*{Quantifying propensity for substructures:\\scaled graph cumulants}

Cumulants are often scaled, as dimensionless quantities allow for interpretable comparisons.\todoprivate{cite nigel renormalization book.gecia quote. dimensionless quantities are the lingua franca of science.}   
For example, the precision of a measurement is often quantified by the relative standard deviation ($\kappa_2^{\sfrac{1\!}{2}}/\kappa_1^{ }$), and the linear correlation between two random variables $X$ and $Y$ is often quantified by the Pearson correlation coefficient ($\text{Cov}(X,Y)/\kappa_2^{\sfrac{1\!}{2}}(X)\kappa_2^{\sfrac{1\!}{2}}(Y)$, i.e., their \mbox{second-order} joint cumulant divided by the geometric mean of their individual \mbox{second-order} cumulants).  
Likewise, we define scaled graph cumulants as \mbox{$\tilde{\kappa}_{rg}^{ } \equiv \kappa_{rg}^{ }/\kappa_{1\oneedge}^{r}$}, and report the ``signed'' $r^{\text{th}}$ root, i.e., the real number with magnitude equal to $\left|\tilde{\kappa}_{rg}^{ }\right|^{\sfrac{1\!}{r}}$ and with the same sign as $\tilde{\kappa}_{rg}^{ }$.\todoprivate{maybe explain our choice principled.} 
These scaled graph cumulants offer principled measures of the propensity (or aversiveness) of a network to display any substructure of interest. 
%

\begin{figure}[H]
\begin{center}
\centerline{\includegraphics[width=1\columnwidth]{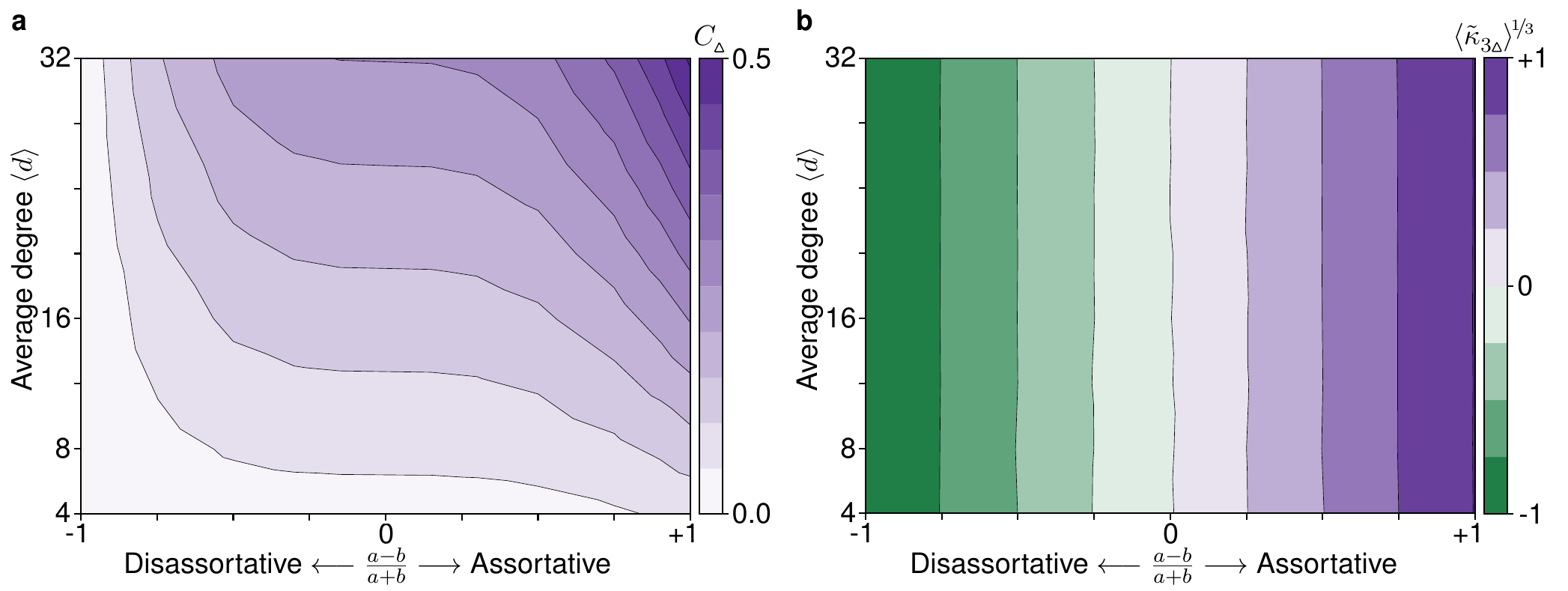}}
\centerline{\includegraphics[width=1\columnwidth]{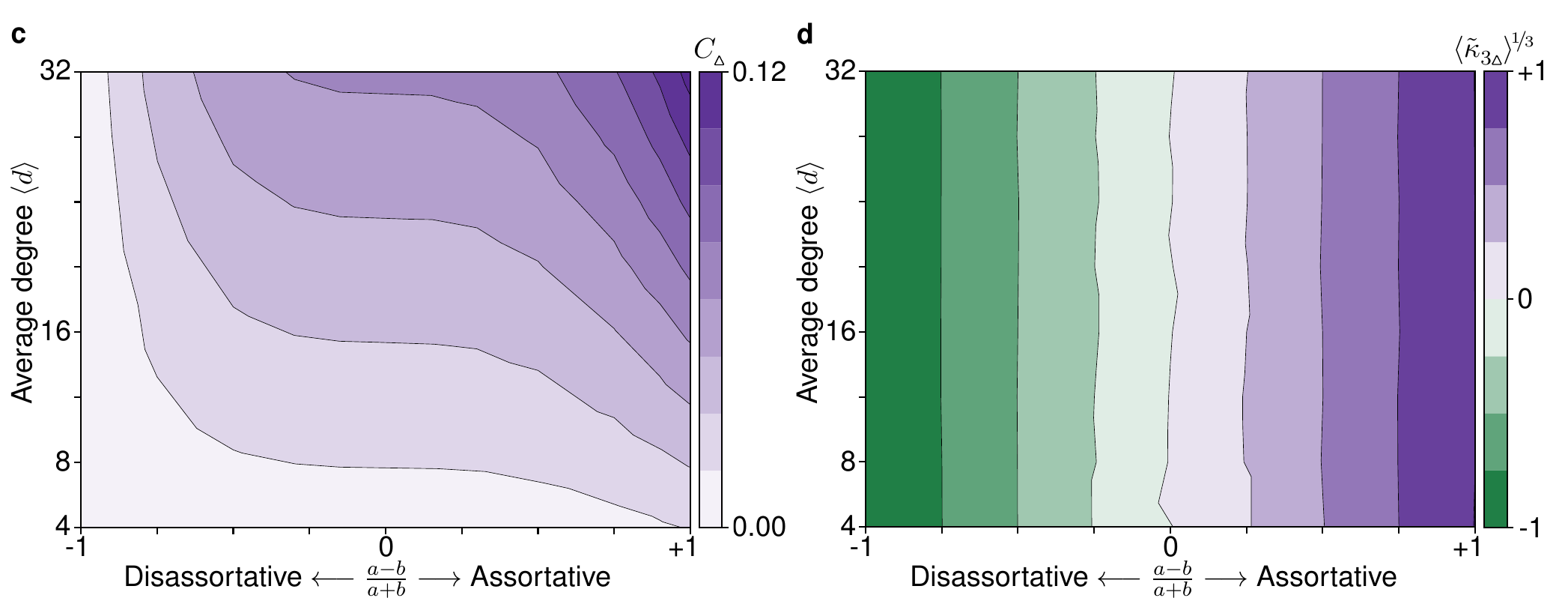}}
\caption{\textbf{The scaled triangle cumulant provides a principled measure of triadic closure.}   
We sampled networks from the symmetric stochastic block model on $n$ nodes with $2$ communities, SSBM$(n,2,a,b)$, for a range of $a$ and $b$, for both $128$ nodes (top) and $512$ nodes (bottom).  
The horizontal axis indicates the level of assortativity: $0$ corresponds to \mbox{Erd\H{o}s--R\'enyi} (ER) graphs (no community structure), \mbox{$+1$} to two disjoint ER graphs (only connections within the two communities), and \mbox{$-1$} to random bipartite graphs (only connections between the two communities, thus, no triangles).  
The vertical axis indicates the edge density, in terms of average degree.  
For each set of parameters, we compute the average of both the global clustering coefficient $C_{\protect\threetriangle}^{ }$ (left) and scaled triangle cumulant $\tilde{\kappa}_{3\protect\threetriangle}^{ }$ (right), displaying the signed fourth root for $\tilde{\kappa}_{3\protect\threetriangle}^{ }$.  
While both measures increase with assortativity (as desired), the clustering coefficient also varies with network size and edge density, whereas the scaled triangle cumulant is invariant to such differences.  
\todoprivate{While both measures increase with assortativity (as desired), the scaled triangle cumulant is invariant to network size and edge density, in contrast with...fuck.  }
\todoscience{Have the same picture for comparing with the local clustering coefficient, add comment about how results are the same here.} 
}
\label{fig:ClusteringTriangle}
\end{center}
\end{figure}

Scaling graph cumulants also allows for meaningful comparisons of the propensity of different networks to exhibit a particular substructure, even when these networks have different sizes and edge densities.  
To illustrate this point, we consider clustering, a hallmark feature of many real networks \cite{cimini2019statistical,barabasi2016network,watts1998collective}.  
This notion is frequently understood as the prevalence of triadic closure \cite{wasserman1994social,newman2018networks}, an increased likelihood that two nodes are connected if they have mutual neighbors, i.e., a propensity for triangles. 
This property is often quantified by the clustering coefficient $C_{\threetriangle}^{ }$, defined as the probability that two neighbors of the same node are themselves connected.\todoprivate{cite.}   
While this quantity is easily expressed within our formalism as \mbox{$C_{\threetriangle}^{ } = \mu_{3\threetriangle}^{ }/\mu_{2\twowedge}^{ }$}, it is neither a graph cumulant nor dimensionless.  
We propose that the scaled triangle cumulant $\tilde{\kappa}_{3\threetriangle}^{ }$ 
is a more appropriate measure of clustering in networks, as demonstrated in \figreftext~\ref{fig:ClusteringTriangle} (for the natural extension to clustering in bipartite networks, see \figreftextSI~\ref{fig:ClusteringInBipartiteNetworks}).

\todoscience{Check how similar to ours it is scaling the clustering coefficient and rooting it. if it is, we can use this in our favor and say look the clustering measures that work best, happened to be our principled derived scaled graph cumulants.}
\todoprivate{maybe put sexy statement about standardized measure here.}

\todoprivate{About the picture: In addition, the scaled triangle cumulant has a natural interpretation: a value of $0$ indicates that the number of triangles is precisely what is expected from the \mbox{lower-order} graph cumulants (manifestly true for the ER model).  
In addition, the scaled triangle cumulant has a natural interpretation: a value of $0$ indicates that the number of triangles is precisely what is expected from pure chance, i.e., an ER model.  
In fact, all scaled graph cumulants have this property. -- things we had before but are difficult to claim. could mention analogy between cumulants of normal distribution being zero and those of ER also being zero.}
\todoprivate{``, and a value of $-1$ indicates that there are no triangles in the network (as is the case in bipartite networks)'' 
Check for non-bipartite triangle-free networks. check how generic the -1 is for triangle free.  Maybe take cumulants of the distribution instead?}
\todoscience{Have the same picture for comparing with the local clustering coefficient, add comment about how results are the same here.} 

\section*{Graph cumulants for networks with additional features} 

While it is possible to treat most networks as undirected, unweighted graphs, real networks frequently contain more information (e.g., node attributes).  
Our formalism naturally incorporates many such augmentations. 
As before, there is a graph moment of order $r$ for each of the unique substructures (now endowed with the additional features) containing $r$ edges. 
The conversion to graph cumulants likewise respects the additional features.  
We now discuss the specific cases of directed edges, node attributes, and edge weights. 
In \SIreftext~\ref{SI:applicationstonetworkwithadditionalproperties}, we illustrate their ability to quantify the propensity for different substructures in a variety of real networks with these additional features. 
In \SIreftext~\ref{SI:networkcumulantslist}, we provide expressions for computing some of these augmented graph moments and graph cumulants.
For clarity, we consider each of these features individually; combining them is relatively straightforward.\todoscience{put tickz graphs downer.}

\subsection*{Directed edges}

When analyzing a directed network, the graph moments must incorporate the orientation of the edges.  
While the \mbox{first-order} graph moment ($\mu_{1\oneedgedir}^{ }$) still simply considers the number of directed edges 
(as edge orientation only carries meaning when considered with respect to the orientation of some other edge), 
there are now five \mbox{second-order} directed graph moments.  
The wedge configuration (two edges sharing one node) is now associated with three moments: one with both edges oriented towards the central node ($\mu_{2\twowedgeii}^{ }$), one with both edges oriented away from it ($\mu_{2\twowedgeoo}^{ }$), and one with an edge towards and the other away ($\mu_{2\twowedgeio}^{ }$).   
The relative orientation of two edges that do not share any node cannot be determined, and therefore is still associated with a single moment ($\mu_{2\twoparalleldir}^{ }$).  
Finally, the configuration of two reciprocal edges (i.e., two nodes that are connected by edges in both directions) is associated with the fifth \mbox{second-order} moment ($\mu_{2\tworeciprocaldir}^{ }$).  
The appropriate normalization is with respect to the counts in the complete directed graph, i.e., that which has every pair of nodes connected by edges in both directions.  
Fig.~\ref{fig:DirectedStructure_Yeast} illustrates how incorporating the directed nature of protein interaction networks reveals additional structure.  


\subsection*{Node attributes} 

Often, nodes of a network have intrinsic attributes (e.g., demographics for social networks).  
The graph moments of such networks are defined by endowing the subgraphs with same attributes. 
For example, consider the case of a network in which every node has one of two possible \mbox{``flavors''}: ``charm'' and ``strange''.\todoscience{cite particle physics for joke. maybe Democrat/republican. belong of one of two categories. or groups. or colors. or affiliation....colors might be the best. make instead of node color put stripes or dotted.}   
There are now three \mbox{first-order} graph moments: an edge between two ``charm'' nodes ($\mu_{1\OnePurpleToPurpleEdge}^{ }$), an edge between two ``strange'' nodes ($\mu_{1\OneGreenToGreenEdge}^{ }$), and an edge between one of each ($\mu_{1\OnePurpleToGreenEdge}^{ }$).  
To compute the moments, we normalize by their counts in the complete graph on the same set of nodes: here, ${n_{\PurpleDot}^{ }}\choose{2}$, ${n_{\GreenDot}^{ }}\choose{2}$, and $n_{\PurpleDot}^{ }n_{\GreenDot}^{ }$, respectively. 
Fig.~\ref{fig:NodeLabel_GenderedPrimaryClass} considers a (binary) gendered network of primary school students \cite{stehle2011high}, illustrating how incorporating node attributes elucidates the correlations between node type and their connectivity patterns. 

A common special case of networks with node attributes are bipartite networks: nodes have one of two types (e.g., authors and publications \cite{bravohermsdorff2019gender}, plants and pollinators \cite{campbell2011network}) and edges are only allowed between nodes of different type.  
As certain subgraphs are now unrealizable, bipartite networks have only one \mbox{first-order} moment (an edge connecting a ``charm'' to a ``strange'', $\mu_{1\OnePurpleToGreenEdge}^{ }$), and two \mbox{second-order} wedge moments: a ``charm'' node connected to two ``strange'' nodes ($\mu_{2\TwoGreenPurpleGreen}^{ }$), and a ``strange'' node connected to two ``charm'' nodes ($\mu_{2\TwoPurpleGreenPurple}^{ }$).

\subsection*{Edge weights}
To compute graph moments for weighted networks, subgraphs should be counted with multiplicity equal to the product of their edge weights \cite{lovasz2012large} (see \SIreftext~\ref{SI:graphcumulantsadd} for a detailed justification).  
The normalization is the same as in the unweighted case, i.e., divide the (weighted) count of the relevant subgraph by the count of this subgraph in the unweighted complete graph with the same number of nodes. 
Note that, unlike unweighted networks, graph moments may be greater than one for weighted networks. 
In \figreftextSI~\ref{fig:WeightedComparison_Infectious}, we analyze a weighted network of social interactions \cite{isella2011whats}, illustrating how allowing for variable connection strength can increase statistical significance and even change the resulting interpretations.\todoscience{like, actually check this shit}  

\todoprivate{***Bef applying to science. Add a section where we do local node and local edge and show application to ml. link prediction and node embedding. maybe a graph neural net.}

\section*{Inference from a single network} 

\subsection*{Unbiased estimators of graph cumulants}
\todoscience{Here/in this section, we need to motivate why one would want to do infer a probability distribution.  
in terms of testing if properties are stat significant. dk-series style.}

Thus far, we have not made a distinction between the graph moments of an observed network $\graphG$ and those of the distribution $\graphGrv$ from which this network was sampled.  
This is because, in a sense, they are the same: \mbox{$\langle \mu_{rg}^{ }(\graphG)\rangle = \mu_{rg}^{ }(\graphGrv)$} (where the angled brackets \mbox{$\langle\,\cdot\,\rangle$} denote expectation with respect to the distribution $\graphGrv$), a property known as ``inherited on the average'' \cite{tukey1950some}. 

However, for cumulants, this distinction is important.  
Cumulants of a distribution are defined by first computing the moments of the distribution, then converting them to cumulants (as opposed to computing the cumulants of the individual samples, then taking the expectation of those quantities).  
Due to the fact that the cumulants are nonlinear functions of the moments, they are not necessarily preserved in expectation, i.e., not inherited on the average \cite{tukey1950some,fisher1930moments}.  
For example, consider the problem of estimating the variance of a distribution over the real numbers given $n$ observations from it.  
Simply using the variance of these observations gives an estimate whose expectation is less than the variance of the underlying distribution, and one should multiply it by the \mbox{well-known} correction factor of $\frac{n}{n-1}$.  
The generalizations of this correction factor for \mbox{higher-order} cumulants are known as the \mbox{$k$-statistics} \cite{fisher1930moments}; 
given a finite number of observations, they are the \mbox{minimum-variance} unbiased estimators of the cumulants of the underlying distribution \cite{halmos1946theory,kendall1940some,kendall1940proof,kenney1951mathematics,james1958moments,fisher1982unbiased}. 

In many applications, one wishes to estimate a probability distribution over networks $\graphGrv$ after observing only a single network $\graphG$.  
Just as with classical cumulants, applying equation~\ref{eq:PartitionSI} to this network also yields biased estimates of the graph cumulants of the underlying distribution, i.e., \mbox{$\langle \kappa_{rg}^{ }(\graphG)\rangle \neq \kappa_{rg}^{ }(\graphGrv)$}. 
In \SIreftext~\ref{SI:unbiasednetworkcumulants}, we describe a procedure to obtain unbiased estimators of graph cumulants \mbox{$\check{\kappa}_{rg}^{ }$}, 
as well as the variance of these estimators \mbox{$\text{Var}(\check{\kappa}_{rg}^{ })$} (\SIreftext~\ref{SI:unbiasedcumulantsergmquantification_variance}).\todoprivate{***maybe we should say obtain these unbiased estimators (under a natural null-hypothesis or something like that), we could do footnote? we assume graphons.}  
In particular, for simple graphs, we provide the expressions for these unbiased estimators up to third order, as well as the variance for first order.  
Obtaining a complete list of the expressions for the unbiased estimators of graph cumulants and their variance would provide a powerful tool for network science (see \SIreftext~\ref{SI:unbiasedcumulantsergmquantification_test} for a detailed discussion). 
This would allow for principled statistical tests of the propensity (or aversiveness) for any substructure without explicitly sampling from some null model, a procedure that is in general quite computationally expensive and often a main obstacle in the analysis of real networks \cite{orsini2015quantifying,butts2017aperfect,robins2007anintroduction,wang2009exponential,veitch2019sampling}. 

However, sometimes one indeed requires samples from the underlying distribution, such as when assessing the statistical significance of some property that cannot be easily expressed in terms of the statistics defining this distribution.  
We now discuss how these unbiased graph cumulants can be used to obtain a principled hierarchical family of network models (see \SIreftext~\ref{SI:fittingergm} and~\ref{SI:geometricdegeneracy} for more details).

\subsection*{A principled hierarchical family of network models}  

A ubiquitous problem, arising in many forms, is that of estimating a probability distribution based on partial knowledge. 
Often, one desires the distribution to have some set of properties, but the problem is typically still highly unconstrained.\todoprivate{Often, one knows of a set of properties that the distribution should have, but the solution is still highly unconstrained. G thinks know is better cause max entropy is really about information, blabla.}  
In such cases, the maximum entropy principle \cite{jaynes1957information1,jaynes1957information2,amari2016information} provides a natural solution: of the distributions that satisfy these desired properties, choose the one that assumes the least amount of additional information, i.e., that which maximizes the entropy.  
For example, when modeling \mbox{real-valued} data, one often uses a Gaussian, the maximum entropy distribution (over the reals) with prescribed mean and variance.

The analogous maximum entropy distributions for networks are known as exponential random graph models (ERGMs). 
These models are used to analyze a wide variety of \mbox{real} networks \cite{saul2007exploring,chakraborty2019exponential,koskinen2018outliers,park2004statistical,buddenhagen2017epidemic}.\todoprivate{remove some of these references.} 
Although it is possible to define an ERGM by prescribing any set of realizable properties, attention is often given to the counts of edges and of other \mbox{context-dependent} substructures (such as wedges or triangles for social networks \cite{wasserman1994social,robins2007anintroduction,strauss1986onageneral,park2004solution}).\todoprivate{remove some of these references.}    
While individual networks sampled from this distribution do not necessarily have the same counts (of edges and of specified substructures), the average values of these counts are required to match those of the observed network.

Unfortunately, ERGMs of this type often result in pathological distributions, exhibiting strong multimodality, such that typical samples have properties far from those of the observed network they were intended to model.  
For example, the networks sampled from such an ERGM might be either very dense or very sparse, despite the fact that averaging over these networks yields the same edge count as that of the observed network.  
%
%
This phenomenon is known as the ``degeneracy problem'', and much effort has gone into understanding it \cite{chatterjee2013estimating,park2004solution,park2005solution}. 
While some remedies have been proposed \cite{horvat2015reducing,snijders2006new,caimo2011bayesian}, a principled and systematic method for alleviating degeneracy has thus far remained elusive. 

\begin{figure}[H]
\begin{center}
\centerline{\includegraphics[width=1\columnwidth]{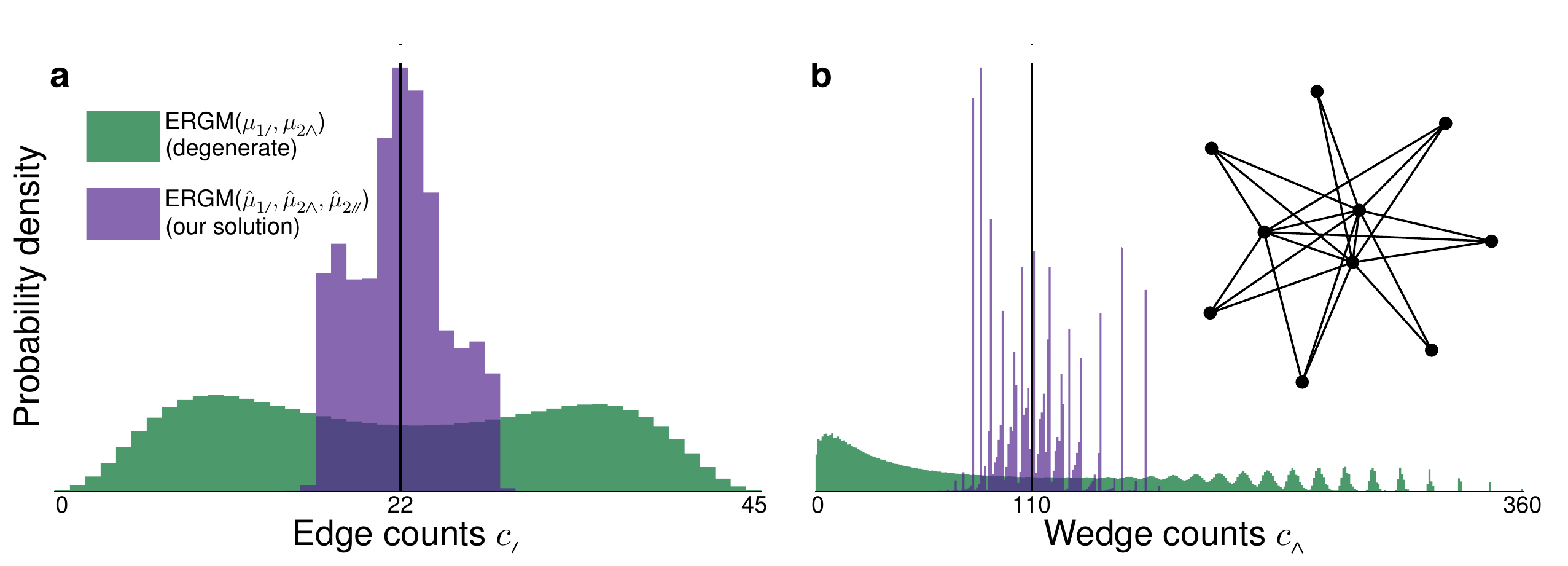}}
\caption{\textbf{Our proposed family of ERGMs does not suffer from the ``degeneracy problem''.}   
We compute the exact distributions over simple graphs with $10$ nodes resulting from fitting two ERGMs to the same observed network (shown on the right), and plot the resulting distributions of: \textbf{a)} edge counts, and \textbf{b)} wedge (i.e., \mbox{$2$-star}) counts.  
Green corresponds to a network model frequently used in the literature: $\text{ERGM}(\mu_{1\protect\oneedge}^{ },\mu_{2\protect\twowedge}^{ })$, the maximum entropy distribution with prescribed expected counts of edges and wedges.  
Purple corresponds to our analogous \mbox{second-order} network model: $\text{ERGM}(\check{\mu}_{1\protect\oneedge}^{ },\check{\mu}_{2\protect\twowedge}^{ },\check{\mu}_{2\protect\twoparallel}^{ })$, with prescribed expected \textit{unbiased} counts of substructures of first and second order, i.e., edges, wedges, \textit{and} two edges that do not share any node. 
We fit the model using the procedure described in \SIreftext~\ref{SI:fittingergm}, with unbiasing parameter \mbox{$\eta=\frac{1}{11}$}.\todoscience{we have only 200 words, create a section on the supp materials that has details of methods in pictures.}   
Black lines denote the counts in the observed network that these distributions are intended to model. 
In sharp contrast to our proposed $\text{ERGM}(\check{\mu}_{1\protect\oneedge}^{ },\check{\mu}_{2\protect\twowedge}^{ },\check{\mu}_{2\protect\twoparallel}^{ })$, the currently used $\text{ERGM}(\mu_{1\protect\oneedge}^{ },\mu_{2\protect\twowedge}^{ })$ can result in a distribution whose typical samples are notably different from the observed network.
This is reflected in the fact that the green distributions have maxima far from their means.  
This undesirable behavior, known as the degeneracy problem, tends to become even more pronounced for larger networks \cite{horvat2015reducing}.\todoprivate{on the supp. mat. detailing this picture, comment about the little two picks on the edge figure.}   
In \SIreftext~\ref{SI:geometricdegeneracy}, we explain why this occurs and why our proposed family of ERGMs does not suffer from this problem.\todoscience{Make callouts in the color of the distributions for the graphs with the picks, showing for example that our wedge picks make sense.} 
}
\label{fig:Degeneracy}
\end{center}

\end{figure}

Based on our framework, we propose a hierarchical 
family of ERGMs that are immune to the degeneracy problem (see \figreftext~\ref{fig:Degeneracy} for an example and \SIreftext~\ref{SI:geometricdegeneracy} for details).
Our hierarchy is based on correlations between an increasing number of individual connections: at $r^{\text{th}}$ order, the proposed ERGM is specified by \textit{all} the graph cumulants of order at most $r$.  
Importantly, when inferring such an ERGM $\graphGrv$ from a single observed network $\graphG$, it is appropriate to use the \textit{unbiased} graph cumulants, such that \mbox{$\kappa_{rg}^{ }(\graphGrv) = \check{\kappa}_{rg}^{ }(\graphG)$} for all subgraphs through the chosen order (see \SIreftext~\ref{SI:fittingergm} for the detailed protocol).

There are two main differences between our proposed family of ERGMs and those typically used in the literature.\todoprivate{cite} 
First, our family prescribes the expected counts of \textit{all} subgraphs (including disconnected) with at most the desired order (see \SIreftext~\ref{SI:spectralinterpretation} for a motivation of this choice based on a spectral representation of distributions over networks).  
In contrast, current ERGMs only consider some (usually connected) subgraphs of interest, perhaps due to the fact that not all of these subgraphs are deemed important to model the observed network, or because disconnected subgraphs are not usually thought of as motifs \cite{milo2002network,benson2016higher}. 
Second, the use of unbiased graph cumulants results in a distribution with expected subgraph counts different from those of the observed network. 
Nevertheless, the ERGM distribution induced by this choice generates samples that are appropriately clustered around the observed network (see \figreftext~\ref{fig:Degeneracy} and \SIreftext~\ref{SI:geometricdegeneracy} for detailed intuition). 

\subsection*{Discussion} 

Over a century ago, Thiele introduced cumulants \cite{thiele1889general}, a concept now fundamental to the field of statistics \cite{rota2000combinatorics,mccullagh2018tensor,gnedenko1949limit,hald2000early,speed1983cumulants,speed1986cumulants}, 
which has justifiably percolated throughout the scientific community \cite{wu2019forecasting,vahamaa2003skewness,blanca2013skewness,kardar2007statistical,kubo1962generalized}. 
In this work, we introduce \textit{graph cumulants}, their generalization to networks (\figreftext~\ref{fig:Partition3Line}).  
This principled hierarchy of network statistics provides a framework to systematically describe and compare networks (\figreftexts~\ref{fig:ClusteringTriangle} and~\ref{fig:ClusteringInBipartiteNetworks}), naturally including those with additional features, such as directed edges (\figreftextSI~\ref{fig:DirectedStructure_Yeast}), node attributes (\figreftextSI~\ref{fig:NodeLabel_GenderedPrimaryClass}), and edge weights (\figreftextSI~\ref{fig:WeightedComparison_Infectious}).  
Moreover, through the lens of the maximum entropy principle, these statistics induce a natural hierarchical family of network models.  
These models are immune to the ``degeneracy problem'', providing a principled prescription for obtaining distributions that are clustered around the properties of the network they intend to model (\figreftexts~\ref{fig:Degeneracy},~\ref{Fig:GeometryManifold}, and~\ref{Fig:Stereoscopic}).


%
To make appropriate predictions, one must acknowledge that the observed data are but one instantiation, inherently incomplete and stochastic, of some underlying process.\todoprivate{maybe remove pothead comment...}   
In network science, one typically has a single network observation, and would like to make inferences about the distribution from which it came.    
This is analogous to characterizing the distribution of a classical random variable given a finite collection of samples.  
However, aside from the mean, the cumulants of a finite sample are not the same (in expectation) as those of the underlying distribution.  
The desired unbiased estimators are known as the \mbox{$k$-statistics} \cite{fisher1930moments,kenney1951mathematics,fisher1982unbiased} (e.g., the $\frac{n}{n-1}$ correction factor for the sample variance).  
Characterizing the distributions of these unbiased estimators allows for a variety of principled statistical tests (e.g., the use of the $\chi^2_{ }$ distribution for analyzing the sample variance).\todoprivate{cite}   
In \SIreftext~\ref{SI:unbiasednetworkcumulants}, we provide a procedure for deriving the analogous unbiased estimators of graph cumulants given a single network observation, and in \SIreftext~\ref{SI:unbiasedcumulantsergmquantification_variance}, we provide a procedure for deriving their variance.   
While the derivations are incredibly tedious, once the expressions are obtained, they could be used to systematically measure the statistical significance of the propensity (or aversiveness) for arbitrary substructures (see \SIreftext~\ref{SI:unbiasedcumulantsergmquantification_test}).
This is an incredibly promising avenue, as it circumvents the need for constructing and sampling from a network null model, a major challenge for many current methods. 
As with the rest of our framework, such an analysis naturally incorporates additional features, such as directed edges, node attributes and edge weights.  

\todoprivate{**compare our work with graphons, dk-series, higher-order clustering, i also added in reference a folder with more examples of related stuff. }
\todoprivate{make more general, our other examples, do comparison of statistical significance with using counts}


Graph cumulants quantify the propensity for substructures throughout the \textit{entire} network.  
However, in some applications, statistics that quantify the propensity of \textit{individual} nodes (or edges) to participate in these substructures are more appropriate.  
%
For example, the local clustering coefficient \cite{watts1998collective} aims to describe the propensity of a node to participate in triangles (as does the edge clustering coefficient for edges).\todoprivate{cite}   
However, as before, there is no general framework for arbitrary substructures.  
Our graph cumulant formalism again offers a systematic prescription.  
Essentially, the local graph cumulants are obtained by giving the node (or edge) of interest a unique identity and computing the graph cumulants of the entire network with this augmented information (see \SIreftext~\ref{SI:localcumulants} for details).  
These local graph cumulants could serve as useful primitives in machine learning tasks such as node classification \cite{kipf2016semi,hamilton2017representation,hamilton2017inductive} and link prediction \cite{liben2007link}.

Just as the scientific community has converged upon the variance as the canonical measure of the spread of a distribution, 
the field of network science could greatly benefit from a similarly standardized measure of the propensity for an arbitrary substructure. 
Inspired by over a century of work in theoretical statistics, the framework of graph cumulants introduced in this paper provides a uniquely principled solution.  


\nocite{csardi2006igraph,mathematicaprogram,danisch2018listing,chiba1985arboricity,kuba2007degree,houbraken2014index,wernicke2006efficient,
 slota2013fast,ribeiro2010efficient,aliakbarpour2018sublinear, miyajima2014continuous,fagiolo2007clustering,barrat2004thearchitecture,
 antoniou2008statistical,robins2004small,brunson2015triadic,lind2005cycles,bhardwaj2010analysis,oeis,liben2007link,nouri2012graph,gao2010survey,
 rucinski1988small,maugis2020testing,bickel2011method}

\bibliography{EdgeCorrelationCitationForScience.bib}

\bibliographystyle{Science}

\todoscience{Include acknowledgements here when submitting to science ``Include acknowledgments of funding, any patents pending, where raw data for the paper are deposited, etc.'' 
}

\newpage

\section*{\Huge{Supplementary Materials}}
\label{SI:main}

\renewcommand*{\thesection}{S\the\value{section}}
\setcounter{section}{0}

\renewcommand{\figurename}{\textbf{Fig.}}
\renewcommand{\thefigure}{S\arabic{figure}}
\setcounter{figure}{0}

\section{Python module for computing graph cumulants}
\label{SI:PythonCode}

We provide a python module that computes graph moments of networks, including networks with directed edges, edge weights, and binary node attributes. 
Conversions to and from graph cumulants (as well as their unbiased counterparts) are also implemented.

We use the python package igraph \cite{csardi2006igraph} to count instances of subgraphs in a network. 
In conjunction with the symbolic computation available in Mathematica \cite{mathematicaprogram}, we automatically derived the expressions for the counts of disconnected subgraphs in terms of the connected counts (see \SIreftextInSI~\ref{SI:ConnectedToDisconnected}), as well as the expressions for converting graph moments to graph cumulants (see  \SIreftextInSI~\ref{SI:networkcumulantslist}).

Our code contains the expressions to obtain graph cumulants up to: 
sixth order for undirected unweighted networks, 
fifth order for undirected weighted networks (as well as all sixth order connected subgraphs), 
fifth order for directed unweighted networks (as well as the sixth order $K_3^{ }$ subgraph), 
subgraphs of $K_{2,2}^{ }$ for bipartite networks, 
and subgraphs of $K_3^{ }$ for networks with binary node attributes. 
Expressions for the unbiased graph cumulants are implemented up to third order for undirected unweighted networks (see \SIreftextInSI~\ref{SI:unbiasednetworkcumulants}).

\subsection{Efficiently computing graph moments}
\label{SI:ConnectedToDisconnected}

The scalability of our framework is determined by the computational time required to count the relevant connected subgraphs. 
This is because the counts of the disconnected subgraphs can be inferred by the counts of the connected subgraphs, and therefore do not need to be explicitly enumerated. 
This is an important fact, as the counts of the disconnected subgraphs are generally orders of magnitude larger than those of the connected ones.  
Moreover, counting connected subgraphs is an active area of research, and much work has gone into their efficient computation \cite{curticapean2017homomorphisms, fomin2012faster, amini2009counting,pinar2017escape}. 

To illustrate how the counts of disconnected subgraphs are derivable from the connected subgraph counts, consider the case of \mbox{second-order} moments for simple graphs. 
From first order, one has the counts of edges in the network, \mbox{$c_{\oneedge}^{ } = {n \choose 2} \mu_{1\oneedge}^{ }$}.  
Consider all unordered pairs of distinct edges; each pair corresponds to a single \mbox{second-order} count: either of a wedge, or of two edges that do not share any node, thus \mbox{${c_{\oneedge}^{ } \choose 2} = c_{\twowedge}^{ } + c_{\twoparallel}^{ }$}.  
Hence, the count of two edges that do not share any node \mbox{$c_{\twoparallel}^{ }$} is directly derivable from the count of edges \mbox{$c_{\oneedge}^{ }$} and the count of wedges \mbox{$c_{\twowedge}^{ }$}. 

A similar argument applies to all orders.  
For instance, at third order, there are two disconnected subgraphs: three edges that do not share any node, and a wedge and an edge that do not share any node.   
By enumerating all triplets of distinct edges, as well as all pairs of a wedge and an edge not contained in that wedge, we obtain the following expressions:  
\begin{align*}
{c_{\oneedge}^{ } \choose 3} &= c_{\threeparallel}^{ } + c_{\threetriangle}^{ } + c_{\threeclaw}^{ } + c_{\threeline}^{ } + c_{\threeedgewedge}^{ }, \\
c_{\twowedge}^{ } (c_{\oneedge}^{ } - 2) &= c_{\threeedgewedge}^{ } + 3c_{\threetriangle}^{ } + 3c_{\threeclaw}^{ } + 2c_{\threeline}^{ }. 
\end{align*}

We now discuss the scalability of counting the instances of a connected subgraph $g$ with $n'$ nodes in a network $\graphG$ with $m$ edges and $n$ nodes.  
The complexity of a na\"{i}ve enumeration of the $\frac{n!}{(n-n')!}$ potential node mappings scales as $\mathcal{O}(n^{n'})$.  
However, there exist notably more efficient algorithms for certain subgraphs (such as triangles, stars, and cliques), especially when $\graphG$ has particular properties, such as sparsity \cite{benson2016higher,danisch2018listing,houbraken2014index}. 
For example, the \mbox{worst-case} computational time complexity for counting \mbox{$n'$-cliques} is known to be at most $\mathcal{O}(n' m^{\sfrac{n'\!}{2}})$ time \cite{chiba1985arboricity}.  
The counts of the $r$-stars can be quickly computed, as they are proportional to the $r^{\text{th}}$ factorial moments of the degree distribution \cite{kuba2007degree,aliakbarpour2018sublinear}. 
Moreover, some of these algorithms can be substantially accelerated through parallel computation \cite{slota2013fast,ribeiro2010efficient} and approximate values can be obtained by stochastic methods \cite{wernicke2006efficient}.\todoprivate{maybe say using sampling techniques?}    
Asymptotics aside, from a pragmatic perspective, Pinar et al. \cite{pinar2017escape} showed that the exact counts of all connected subgraphs with up to 5 nodes can be obtained for networks with tens of millions of edges in minutes on a commodity machine (64GB memory). 

\todoprivate{Note that $c_g^{ }$ is the number of \textit{unique} instances of $g$ in $\graphG$; if one instead counts the number of subisomorphisms (i.e., mappings from the nodes of $g$ to a subset of nodes in $\graphG$), different node mappings may give the same subgraph, and each instance will be counted $\left|\text{Aut}(g)\right|$ times, where $\text{Aut}(g)$ is the automorphism group of the subgraph $g$. ... *** if include -- subisomorphisms -- we can look at some math papers to see what they mean in terms of embedding.} 

\section{Applications to networks with additional features}
\label{SI:applicationstonetworkwithadditionalproperties}

In this section, we first demonstrate how our framework provides a natural notion of clustering in bipartite networks.  
We then illustrate the utility of incorporating additional network features by analyzing real networks with directed edges, node attributes, and weighted edges.\todoprivate{mention that here is multiple substructures.}   

\subsection{Quantifying clustering} 
\label{SI:subclustering}

\begin{figure}[H]
\begin{center}
\centerline{\includegraphics[width=1\columnwidth]{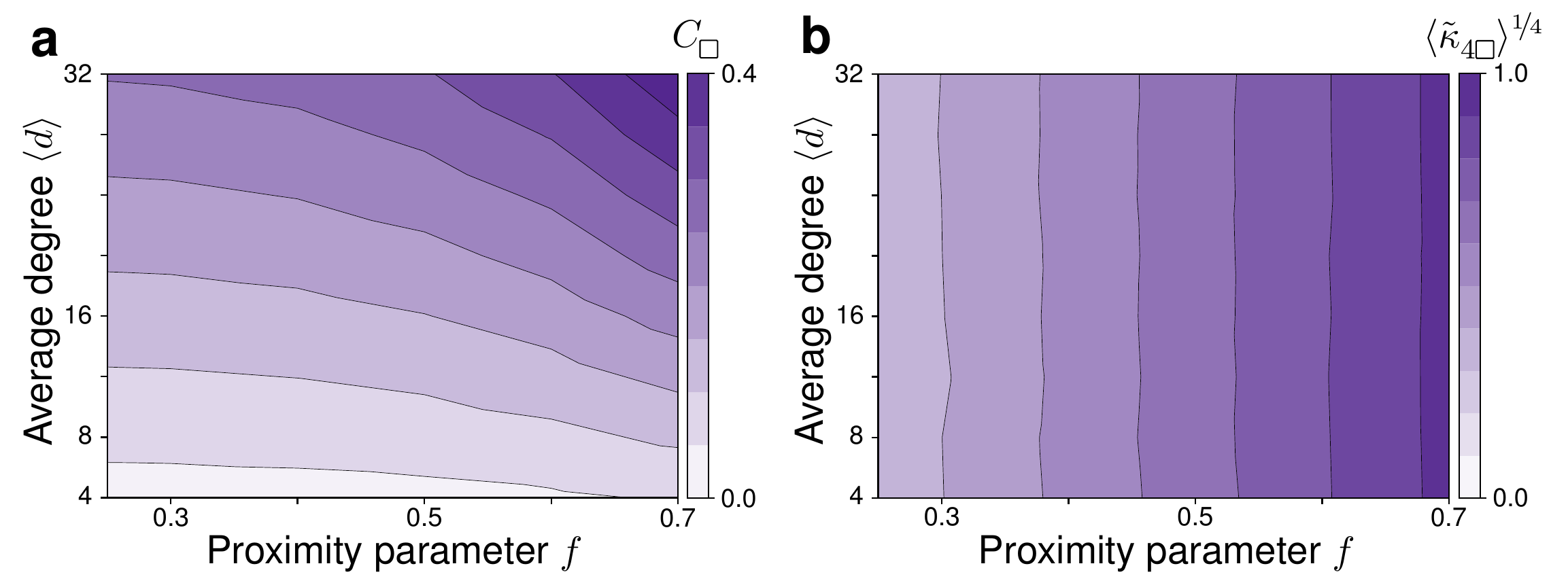}}
\caption{\textbf{A principled measure of clustering in bipartite networks.}  
We simulated a bipartite geometric graph model with two parameters, $f$ and $\langle d \rangle$, which determine the propensity for clustering (horizontal axis) and the edge density (vertical axis), respectively. 
The model first divides the nodes (here, $256$) into two groups of equal size and randomly places them on the unit sphere.  
Each node may only connect to nodes from the other group, and only when they are within a certain radius, such that the area it contains is a fraction \mbox{$1-f$} of the entire unit sphere.  
Among these possible connections, a random subset is chosen, so as to match the desired average degree $\langle d \rangle$.  
For each set of parameters, we compute the average of both the global bipartite clustering coefficient $C_{\protect\foursquare}^{ }$ from \cite{robins2004small} (left) and the scaled square cumulant for bipartite networks  $\tilde{\kappa}_{4\protect\foursquare}^{ }$ (right), displaying the signed fourth root for $\tilde{\kappa}_{4\protect\foursquare}^{ }$.  
%
While both clustering measures increase with $f$, the bipartite clustering coefficient $C_{\protect\foursquare}^{ }$ also notably increases with average degree, whereas $\tilde{\kappa}_{4\protect\foursquare}^{ }$ is insensitive to such changes in edge density.  
\todoscience{In addition, as $f$ approaches $0$ (the fully random bipartite limit), $\tilde{\kappa}_{4\protect\foursquare}^{ }$ also approaches $0$, indicating no propensity for clustering.} 
\todoscience{Add picture with different number of nodes too. 
Also add picture for the \mbox{6-cycle} scaled graph cumulant.  
Make distributions and get cumulants of those distributions?
put color on bipartite kappa and mus.}
}
\label{fig:ClusteringInBipartiteNetworks}
\end{center}
\end{figure}

Quantifying clustering in networks with additional features is an active domain of research that has arguably not reached a consensus.  
For example, there have been multiple proposals for weighted \cite{opsahl2009clustering,barrat2004thearchitecture,antoniou2008statistical}, directed \cite{miyajima2014continuous,fagiolo2007clustering} and bipartite networks \cite{robins2004small,zhang2008clustering,brunson2015triadic,lind2005cycles}.  
For all cases, our framework provides a principled measure of clustering, viz., the relevant scaled graph cumulant. 
For example, for directed networks, the two \mbox{third-order} scaled triangle cumulants provide two measures of clustering: one with cyclic orientation \mbox{$\tilde{\kappa}_{3\protect\threetrianglelinear}^{ }$} and one with transitive \mbox{$\tilde{\kappa}_{3\protect\threetrianglecyclic}^{ }$}. 
For bipartite networks, extensions are somewhat less straightforward, as triangles are now excluded. 
Several proposed measures consider the appearance of \mbox{4-cycles}, similarly compared to the number of incomplete cycles.
In \figreftextSI~\ref{fig:ClusteringInBipartiteNetworks}, we compare the scaled graph cumulant of the \mbox{4-cycle} subgraph \mbox{$\tilde{\kappa}_{4\protect\foursquare}^{ }$} with the clustering coefficient proposed by \cite{robins2004small}, expressed in our framework as \mbox{$C_{\protect\foursquare}^{ } = \mu_{4\foursquare}/{\mu_{3\threeline}}$}.   
Again our measure is more directly sensitive to the propensity for clustering. 
\todoscience{Say scaling with node number too when the picture with more nodes is added. Add the figure with 6-cycle scaled graph cumulant. add citation for 6-cycle  (or \mbox{6-cycles})} 
\todoscience{Put color in the bipartite kappa and mus.}

\todoprivate{*******stopped with bib checking and notes here.}

\newpage
\subsection{Networks with node attributes}
\label{SI:NodeAttributeExample}
\todoprivate{Bipartite networks are a more general case of network with node attributes...maybe write something in the beginning of the session}  

\begin{figure}[H]
\begin{center}
\centerline{\includegraphics[width=1\columnwidth]{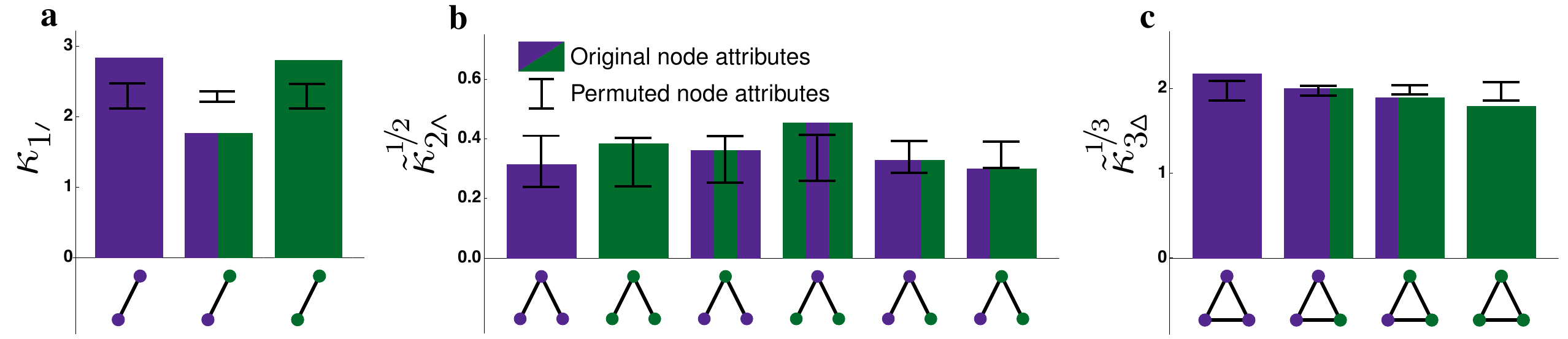}}
\caption{\textbf{Including node attributes reveals additional structure.} 
We use the scaled graph cumulants for weighted networks with a binary node attribute to analyze a social network ($222$ nodes and $5364$ edges) of interactions between primary school students from \cite{stehle2011high}, with edge weights proportional to the number of interactions between pairs of students and node attributes corresponding to their (binary) sex.\todoprivate{during one school day - maybe put back.}   
Purple nodes indicate female students, and green nodes indicate male students.  
Colored bars denote the value of (scaled) graph cumulant in the original network, and error bars denote the mean plus or minus one standard deviation for randomly shuffled node attributes
(64 runs).\todoscience{include number of runs. compare with doing same thing but for counts, run more.}  
\textbf{a)} The three \mbox{first-order} cumulants.  
\textbf{b)} The six \mbox{second-order} scaled wedge cumulants.  
\textbf{c)} The four \mbox{third-order} scaled triangle cumulants.  
The \mbox{first-order} graph cumulants (i.e., the density of edges between nodes of the indicated type) reveal a symmetric preference for homophily between the sexes, an effect \mbox{well-documented} in the social science literature \cite{wasserman1994social}. 
As all \mbox{second-order} scaled wedge cumulants are positive, we can infer a preference for hubs (i.e., nodes with degree notably larger than the average).  
Likewise, as all \mbox{third-order} scaled triangle cumulants are notably positive, we can infer a preference for triadic closure, a feature also present in many social networks \cite{watts1998collective,robins2007anintroduction}.\todoprivate{add more citation}   
While there does not appear to be much of a difference between the sexes at second order, the \mbox{third-order} scaled cumulants suggest that triadic closure is more prevalent when more participants are female.  
\todoscience{Check the lemon cumulants too.}
}
\label{fig:NodeLabel_GenderedPrimaryClass}
\end{center}
\end{figure}

\newpage
\subsection{Networks with directed edges}
\label{SI:DirectExample}

\begin{figure}[H]
\begin{center}
\centerline{\includegraphics[width=1\columnwidth]{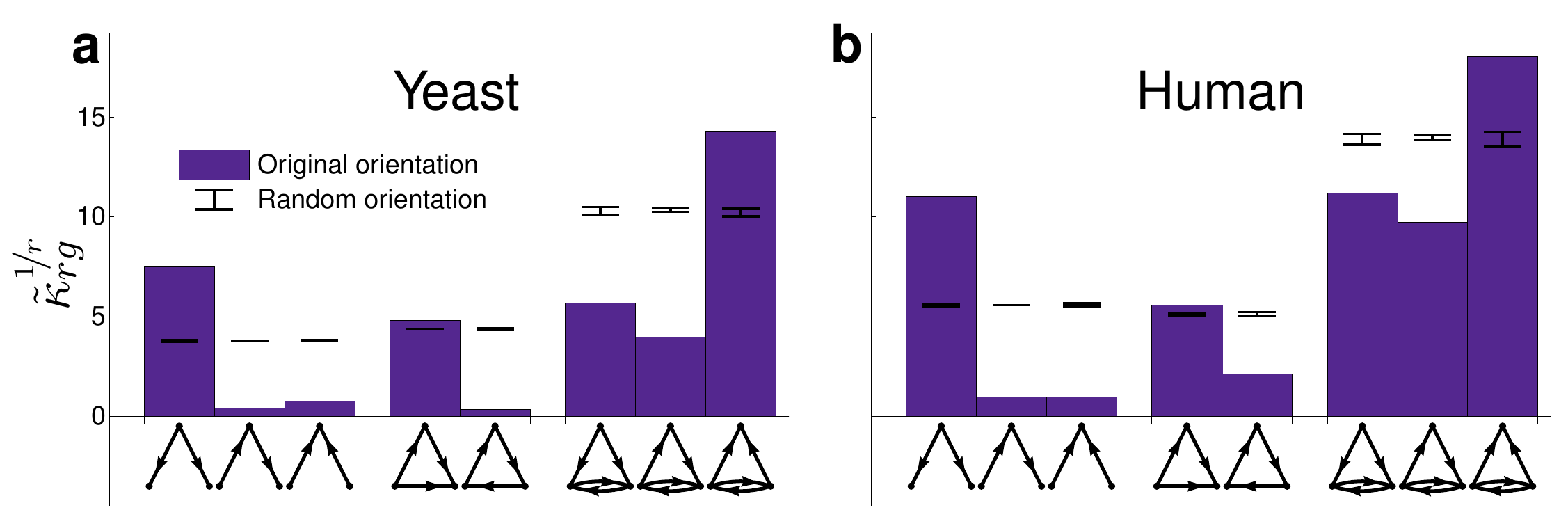}}
\caption{\textbf{The graph cumulant formalism naturally incorporates directed structure.}  
We use the scaled cumulants for directed networks to analyze regulatory networks of protein interactions from \cite{bhardwaj2010analysis} of: 
\textbf{a}) yeast ($4441$ nodes and $12873$ edges), and 
\textbf{b}) humans ($3197$ nodes and $6896$ edges).
\todoprivate{to do understand exactly what these networks were. transcriptional.}  
\todoprivate{add
%
we only display the scaled graph cumulants for substructures that are believed to be particularly important in such networks\todoprivate{citation}, 
namely, all the orientations associated with the $\tilde{\kappa}_{2\protect\twowedge}^{ }$, $\tilde{\kappa}_{3\protect\threetriangle}^{ }$, and $\tilde{\kappa}_{4\protect\fourdirecttriangle}^{ }$. 
Here, we display the scaled graph cumulants for all orientations associated with the $\tilde{\kappa}_{2\protect\twowedge}^{ }$, $\tilde{\kappa}_{3\protect\threetriangle}^{ }$, and $\tilde{\kappa}_{4\protect\fourdirecttriangle}^{ }$.}   
Error bars denote the mean and one standard deviation for the same network with randomized edge orientations.\todoprivate{include number of runs.} 
Despite the marked phenotypical differences between these two species, their regulatory networks display notable similarities.  
In particular, of the triangular substructures ($\tilde{\kappa}_{3\protect\threetriangle}^{ }$), the feedforward (or transitive) substructure ($\tilde{\kappa}_{3\protect\threetrianglelinear}^{ }$) is significantly more prevalent than the cyclic ($\tilde{\kappa}_{3\protect\threetrianglecyclic}^{ }$).\todoprivate{check literature and comment about that, cite.}   
Interestingly, while there is a higher prevalence of a central protein regulating many others ($\tilde{\kappa}_{2\protect\twowedgeoo}^{ }$), proteins that regulate each other display a propensity to both regulate the same other protein ($\tilde{\kappa}_{4\protect\fourdirectedtriangleii}^{ }$).  
}
\label{fig:DirectedStructure_Yeast}
\end{center}
\end{figure}
\todoprivate{compare what would get for directed and node attribute if did count instead*****}

\newpage
\subsection{Networks with weighted edges}
\label{SI:WeightExample}

\begin{figure}[H]
\begin{center}
\centerline{\includegraphics[width=0.80\columnwidth]{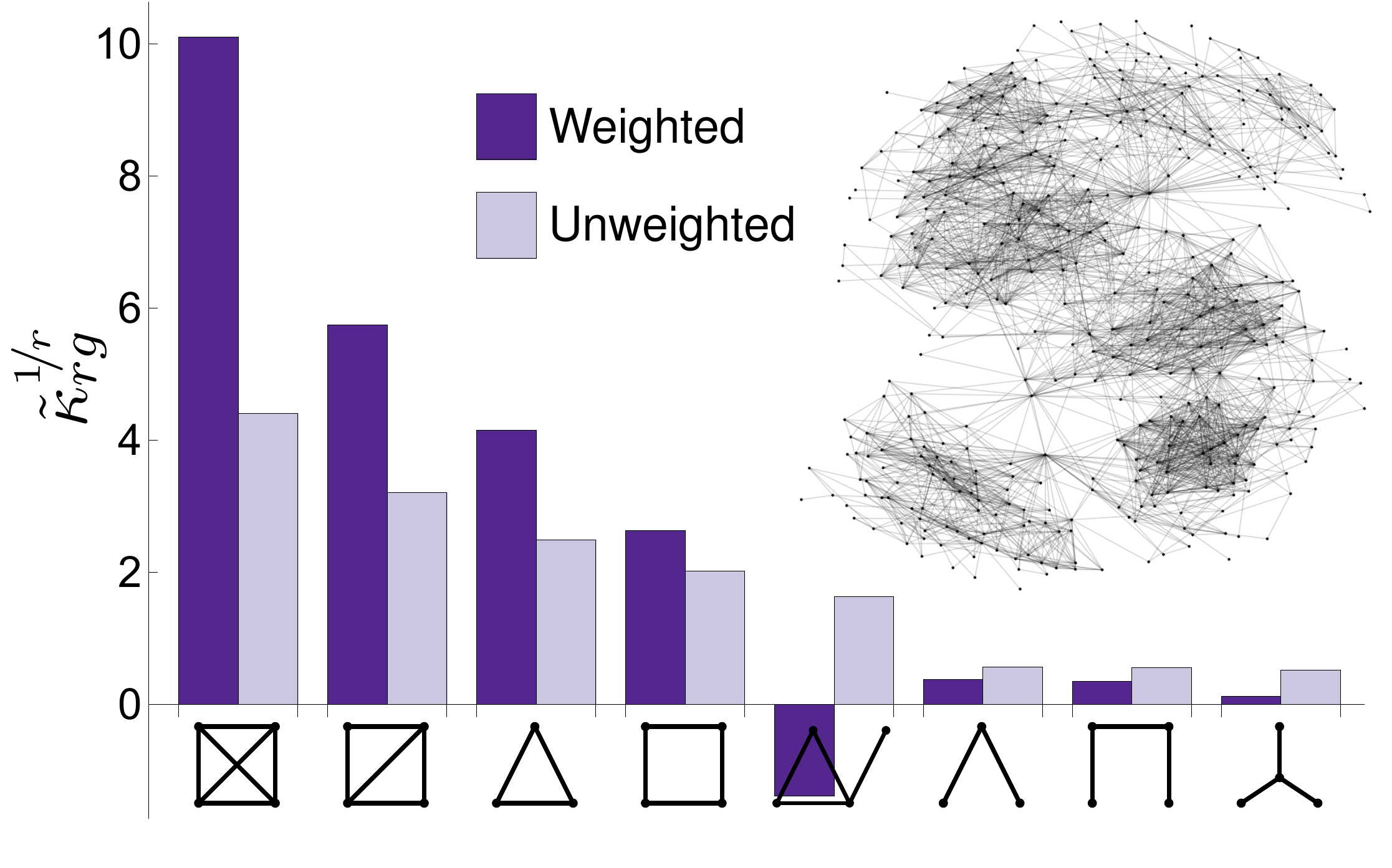}}
\caption{\textbf{Incorporating edge weights can increase the signal and possibly change the resulting interpretations.}  
We use the scaled graph cumulants for weighted networks to analyze a social network ($410$ nodes and $2765$ edges) of \mbox{face-to-face} interactions during an exhibition on infectious diseases from \cite{isella2011whats}, with edge weights proportional to the total time a pair of people spent interacting.\todoprivate{\mbox{$12$-hour}.The vertical axis corresponds to the value of the scaled graph cumulants for the weighted network (in dark purple) and for the same network but unweighted (in light purple).}  
The scaled cumulants associated with clustering (i.e., $\tilde{\kappa}_{3\protect\threetriangle}^{ }$ and $\tilde{\kappa}_{6\protect\sixKfour}^{ }$) increase when using the true edge weights (dark purple) as opposed to using the corresponding unweighted network (light purple).  
Conversely, many of the others become smaller, in particular, $\tilde{\kappa}_{2\protect\twowedge}^{ }$ and $\tilde{\kappa}_{3\protect\threeclaw}^{ }$, both of which are associated with \mbox{power-law} properties of the degree distribution.  
The most notable deviation occurs for $\tilde{\kappa}_{4\protect\fourtriangleedge}^{ }$, which is positive for the unweighted network, but negative for the weighted network.  
This negative cumulant could be interpreted as an anticorrelation between participating in triadic closure and interacting with many others, not unreasonable for such an exhibition: hosts are likely to talk with many different people, whereas groups of visitors tend to interact amongst themselves.  
\todoscience{ Maybe would be nice if we can find a weighted network for which the gradient of weight does work. Create error bars by randomizing the weights.} 
}
\label{fig:WeightedComparison_Infectious}
\end{center}
\end{figure}

\section{Unbiased graph cumulants} 
\label{SI:unbiasednetworkcumulants}

In this section, we discuss the desired properties of the unbiased estimators $\check{\kappa}_{rg}^{ }$ of graph cumulants and how to obtain them. 

\todoscience{Let's look at DCI's notes on this paragraph, could easily be consolidated when describing the new method: }
In the spirit of \mbox{$k$-statistics} \cite{fisher1930moments}, imagine a large network $\graphG_\infty^{ }$ with many nodes (the ``population''), from which one randomly subsamples \mbox{$n$} nodes, and observes the induced subgraph $\graphG$ (the ``sample'').  
We require the expectations of the unbiased graph cumulants $\check{\kappa}_{rg}^{ }$ to be invariant under this node subsampling (i.e., \mbox{$\langle \check{\kappa}_{rg}^{ } (\graphG) \rangle = \kappa_{rg}^{ } (\graphG_\infty^{ })$}), and to have the appropriate limit (i.e., \mbox{$\check{\kappa}_{rg}^{ } (\graphG) \rightarrow \kappa_{rg}^{ } (\graphG)$} as \mbox{$n \rightarrow \infty$}).  
Analogous to \mbox{real-valued} random variables, the expressions for the unbiased graph cumulants $\check{\kappa}_{rg}^{ } (\graphG)$ will be polynomials in the graph moments $\mu_{rg}^{ } (\graphG)$.  


Graph moments are preserved in expectation under random subsampling of the nodes \cite{maugis2020testing,bickel2011method}, i.e., 
\begin{align}
\langle \mu_{rg}(G) \rangle = \mu_{rg}^{ }(G_\infty^{ }). \label{Eq:MomentsPreservedInExpectation}
\end{align}
As an example, consider the expectation of $\mu_{1\oneedge}^{ }$ when removing a single random node $i$.  
Let $c_{\oneedge}^{ }$ be the counts of edges in a network with $n$ nodes, $d_i^{ }$ be the degree of node $i$, and $c'_{\oneedge}$ be the counts of edges in this network after removing node $i$ (and all its connections).  
Clearly, \mbox{$c'_{\oneedge} = c_{\oneedge}^{ } - d_i^{ }$}, so \mbox{$\langle c'_{\oneedge} \rangle = c_{\oneedge}^{ } - \langle d_i^{ } \rangle$}.  
As \mbox{$\langle d_i^{ } \rangle = \frac{2 c_{\oneedge}^{ }}{n}$}, we obtain \mbox{\smash{$\langle c'_{\oneedge} \rangle = (1-\frac{2}{n}) c_{\oneedge}^{ }$}}.  
Dividing by the edge counts in the corresponding complete networks, we have \mbox{$\langle \mu'_{1\oneedge} \rangle = \frac{{n \choose 2}}{{n-1 \choose 2}} (1-\frac{2}{n}) \mu_{1\oneedge}^{ } = \mu_{1\oneedge}^{ }$}.  
By induction, the expectation of $\mu_{1\oneedge}^{ }$ is also preserved under removal of any number of random nodes. 

Products of moments, however, are generally not preserved in expectation (under node subsampling).  
Fortunately, they can be expressed in terms of a linear combination of individual moments \cite{rucinski1988small,maugis2020testing}, 
which, as mentioned above, \textit{are} preserved in expectation.  
For example, consider $\mu_{1\oneedge}^{2}$.  
The squared counts of edges satisfies the following relation: 
\begin{align*}
{c_{\oneedge}^{ } \choose 2} = c_{\twowedge}^{ } + c_{\twoparallel}^{ },
\end{align*}
which, in turn, implies 
\begin{align}
\mu^2_{\oneedge} &= \frac{1}{\numsymbol_{\oneedge}^{2}} \Big( \numsymbol_{\oneedge} \mu^{ }_{\oneedge} + 2\numsymbol_{\twowedge} \mu^{ }_{\twowedge} + 2\numsymbol_{\twoparallel}\mu^{ }_{\twoparallel} \Big) \nonumber\\
&= \frac{2}{n(n-1)} \mu^{ }_{\oneedge} + \frac{4(n-2)}{n(n-1)} \mu^{ }_{\twowedge} + \frac{(n-2)(n-3)}{n(n-1)} \mu^{ }_{\twoparallel},\label{Eq:MuEdgeSquaredRelation}
\end{align}
where $\numsymbol_{g}^{ }$ is the count of subgraph $g$ in the complete network with $n$ nodes (e.g., the denominators of equations~\ref{eq:countedge}--\ref{eq:countk4} in \SIreftext~\ref{SI:networkcumulantslist}).  

A s moments are preserved in expectation, taking the limit as \mbox{$n\rightarrow\infty$} in equation~\ref{Eq:MuEdgeSquaredRelation} yields 
\begin{align*}
\mu^2_{\oneedge} \xrightarrow[n\rightarrow\infty]{} \mu^{ }_{\twoparallel}. 
\end{align*}
In fact, this is general, 
when a graph distribution is obtainable via sampling from a single graphon (the natural limit of a sequence of graphs with an increasing number of nodes \cite{borgs2006graph}), products of graph moments limit to the graph moment associated to their disjoint union \cite{lovasz2012large}:  
\begin{align}
\prod_i \mu_{r_i^{ }g_i^{ }}^{ } \xrightarrow[n\rightarrow\infty]{} \mu_{(\sum_i r_i^{ }) (\bigcup_i g_i^{ })}^{ }. \label{Eq:GeneralDisconnectedExpression}
\end{align}
In particular, for a single network observation, this implies that the unbiased graph cumulants are zero for all disconnected subgraphs.  

By combining relations \ref{Eq:MomentsPreservedInExpectation}~and~\ref{Eq:GeneralDisconnectedExpression}, we obtain a succinct derivation of the unbiased estimators: 
begin with the combinatorial definition of the cumulant (equation~\ref{eq:PartitionSI}), and replace all products of graph moments with a single graph moment associated with the disjoint union of the individual graphs, e.g.,
\begin{align}
\check{\kappa}_{1\oneedge}^{ } &= \mu_{1 \oneedge}^{ },\label{Eq:UnbiasedEdge}\\
\check{\kappa}_{2\twowedge}^{ } &= \mu_{2 \twowedge}^{ } - \mu_{2\twoparallel}^{ },\label{Eq:UnbiasedWedge}\\
\check{\kappa}_{2\twoparallel}^{ } &= 0,\label{Eq:UnbiasedTwoParallel}\\
\check{\kappa}_{3\threetriangle}^{ } &= \mu_{3 \threetriangle}^{ } - 3\mu_{3\threeedgewedge}^{ } + 2\mu_{3\threeparallel}^{ },\label{Eq:UnbiasedTriangle}\\
\check{\kappa}_{3\threeclaw}^{ } &= \mu_{3 \threeclaw}^{ } - 3\mu_{3\threeedgewedge}^{ } + 2\mu_{3\threeparallel}^{ },\label{Eq:UnbiasedClaw}\\
\check{\kappa}_{3\threeline}^{ } &= \mu_{3 \threeline}^{ } - 2\mu_{3\threeedgewedge}^{ } + \mu_{3\threeparallel}^{ },\label{Eq:UnbiasedThreeLine}\\
\check{\kappa}_{3\threeedgewedge}^{ } &= 0,\label{Eq:UnbiasedThreeEdgeWedge}\\
\check{\kappa}_{3\threeparallel}^{ } &= 0.\label{Eq:UnbiasedThreeParallel}
\end{align}

\todoscience{This graph limit $\graphG_{N\rightarrow\infty}^{ }$ essentially defines a graphon, where sampling from this graphon is equivalent to sampling from the desired distribution $\graphGrv$ over $n$ nodes.}  
\todoscience{This procedure naturally generalizes for networks with additional properties. When this is done for the pictures cite them.}
\todoscience{say doing inference with multiple graphs if we have them. cite our arxiv print about that.}

\section{Fitting ERGMs using unbiased graph cumulants}
\label{SI:fittingergm} 

\todoprivate{check consistency when saying constraint, target, etc.
the spectral motivation can be used to argue why all corresponding subgraphs should be taken.
maybe somewhere say We remark that fitting subgraph counts is equivalent to fit moments. We do not fit the cumulants directly, instead converting them to moments, due to the nonlinearity that it. } 

\todoprivate{Suppose one observes a single network $\graphG$, then generates a maximum entropy distribution with graph cumulants $\kappa_{rg}^{ }(\graphGrv)$ equal to these $\check{\kappa}_{rg}^{ }(\graphG)$ up to some choice of order $r'$ (i.e., an ERGM of order $r'$).  
Consider sampling single networks from this distribution and computing their $\check{\kappa}_{rg}^{ }$.  
For \mbox{$r\leq r'$}, the expectations $\langle \check{\kappa}_{rg}^{ } \rangle$ are equal to the cumulants $\kappa_{rg}^{ }$ of the distribution itself.  
Thus, the $\check{\kappa}_{rg}^{ }$ provide unbiased estimators for the graph cumulants of the ERGM distribution $\graphGrv$ from which $\graphG$ was sampled.}

We now describe how to infer a model from our proposed hierarchical family of ERGMs using a \textit{single} observed network $\graphG$.  
In particular, we consider ERGMs with prescribed expected graph moments of at most order $r$ (or, equivalently, the associated subgraph counts, as the number of nodes $n$ is fixed).\todoprivate{should we explain what a network model is?}  
These distributions have the following form:  
\begin{align} 
p(\graphG) &= \frac{1}{Z} \; \text{ER}_{n,\sfrac{1\!}{2}}(\graphG) \: \exp\!\left(\sum_{g} \beta_g^{ } c_{g}^{ }(\graphG)\right),  \label{eq:ERGM}\\
Z &= \sum_{\graphG'  \in \Omega} \left[ \text{ER}_{n,\sfrac{1\!}{2}}(\graphG') \: \exp\!\left(\sum_{g} \beta_g^{ } c_{g}^{ }(\graphG')\right) \right], 
\end{align}
where 
$\text{ER}_{n,p}(\graphG)$ is the probability of the network $\graphG$ in an \mbox{Erd\H{o}s--R\'enyi} random graph model 
(i.e., the presence of an edge between any pair of nodes occurs independently with the given probability $p$)\footnote{Just as a biased coin has maximal entropy for \mbox{$p=\sfrac{1\!}{2}$}, here, the maximum entropy distribution is given by $\text{ER}_{n,\sfrac{1\!}{2}}$, which is uniform over all labeled simple graphs with $n$ nodes, i.e., uniform over all their associated adjacency matrices.}; 
$\beta_g^{ }$ is the parameter (or Lagrange multiplier) associated with subgraph $g$; 
$c_g^{ }(\graphG)$ is the count of subgraph $g$ in the network $\graphG$; 
$\Omega$ is the space of all simple graphs with $n$ nodes; 
and $Z$ is the normalization constant (or partition function).  

\todoprivate{add old footnote: The distribution over simple graphs with $n$ nodes, in which each edge occurs independently with the same probability $p$.  
In particular, $\text{ER}_{n,\sfrac{1\!}{2}}$ is the uniform (i.e., maximum entropy) distribution over all \mbox{$n \times n$} binary symmetric traceless matrices, and thus proportional to $\frac{n!}{\text{Aut}(\graphG)}$. 
Just as a biased coin has maximal entropy for \mbox{$p=\sfrac{1\!}{2}$}, the maximum entropy distribution is given by $\text{ER}_{n,\sfrac{1\!}{2}}$ (i.e., uniform over all labeled graphs), hence its use as a base measure in our $\text{ERGM}$. add fucking the word label and stop. 
(where the presence of an edge between any pair of nodes is determined by an independent biased coin flip with the given probability $p$)
} 

While an ERGM is typically specified by the desired expectations of the statistics of interest (in this case, subgraph counts/moments), 
the parameters needed to compute the distribution are the $\beta_g^{ }$, 
and in general must be determined numerically.\todoprivate{(e.g., by maximum likelihood estimation)}  
Moreover, as the number of unique graphs grows \mbox{super-exponentially} in the number of nodes (e.g., there are $12005168$ simple graphs with $10$ nodes \cite{oeis}), the partition function $Z$ often cannot be exactly computed.   
However, there exists a large body of literature on sampling and variational techniques for efficiently approximating $Z$ as a function of $\beta_g^{ }$ \cite{snijders2006new}.\todoprivate{citemore} 



Given a single observed network, the protocol for inferring an ERGM from our hierarchical family is as follows: \\[-20pt]
\begin{enumerate}
\item Choose the order of the desired ERGM.  
Use the graph moments of observed network $\mu_{rg}^{ }$ to compute the unbiased estimators $\check{\kappa}_{rg}^{ }$ of all graph cumulants up to and including this order (see \SIreftextInSI~\ref{SI:unbiasednetworkcumulants}).\\[-20pt] 
\item Substitute these unbiased cumulants into the combinatorial formula (\mbox{equation~\ref{eq:PartitionSI}}) to obtain the desired unbiased graph moments $\check{\mu}_{rg}^{ }$ for this ERGM.\\[-20pt] 
\item Fit the parameters $\beta_g^{ }$ such that the resulting ERGM distribution has expected graph moments equal to these $\check{\mu}_{rg}^{ }$.\\[-20pt] 
\end{enumerate}

\subsection{Partially unbiased ERGMs} 
\label{SI:ergmpartiallyunbiasing}

The expressions for the unbiased graph cumulants are derived by assuming that the nodes were sampled randomly from a much larger underlying network $\graphG_\infty^{ }$.  
However, this assumption may not always be appropriate, such as when the observed network is small, or when it could feasibly represent a significant fraction of the system of interest.  
For these cases, we introduce an adjustable unbiasing parameter \mbox{$\eta \in [0,1]$}, where $1$ corresponds to the aforementioned ``fully'' unbiased case, and $0$ corresponds to using the original moments $\mu_{rg}^{ }$ of the observed network (``no unbiasing'').  
Essentially, instead of assuming that the underlying network is infinite, \mbox{$\eta = 1-\frac{n}{N}$} controls its size ($N$ nodes) relative to that of the observed network ($n$ nodes).  

The procedure is similar to before, with a modified second step.  
Instead of using the combinatorial expressions to convert $\check{\kappa}_{rg}^{ }$ to the desired moments for the ERGM $\check{\mu}_{rg}^{ }$, one inverts the expressions for the unbiased cumulants assuming a graph with \mbox{$N$} nodes.  
As the unbiased cumulants associated to disconnected subgraphs are always zero, their inverse expressions are not completely determined by the forward expressions; one must also use the expressions relating the products of graph moments for a single network, yielding, e.g.,  
\begin{align*}
\check{\mu}_{1\oneedge}^{ } &= \check{\kappa}_{1\oneedge}^{ },\\ 
\check{\mu}_{2\twowedge}^{ } &= \frac{\numsymbol_{\twoparallel}^{ }}{\numsymbol_{\twowedge}^{ }+\numsymbol_{\twoparallel}^{ }}\check{\kappa}_{2\twowedge}^{ } + \frac{\numsymbol_{\oneedge}^{2}}{2(\numsymbol_{\twowedge}^{ }+\numsymbol_{\twoparallel}^{ })}\check{\kappa}_{1\oneedge}^{2} - \frac{\numsymbol_{\oneedge}^{ }}{2(\numsymbol_{\twowedge}^{ }+\numsymbol_{\twoparallel}^{ })}\check{\kappa}_{1\oneedge}^{ },\\
 &= \frac{(N-2)(N-3)}{(N+1)(N-2)}\check{\kappa}_{2\twowedge}^{ } + \frac{N(N-1)}{(N+1)(N-2)}\check{\kappa}_{1\oneedge}^{2} - \frac{2}{(N+1)(N-2)}\check{\kappa}_{1\oneedge}^{ },\\
\check{\mu}_{2\twoparallel}^{ } &= -\frac{\numsymbol_{\twowedge}^{ }}{\numsymbol_{\twowedge}^{ }+\numsymbol_{\twoparallel}^{ }}\check{\kappa}_{2\twowedge}^{ } + \frac{\numsymbol_{\oneedge}^{2}}{2(\numsymbol_{\twowedge}^{ }+\numsymbol_{\twoparallel}^{ })}\check{\kappa}_{1\oneedge}^{2} - \frac{\numsymbol_{\oneedge}^{ }}{2(\numsymbol_{\twowedge}^{ }+\numsymbol_{\twoparallel}^{ })}\check{\kappa}_{1\oneedge}^{ }\\
 &= -\frac{4(N-2)}{(N+1)(N-2)}\check{\kappa}_{2\twowedge}^{ } + \frac{N(N-1)}{(N+1)(N-2)}\check{\kappa}_{1\oneedge}^{2} - \frac{2}{(N+1)(N-2)}\check{\kappa}_{1\oneedge}^{ }.
\end{align*}
Note that the combinatorial definitions are recovered as \mbox{$N\rightarrow\infty$}.  
\section{A geometric understanding\\of the degeneracy problem}
\label{SI:geometricdegeneracy}

The degeneracy problem refers to the appearance of undesirable \mbox{large-scale} multimodality in the distribution induced by an ERGM; 
despite the fact that averaging over this distribution gives expected subgraph counts equal to those of the observed network (as desired), 
typical samples from it have counts vastly different from these average values.\todoprivate{cite. maybe say in the ERGM derivable by a single network.maybe say for the constrained.}  

\todoprivate{** check consistency when saying constraint, target, etc}

Essentially, this arises due to the shape of the base distribution (i.e., $\text{ER}_{n,\sfrac{1\!}{2}}$) as a function of the statistics whose expected values are constrained \cite{horvat2015reducing} (here the subgraph counts, or equivalently, the corresponding graph moments).  
Recall from equation~\ref{eq:ERGM} that these ERGM distributions have the following form: 
\begin{align}
p(\graphG) &\propto \text{ER}_{n,\sfrac{1\!}{2}}(\graphG) \: \exp \!\left( \sum_{g} \beta_g^{ } c_{g}^{ }(\graphG) \right). \label{eq:ERGM2}
\end{align}
Projecting this distribution to the space of the relevant subgraph counts (i.e., summing the probability of all networks for which these counts are the same), and taking its logarithm yields: 
\begin{align} 
\ln p(\vec{c}) &= \ln\! \big( \text{ER}_{n,\sfrac{1\!}{2}}(\vec{c}) \big) + \vec{\beta} \cdot \vec{c},  \label{eq:ERGM2log}
\end{align}
where $\vec{c}\,$ is the vector of relevant subgraphs counts, $\vec{\beta}\,$ is the vector of their associated parameters, and we have dropped the term associated with the partition function (as it does not depend on $\vec{c}$). 
Thus, to understand the behavior of $p(\graphG)$, it is geometrically instructive to look at the shape of \mbox{$\,\ln\!\big(\text{ER}_{n,\sfrac{1\!}{2}}\big)$} as a function of $\vec{c}$. 

To provide intuition about the degeneracy problem and our proposed solution, here we give attention to a commonly used (and easily visualizable) 2D model, denoted by \mbox{$\text{ERGM}(\mu_{1\oneedge}^{ },\mu_{2\twowedge}^{ })$}, which prescribes the expected counts of edges and wedges in the distribution to be equal to those of the observed network.    
For comparison, we consider our \mbox{second-order} ERGM, which additionally prescribes the expected counts of pairs of edges that do not share any node. 
We will discuss both the case when the expectations of these three subgraph counts are prescribed to be those of the observed network, denoted by $\text{ERGM}(\mu_{1\protect\oneedge}^{ },\mu_{2\protect\twowedge}^{ },\mu_{2\protect\twoparallel}^{ })$, and when they are prescribed to be equal to the unbiased values (see \SIreftextInSI~\ref{SI:fittingergm}), denoted by $\text{ERGM}(\check{\mu}_{1\protect\oneedge}^{ },\check{\mu}_{2\protect\twowedge}^{ },\check{\mu}_{2\protect\twoparallel}^{ })$.

\begin{figure}[H]
\begin{center}
\centerline{\includegraphics[width=1\columnwidth]{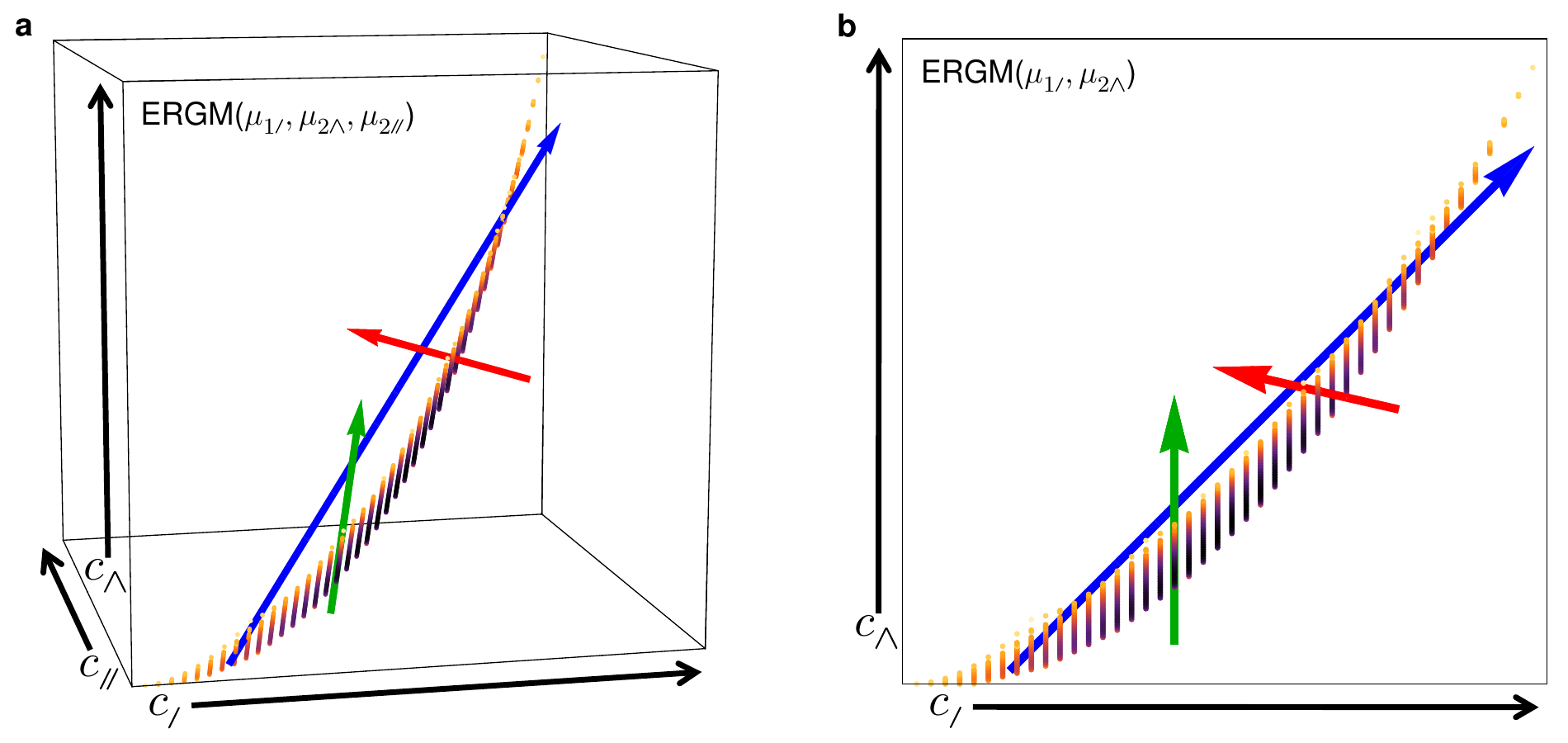}}
\caption{\textbf{Prescribing the expected counts of all subgraphs of first and second order allows for greater control over the resulting distribution.}  
Here, we explicitly enumerate all graphs with $10$ nodes. 
Each point corresponds to a tuple of subgraph counts, and the color corresponds to $\ln\!\big(\text{ER}_{10,\sfrac{1\!}{2}}\big)$, with darker colors denoting higher probability.  
\textbf{a)} When represented as a function of \mbox{($c_{\protect\oneedge}^{ },c_{\protect\twowedge}^{ },c_{\protect\twoparallel}^{ }$)}, the density of $\text{ER}_{10,\sfrac{1\!}{2}}$ lies on a 2D submanifold with extrinsic curvature.  
Three orthogonal directions can independently control the edge counts (blue arrow), the wedge counts (green arrow), and the spread in the edge counts (red arrow). 
\textbf{b)} In contrast, when representing the distribution as a function of \mbox{($c_{\protect\oneedge}^{ },c_{\protect\twowedge}^{ }$)}, these three quantities cannot be independently controlled.
This can lead to the ``degenerate'' distributions with significant bimodality observed in certain \mbox{$\text{ERGM}(\mu_{1\protect\oneedge}^{ },\mu_{2\protect\twowedge}^{ })$}. 
\todoprivate{change axes to mu, add hat. in 2d fig add center of ER and example of graph in the concave boundary. showing that the average value of er has mus that are not in the manifold.}
\todoprivate{Each point corresponds to a tuple of subgraph counts, and the color corresponds to $\ln\!\big(\text{ER}_{10,\sfrac{1\!}{2}}\big)$ (i.e., it is proportional to the number of labeled graphs on $10$ nodes with exactly these counts), with darker colors denoting higher probability.}
\todoprivate{more continuous ``orthogonal'' stuff... arrow that is plane that pushes in the entire thing.}
}
\label{Fig:GeometryManifold}
\end{center}
\end{figure}

Consider all tuples representing realizable subgraph counts of a single network in their respective 2D (for $\text{ERGM}(\mu_{1\protect\oneedge}^{ },\mu_{2\protect\twowedge}^{ })$) or 3D (for $\text{ERGM}(\mu_{1\protect\oneedge}^{ },\mu_{2\protect\twowedge}^{ },\mu_{2\protect\twoparallel}^{ })$ and $\text{ERGM}(\check{\mu}_{1\protect\oneedge}^{ },\check{\mu}_{2\protect\twowedge}^{ },\check{\mu}_{2\protect\twoparallel}^{ })$) spaces (\figreftextSI~\ref{Fig:GeometryManifold}).  
Any point within the convex hull formed by these points may serve as the prescribed expected values of some ERGM.  
However, some of these choices \textit{require} degenerate distributions. 
For example, consider \mbox{$\langle c_{\oneedge}^{ } \rangle = \frac{1}{2}\numsymbol_{\oneedge}^{ }$}, \mbox{$\langle c_{\twowedge}^{ } \rangle = \frac{1}{2}\numsymbol_{\twowedge}^{ }$} (where $\numsymbol_{g}^{ }$ is the count of subgraph $g$ in the complete graph with $n$ nodes). 
Indeed, the \textit{only} distribution with these expected values is an equal mixture of the empty and complete networks --- in a sense, the most ``degenerate'' distribution possible!  

Even if one restricts attention to tuples of subgraph counts that are realizable by a single network, \mbox{$\text{ERGM}(\mu_{1\oneedge}^{ },\mu_{2\twowedge}^{ })$} still does not always concentrate around these values.\todoprivate{prescribed expected counts.}  
In particular, this occurs when one chooses a network that lies along the concave boundary of the support of $\ln(\text{ER}_{n,\sfrac{1\!}{2}}(c_{\oneedge}^{ },c_{\twowedge}^{ }))$ (i.e., the region in \figreftextSI~\ref{Fig:GeometryManifold}b, where $c_{\twowedge}^{ }$ is large for a given number of edges).\todoprivate{***maybe a picture with the location of the concave boundary?}  
This can be understood by considering equation~\ref{eq:ERGM2log}: the $\vec{\beta}$ term (which serves to enforce the prescribed expected subgraph counts) is linear and  
essentially ``pushes'' on the distribution with the same direction and magnitude everywhere. \todoprivate{can only ``push'' the distribution linearly to enforce these constraints.}  
Thus, increasing the expected counts of wedges is inevitably coupled with a motion of the probability density toward the ``tips'' of this \mbox{crescent-shaped} domain.\todoprivate{(which contains networks with subgraph counts that are quite far from the desired average).}   
Hence, the expected counts of wedges and the \textit{spread} in the counts of edges cannot be independently controlled, and the distribution can become degenerate.

In contrast, the $\vec{\beta}$ in $\text{ERGM}(\mu_{1\protect\oneedge}^{ },\mu_{2\protect\twowedge}^{ },\mu_{2\protect\twoparallel}^{ })$ has an additional degree of freedom.  
Thus, it is able to independently control the expected counts of edges and wedges, as well as the spread in the counts of edges.  
%
%
%
However, if one requires that the expected counts $(c_{\protect\oneedge}^{ }, c_{\protect\twowedge}^{ }, c_{\protect\twoparallel}^{ })$ are exactly equal to those of the observed network, the resulting distribution necessarily concentrates on networks with precisely this edge count.\todoprivate{make consistent the counts of the.}   
%
%
Essentially, this occurs because the triplet $(c_{\protect\oneedge}^{ }, c_{\protect\twowedge}^{ }, c_{\protect\twoparallel}^{ })$ of any individual network lies on the boundary of the convex hull formed by all such triplets. 
\todoscience{we need a picture with callouts for this explanation (and the previous one too)}   
For a fixed number of edges, the relationship between the \mbox{second-order} moments is linear: \mbox{$\numsymbol_{\twowedge}^{ }\mu_{2\twowedge}^{ } + \numsymbol_{\twoparallel}^{ }\mu_{2\twoparallel}^{ } = C$} (see \figreftextSI~\ref{Fig:GeometryManifold}a).   
Additionally, the relationship between the counts of edges and this invariant sum $C$ has a curvature that does not change sign: \mbox{$C = \frac{1}{2} (\numsymbol_{\oneedge}^{2}\mu_{1\oneedge}^{2} - \numsymbol_{\oneedge}^{ }\mu_{1\oneedge}^{ })$}.\todoprivate{this is not as clear as we want...} 
%
\todoprivate{Thus, for any observed network, all distributions with expected counts that match those of this network must only have support in the linear direction corresponding to the set of triplets with the same invariant sum $C$, i.e., networks with the same number of edges as the observed network.Thus, for any observed network, all distributions with expected counts that match those of this network must only have support in the linear direction given by the set of triplets with the same invariant sum $C$. 
As this set of triplets corresponds to the networks with the same number of edges as the observed network, we have ``solved'' the degeneracy problem by essentially fixing the number of edges in the ERGM.  }   
Thus, any distribution with expected counts equal to those of an observed network must have support only in the linear direction given by the set of triplets with the same invariant sum $C$ (and therefore the same number of edges).
Thus, we have ``solved'' the degeneracy problem by essentially fixing the number of edges in the ERGM.  
However, such a solution is not satisfactory for many applications.  

\begin{figure}[H]
\begin{center}
\centerline{\includegraphics[width=0.25\columnwidth,trim={4.6cm 4.4cm 5.2cm 5.6cm},clip]{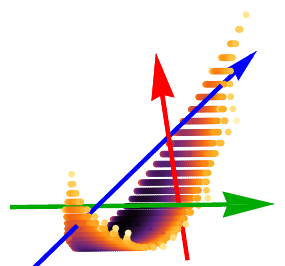}\includegraphics[width=0.25\columnwidth,trim={4.6cm 4.4cm 5.2cm 5.6cm},clip]{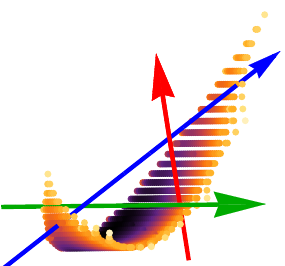}} 
\caption{\textbf{A stereographic image of the base distribution for graphs with 10 nodes.} 
Displayed is the logarithm of the base distribution $\ln(\text{ER}_{10,\sfrac{1\!}{2}}(c_{\protect\oneedge}^{ },c_{\protect\twowedge}^{ },c_{\protect\twoparallel}^{ }))$ embedded in the space of all subgraph counts up to and including second order, i.e., edges, wedges, and two edges that do not share any node. 
Darker colors indicate higher probability.  
To fully appreciate the stereographic effect, print the image using either standard A4 paper or US letter size.  
Begin with the paper close to your eyes.  
Allow your left eye to focus on the left image, while your right eye focuses on the right image.  
Slowly move the paper away from your eyes, while maintaining focus on the middle of the ``three'' images, until it is approximately \mbox{$20$--$40$ cm} (\mbox{$8$--$16$ in}) away.  
The ``middle'' image should be more apparent than those on either side, and its upper right corner should appear farther away than its bottom left.\todoprivate{We have assumed an average pupillary distance of $63$ mm, and remark that the parallax has been decreased by a factor of $\sim$$2$ to make it easier to overlay the two images.}  
See \href{https://www.youtube.com/watch?v=z9YhaK_P9d0}{\underline{here}} for an animation.
\todoprivate{explicitly chop the tail of the blue arrow} 
}
\label{Fig:Stereoscopic}
\end{center}
\end{figure}

In order to obtain a \mbox{non-degenerate} distribution containing networks with different numbers of edges, the triplet of expected counts must be slightly in the interior of the convex hull, in the direction of the red arrow in \figreftextSIs~\ref{Fig:GeometryManifold}a and~\ref{Fig:Stereoscopic}.  
The unbiased graph cumulants derived in \SIreftextInSI~\ref{SI:unbiasednetworkcumulants} provide a natural and consistent prescription for obtaining such modified triplets of expected counts (and, more generally, modified tuples of expected counts for higher order ERGMs).  
While this may seem to be an unusual choice (as such tuples are not realizable by any individual network), it is indeed quite natural: even the $\text{ER}_{n,p}$ 
distributions have tuples of expected counts that lie in this direction.\todoscience{have the picture adding ER point.}  

For a few extremal networks, our prescription for obtaining the modified expected counts may result in tuples that lie outside the convex hull, and thus do not lead to realizable ERGMs.  
This tends to occur for networks that are unlikely to be observed when subsampling nodes from a large network (such as regular or \mbox{nearly-regular} graphs).\todoprivate{(i.e., all nodes have exactly the same number of edges) (i.e., graphs in which all nodes have exactly the same number of edges)}  
From a pragmatic perspective, this is unlikely to be an issue, as \mbox{real} networks tend not to have such properties.  
Moreover, if one does observe a network for which this is the case, the issue is often alleviated by using an intermediate choice of the unbiasing parameter $\eta$ to obtain the modified tuples of expected counts (see \SIreftextInSI~\ref{SI:ergmpartiallyunbiasing}).  
\section{Statistical inference without\\constructing an explicit null model} 
\label{SI:unbiasedcumulantsergmquantification}

In order to assess the statistical significance of a network's propensity for substructures, one needs to compare the observed cumulants with the distribution obtained from some appropriate null model. 
While our proposed hierarchical family of ERGMs is a principled option, unfortunately, obtaining the parameters $\vec{\beta}$ is often computationally prohibitive. 
Fortunately, our procedure to derive the unbiased graph cumulants $\check{\kappa}_{rg}^{ }$ (\SIreftextInSI~\ref{SI:unbiasednetworkcumulants}) can also be used to derive their variance $\text{Var}(\check{\kappa}_{rg}^{ })$, allowing for statistical tests of a network's propensity for substructures \textit{without} explicitly constructing a null model.  

We first explain how to perform such a statistical test, 
then we describe how to obtain the variance of the unbiased graph cumulants. 
%

\subsection{Statistical test using \protect{$\check{\kappa}_{rg}^{ }$} and \protect{$\text{Var}(\check{\kappa}_{rg}^{ })$}} 
\label{SI:unbiasedcumulantsergmquantification_test}

\todoprivate{say we dont use scaled -- we are doing z score so it doesnt matter}
To analyze a substructure $g$ with $r$ edges, first measure the moments of the observed network up to order $2r$, and use these to compute the unbiased cumulant $\check{\kappa}_{rg}^{ }$, as well as its variance $\text{Var}(\check{\kappa}_{rg}^{ })$.  
If \mbox{$\check{\kappa}_{rg}^{ }\neq0$}, this potentially indicates a propensity (or aversiveness) for the substructure $g$, depending on the sign of $\check{\kappa}_{rg}^{ }$.\todoscience{figure out what is wrong and how.}     
To determine if such an assessment is statistically significant, one should compute the (squared) $Z$-score associated with the null hypothesis that \mbox{$\check{\kappa}_{rg}^{ }\sim0$}:  \mbox{$Z_{rg}^2 = \frac{\check{\kappa}_{rg}^{2}}{\text{Var}(\check{\kappa}_{rg}^{ })}$}. 
If \mbox{$Z_{rg}^2\gg1$}, one can be reasonably confident that the observed network has a propensity (or aversiveness) for the substructure $g$.  
This procedure can also be used to measure the similarity between two networks, by applying a \mbox{two-sample} \mbox{t-test} to the pair of unbiased cumulants associated to each particular substructure.\todoscience{look this more}   

The standard conversion from a \mbox{$Z$-score} to a \mbox{$p$-value} tacitly assumes normality.  
While this does not necessarily hold in general, the distribution of \mbox{$\check{\kappa}_{rg}^{ }$} is indeed asymptotically normal as \mbox{$n \rightarrow \infty$} \cite{bickel2011method}.\todoprivate{cite correct paper, i think sophie something has it.}  

\subsection{Deriving \protect{$\text{Var}(\check{\kappa}_{rg}^{ })$}} 
\label{SI:unbiasedcumulantsergmquantification_variance}

The variance of the unbiased graph cumulants can be obtained by exploiting the known \cite{maugis2020testing,bickel2011method} constraints on the products of graph moments (\SIreftextInSI~\ref{SI:unbiasednetworkcumulants}).   
In general, the expressions for $\text{Var}(\check{\kappa}_{rg}^{ })$ require moments up to order $2r$, as do the analogous expressions for the variance of the unbiased estimators for the cumulants of \mbox{real-valued} random variables.

We now describe this procedure, using the derivation of \mbox{$\text{Var}(\check{\kappa}_{1\oneedge}^{ })$} as an example.  
The variance of the unbiased graph cumulant associated to the edge (\mbox{$\check{\kappa}_{1\oneedge} = \mu_{1\oneedge}$} i.e., the edge density) is given by 
\begin{align}
\text{Var}(\check{\kappa}_{1\oneedge}^{ }) &=  \langle \check{\kappa}_{1\oneedge}^{2} \rangle - \langle \check{\kappa}_{1\oneedge}^{ } \rangle^2. \label{Eq:VarianceEdge} 
\end{align}
The second term is simply $\check{\mu}_{1\oneedge}^{2}$ of the distribution.  
Evaluation of the first term requires the ``product'' rule for the graph moments (see \SIreftextInSI~\ref{SI:unbiasednetworkcumulants}), in particular, we have that \mbox{$\langle\mu_{1\oneedge}^{2}\rangle = \frac{\numsymbol_{\oneedge}^{ }}{\numsymbol_{\oneedge}^{2}} \check{\mu}_{1\oneedge}^{ } + \frac{2\numsymbol_{\twowedge}^{ }}{\numsymbol_{\oneedge}^{2}} \check{\mu}_{2\twowedge}^{ } + \frac{2\numsymbol_{\twoparallel}^{ }}{\numsymbol_{\oneedge}^{2}} \check{\mu}_{2\twoparallel}^{ }$} (i.e., equation~\ref{Eq:MuEdgeSquaredRelation}).
Thus, we have:
\begin{align}
\text{Var}\big(\check{\kappa}_{1\oneedge}^{ }\big) &= \frac{\numsymbol_{\oneedge}^{ }}{\numsymbol_{\oneedge}^{2}} \check{\mu}_{1\oneedge} + \frac{2\numsymbol_{\twowedge}^{ }}{\numsymbol_{\oneedge}^{2}} \check{\mu}_{2\twowedge} + \frac{2\numsymbol_{\twoparallel}^{ }}{\numsymbol_{\oneedge}^{2}} \check{\mu}_{2\twoparallel} - \check{\mu}_{1\oneedge}^{2} \nonumber\\
&= \frac{2}{n(n-1)} \check{\mu}_{1\oneedge}^{ } + \frac{4(n-2)}{n(n-1)} \check{\mu}_{2\twowedge}^{ } + \frac{(n-2)(n-3)}{n(n-1)} \check{\mu}_{2\twoparallel}^{ } - \check{\mu}_{1\oneedge}^{2}.\label{Eq:VarEdgeLine2}
\end{align}
We note that, in fully unbiased case, \mbox{$\check{\mu}_{1\oneedge}^2 = \check{\mu}_{2\twoparallel}^{ }$}, as \mbox{$\check{\kappa}_{2\twoparallel}^{ } = 0$}.

\section{Local graph cumulants}
\label{SI:localcumulants}

Graph cumulants are statistics of the \textit{entire} network, quantifying its overall propensity for a given substructure.  
However, in some applications, such as node classification \cite{kipf2016semi,hamilton2017representation,hamilton2017inductive} and link prediction \cite{liben2007link}, one often desires statistics of the propensity of an \textit{individual} node or edge to participate in a given substructure. 
The graph cumulant framework naturally incorporates both of these ``local'' cases.\todoprivate{quotation? yes/no?}  
In this section, we describe how to derive these local graph moments and cumulants for both nodes and edges, providing the expressions necessary to compute both local triangle cumulants.

\subsection{Node local graph cumulants}
\label{SI:localcumulants_node}

The node local graph moments and cumulants are defined by giving a unique identity to the node of interest (here, symbolically distinguished by an empty circle), and applying the equations for general node attributes (see \SIreftextInSI~\ref{SI:networkcumulantslist_NodeAttributeNetworks}).  
For example, for simple graphs, there are now two \mbox{first-order} moments.  
One is defined as the count of edges between the distinguished node and any other node (i.e., the degree of the distinguished node), again normalized by the corresponding count in the associated complete graph:  
\begin{align*}
\mu_{1\oneedgelocalnodeself}^{ } &= \frac{c_{\oneedgelocalnodeself}^{ }}{n-1}. 
\end{align*}
The other \mbox{first-order} moment is defined as the count of edges that do not use the distinguished node, normalized by the corresponding count in the associated complete graph: 
\begin{align*}
\mu_{1\oneedgelocalnodeother}^{ } &= \frac{c_{\oneedgelocalnodeother}^{ }}{{n-1\choose2}}.
\end{align*} 
Likewise:
\begin{align*}
\mu_{2\twowedgelocalnodecenter}^{ } &= \frac{c_{\twowedgelocalnodecenter}^{ }}{{n-1\choose2}},\\
\mu_{2\twowedgelocalnodeend}^{ } &= \frac{c_{\twowedgelocalnodeend}^{ }}{2{n-1\choose2}}, \\
\mu_{3\threetrianglelocalnode}^{ } &= \frac{c_{\threetrianglelocalnode}^{ }}{{n-1\choose2}}.
\end{align*}

The definition of node local graph cumulants follows the same procedure as before, now taking care to incorporate the presence of this distinguished node, e.g., 
\begin{align}
\kappa_{3\threetrianglelocalnode}^{ } &= \mu_{3\threetrianglelocalnode}^{ } - \mu_{2\twowedgelocalnodecenter}^{ }\mu_{1\oneedgelocalnodeother}^{ } - 2\mu_{2\twowedgelocalnodeend}^{ }\mu_{1\oneedgelocalnodeself}^{ } + 2\mu_{1\oneedgelocalnodeself}^2\mu_{1\oneedgelocalnodeother}^{ }. \label{Eq:LocalNodeTriangleCumulant}
\end{align}
The same care must be taken when scaling the node local graph cumulants, e.g., 
\begin{align}
\tilde{\kappa}_{3\threetrianglelocalnode}^{ } &= \frac{\kappa_{3\threetrianglelocalnode}^{ }}{\mu_{1\oneedgelocalnodeself}^{2}\mu_{1\oneedgelocalnodeother}^{ }}. \label{Eq:LocalNodeTriangleCumulantScaled}
\end{align}

\subsection{Edge local graph cumulants}
\label{SI:localcumulants_edge}

A similar procedure can be used to obtain edge local graph cumulants, where instead of distinguishing a node, one now distinguishes an edge (here, represented by a \mbox{four-pointed} star at the midpoint of that edge). 
Again, there are two \mbox{first-order} edge local graph cumulants, although the one associated with the distinguished edge itself is trivial:
\begin{align}
\mu_{1\oneedgelocaledge}^{ } &\equiv 1. \nonumber
\end{align}
The other, associated with the remaining edges, is given by
\begin{align}
\mu_{1\oneedgedetachedlocaledge}^{ } &= \frac{c_{\oneedgedetachedlocaledge}^{ }}{{n\choose2} - 1},  \nonumber
\end{align}
where the floating star indicates that the distinguished edge is not included in the illustrated subgraph.  
In particular, as the nodes associated with the distinguished edge are equivalent to any other, they are neither required in nor excluded from the illustrated subgraph.  
Thus, \mbox{$c_{\oneedgedetachedlocaledge}^{ } = c_{\oneedge}^{ } - 1$}, i.e., the count of edges in the network minus the one distinguished edge.  

Likewise,
\begin{align}
\mu_{2\twowedgelocaledge}^{ } &= \frac{c_{\twowedgelocaledge}^{ }}{2(n-2)}, \nonumber\\
\mu_{2\twowedgedetachedlocaledge}^{ } &= \frac{c_{\twowedgedetachedlocaledge}^{ }}{3{n\choose3} - 2(n-2)}, \nonumber\\
\mu_{3\threetrianglelocaledge}^{ } &= \frac{c_{\threetrianglelocaledge}^{ }}{n-2}. \nonumber
\end{align}


Again, our procedure requires no modification; the edge local graph cumulants are given by a straightforward application of the combinatorial definition, e.g., 
%
\begin{align}
\kappa_{3\threetrianglelocaledge}^{ } &= \mu_{3\threetrianglelocaledge}^{ } - 2\mu_{2\twowedgelocaledge}^{ }\mu_{1\oneedgedetachedlocaledge}^{ } - \mu_{2\twowedgedetachedlocaledge}^{ }\mu_{1\oneedgelocaledge}^{ } + 2\mu_{1\oneedgelocaledge}^{ }\mu_{1\oneedgedetachedlocaledge}^2. \label{Eq:LocalEdgeTriangleCumulant}
\end{align}
Scaling the edge local graph cumulants also incorporates the distinguished edge, e.g., 
\begin{align}
\tilde{\kappa}_{3\threetrianglelocaledge}^{ } &= \frac{\kappa_{3\threetrianglelocaledge}^{ }}{\mu_{1\oneedgelocaledge}^{ }\mu_{1\oneedgedetachedlocaledge}^2}. \label{Eq:LocalEdgeTriangleCumulantScaled}
\end{align}

\todoscience{We remark that, by definition, the edge local graph moments and cumulants associated with the absence of an edge are zero. maybe talk about case of link prediction.}

\section{Graph cumulants are additive}
\label{SI:graphcumulantsadd}




\todoprivate{Entropy is additive, once realized, has become pervasive and incredibly useful (better) throughout physics and other places.  -Lee (unfinished)
Likewise, the additive nature of cumulants could bring about a similar revolution in the field of network science.  }
Essentially, the defining property of cumulants is their unique additive nature when applied to sums of independent random variables \cite{thiele1903theory,rota2000combinatorics} (e.g., \mbox{$\text{Var}(X+Y) = \text{Var}(X) + \text{Var}(Y)$} when $X$ and $Y$ are independent).  
This property is integral to foundational results in probability and statistics, such as the central limit theorem and its generalizations \cite{gnedenko1949limit,hald2000early}.  
In this section, we first define a natural notion of ``summing'' (denoted by $\oplus$) \mbox{graph-valued} random variables with the same number of nodes. 
We then show that the graph cumulants of these distributions sum when they are independent.\todoscience{make this phrase more clear}  

\begin{figure}[H]
\begin{center}
\centerline{\includegraphics[width=1\columnwidth]{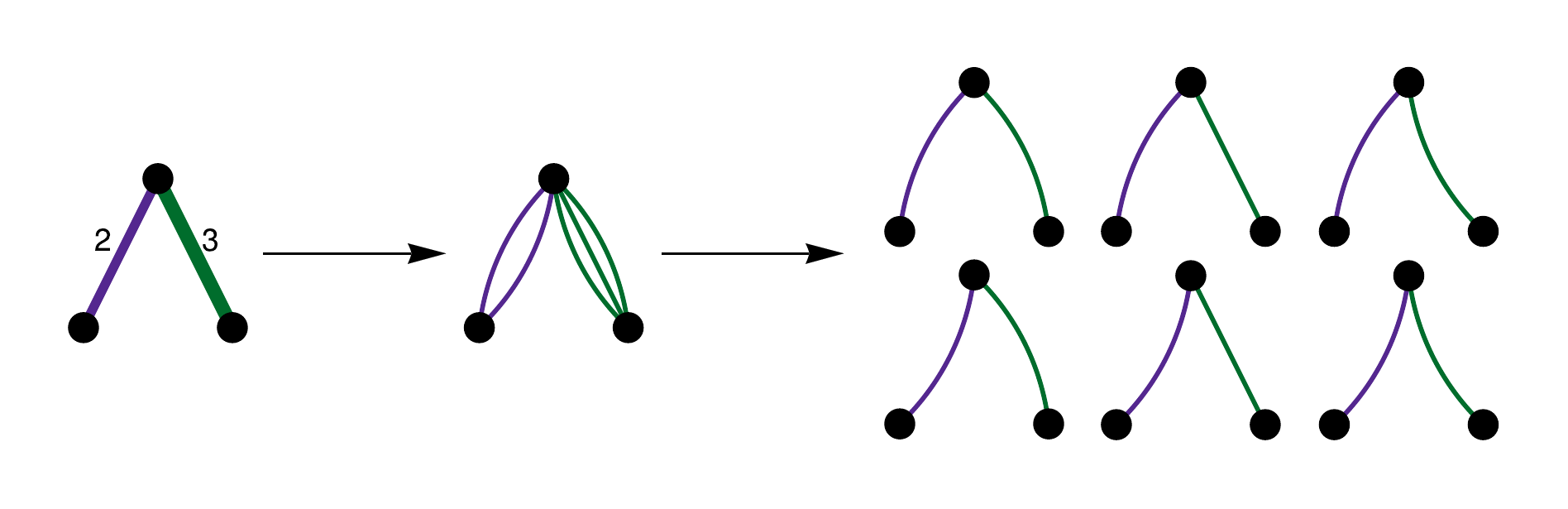}}
\caption{\textbf{Subgraphs are counted with multiplicity equal to the product of their edge weights.} 
To motivate this prescription, consider a network with integer edge weights, and represent each \mbox{integer-weighted} edge as that number of \mbox{unit-weighted} edges.    
For the weighted wedge graph shown here, by counting the number of pairs of \mbox{unit-weighted} edges that share one node, we find that there are \mbox{$2\times3$} unweighted wedges in this weighted graph. 
Indeed, for \textit{any} weighted network (including those with \mbox{real-valued} edge weights), each subgraph should be counted with weight equal to the product of its edge weights. 
}
\label{fig:SchemeWeightsWedge}
\end{center}
\end{figure}

\begin{figure}[H]
\begin{center}
\centerline{\includegraphics[width=1\columnwidth]{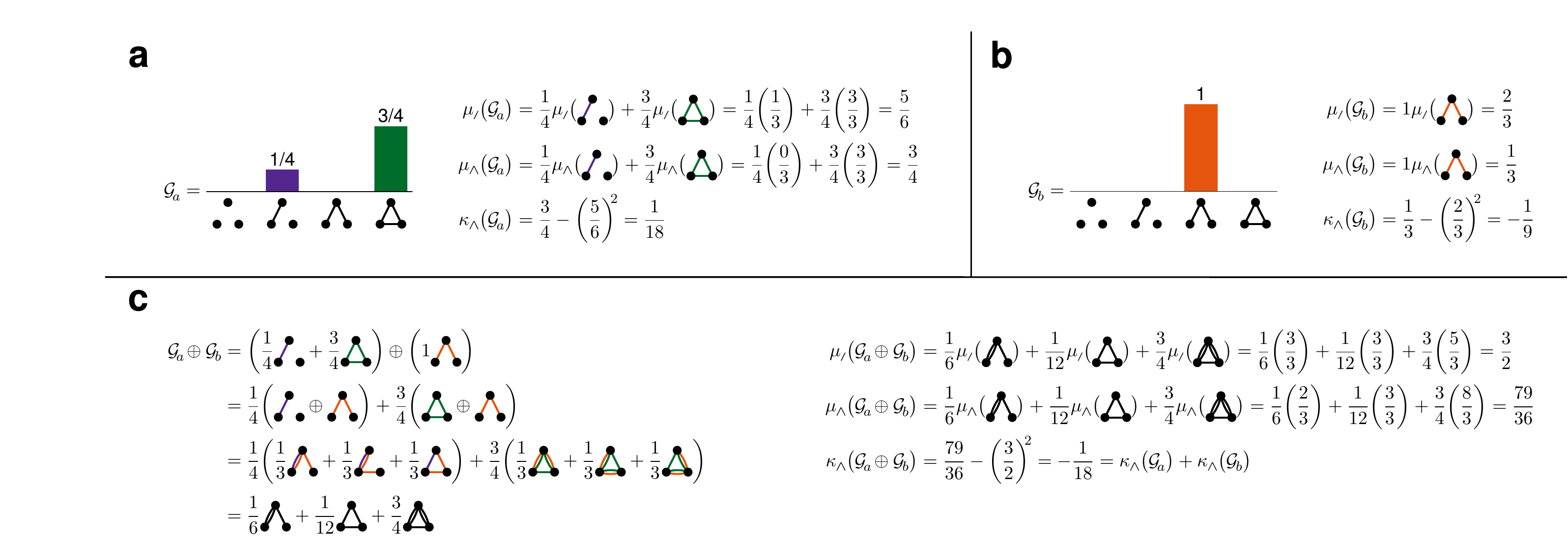}}
\caption{
\textbf{Graph cumulants have the same additive property as classical cumulants.}
\textbf{a)}~A~\mbox{graph-valued} random variable $\graphGrv$ corresponds to a probability distribution over graphs with $n$ nodes.  
Here, \mbox{$n=3$} and $\graphGrv_{\!a}^{ }$ yields a single edge with probability \mbox{$\sfrac{1\!}{4}$} and a triangle with probability \mbox{$\sfrac{3\!}{4}$}.  
\textbf{b)}~A~\mbox{graph-valued} random variable may also be concentrated on a single graph, as is the case for $\graphGrv_{\!\!\:b}^{ }$, which always yields a wedge.  
\textbf{c)}~Our graph sum $\oplus$ is defined for \mbox{graph-valued} random variables with the same number of nodes.  
Specifically, the adjacency matrices of each network are summed \mbox{component-wise}, with probability uniformly distributed over all permutations of the nodes. 
If the \mbox{graph-valued} random variables being summed are independent, the resulting (weighted) \mbox{graph-valued} random variable has cumulants equal to the sum of those of the original distributions.  
}
\label{fig:SchemeSummingGraphDistributions}
\end{center}
\end{figure}

There are a variety of operations that compose two graphs, such as the disjoint union\todoprivate{cite} and a variety of graph products \cite{nouri2012graph}.  
Here, we consider the sum of two graphs \mbox{$\graphG_{\!a}^{ }\!\!\: \oplus \graphG_{\!\!\;b}^{ }$} to be at the level of their adjacency matrices, defined by simply adding the entries \mbox{component-wise}.  
In general, as the same graph can be represented by many adjacency matrices, we must assign equal probability to each.  
In particular, for a graph $\graphG$ with $n$ nodes represented by an adjacency matrix \mbox{$A=\{a_{i\!j}^{ }\}$}, then one distributes the probability associated to this graph uniformly over all matrices \mbox{$A'=\{a_{\sigma(i)\sigma(j)}\}^{ }$} for all permutations $\sigma$ of \mbox{$\{1,\ldots, n\}$}.  
Thus, when summing two graphs, one considers all the ways that their sets of representative adjacency matrices could sum. 
The result is a \mbox{graph-valued} random variable over weighted graphs (see \figreftextSIs~\ref{fig:SchemeWeightsWedge} and~\ref{fig:SchemeSummingGraphDistributions}).\todoprivate{, where a graph-valued random variable induces a probability distribution over all adjacency matrices over $n$ nodes. 
In this sense, a single graph} 
This notion extends to \mbox{graph-valued} random variables by the distributive property,
\begin{align}
\bigg(\sum_{\graphG \in \Omega} p_{\!a}^{ }\!\!\:(\graphG)\, \graphG\bigg) \oplus \bigg(\sum_{\graphG' \in \Omega} p_b^{ }\!\!\:(\graphG')\, \graphG'\bigg) = \sum_{\graphG^{ }\in \Omega} \sum_{\graphG'\in \Omega} p_{\!a}^{ }\!\!\:(\graphG)\, p_b^{ }\!\!\:(\graphG')\, \graphG\oplus \graphG'.
\end{align}
Moreover, as $\oplus$ is clearly commutative and associative, it is also \mbox{well-defined} for multiple \mbox{graph-valued} random variables.

\todoprivate{We removed this from main text: To motivate this definition, consider an unweighted network as equivalent to a weighted complete graph with edge weight $1$ if there is an edge, and $0$ if there is not.  
The subgraph counts of this weighted complete graph are equal to those of the original unweighted network.  
When summing such networks (see \SIreftext~\ref{SI:graphcumulantsadd} for details), we essentially superimpose multiple networks on the same set of nodes, leading to the appearance of multiple edges between the same pair of nodes.  
The counts of a subgraph $g$ with $r$ edges in the resulting network is given by the number of sets of $r$ edges with the same connectivity as $g$.  
This is equivalent to treating each collection of multiple edges as a single edge with an integer weight equal to the number of edges that it represents (see 
\figreftextSI~\ref{fig:SchemeWeightsWedge}).  
This treatment naturally extends to \mbox{real-valued} edge weights \cite{lovasz2012large}.  }

Even when summing unweighted \mbox{graph-valued} random variables, the result is a \mbox{graph-valued} random variable over weighted graphs.\todoprivate{i would prefer saying graph-valued random variable. Even when unweighted \mbox{graph-valued} random variables are summed, the result is a weighted \mbox{graph-valued} random variable.}   
Thus, to obtain the graph cumulants of the resulting distribution, we must generalize the notion of subgraph density to weighted networks (itself a useful extension). 
\todoprivate{Thus, to incorporate this notion in our graph cumulants framework, we must generalize the notion of subgraph density to weighted networks (itself a useful extension). As the summation of two unweighted graph-valued random yields a distribution over weighted networks, to show that their cumulants sum, we must generalize the notion of subgraph density to weighted networks (itself a useful extension).} 
Several ways have been proposed to generalize counts of subgraphs to weighted networks \cite{miyajima2014continuous,opsahl2009clustering,barrat2004thearchitecture,antoniou2008statistical}.  
Within our framework, the consistent prescription is to treat a weighted edge as a collection of multiple edges that sum to its weight.  
Hence, when counting subgraphs, one should consider each instance with multiplicity equal to the product of its edge weights \cite{lovasz2012large}. 
The normalization for graph moments is the same as before, i.e., the counts of the subgraphs in the unweighted 
complete network (thus, the graph moments of weighted networks may be greater than one). 
Likewise, the conversion from graph moments to graph cumulants remains identical (see expressions in \SIreftextInSI~\ref{SI:networkcumulantslist_WeightedNetworks}). 

With the definitions for summing \mbox{graph-valued} random variables and for computing moments and cumulants of weighted networks, we can now state the main result of this section (see \figreftextSI~\ref{fig:SchemeSummingGraphDistributions}): 
\textit{For two independent \mbox{graph-valued} random variables over $n$ nodes, \mbox{$\graphGrv_{\!a}^{ }$ and $\graphGrv_{\!\!\:b}^{ }$}, the graph cumulants of their sum is the sum of their cumulants:}
\begin{align}
\kappa_{rg}^{ }(\graphGrv_{\!a}^{ } \!\!\:\oplus \graphGrv_{\!\!\:b}^{ }) = \kappa_{rg}^{ }(\graphGrv_{\!a}^{ }) + \kappa_{rg}^{ }(\graphGrv_{\!\!\:b}^{ }). \label{Eq:AdditiveCumulants} 
\end{align}
By induction, this holds for the sum of any number of independent \mbox{graph-valued} random variables.\todoscience{make the italic statement more precise.} 

\todoscience{comment what it means to not be independent.} 

To demonstrate how to verify this property, we first consider the specific cases of \mbox{$\kappa_{1\oneedge}^{ }$, $\kappa_{2\twowedge}^{ }$, and $\kappa_{3\threetriangle}^{ }$}, and then give the combinatorial argument for the general case. 
Clearly, for the first moment: 
\begin{align}
\mu_{1\oneedge}^{ }(\graphGrv_{\!a}^{ } \oplus \graphGrv_{\!\!\:b}^{ }) = \mu_{1\oneedge}^{ }(\graphGrv_{\!a}^{ }) + \mu_{1\oneedge}^{ }(\graphGrv_{\!\!\:b}^{ }), \label{eq:muedgesum}
\end{align}
as edge weights simply sum and the normalization remains the same.  

For \mbox{$\mu_{2\twowedge}^{ }(\graphGrv_{\!a}^{ } \!\!\:\oplus \graphGrv_{\!\!\:b}^{ })$}, one must consider the $2^2_{ }$ ways in which a wedge could be formed: both edges from $\graphGrv_{\!a}^{ }$, giving $\mu_{2\twowedge}^{ }(\graphGrv_{\!a}^{ })$; both edges from $\graphGrv_{\!\!\:b}^{ }$, giving $\mu_{2\twowedge}^{ }(\graphGrv_{\!\!\:b}^{ })$; a ``left'' edge from $\graphGrv_{\!a}^{ }$ and a ``right'' edge from $\graphGrv_{\!\!\:b}^{ }$, giving $\mu_{1\oneedge}^{ }(\graphGrv_{\!a}^{ })\mu_{1\oneedge}^{ }(\graphGrv_{\!\!\:b}^{ })$; and a ``left'' edge from $\graphGrv_{\!\!\:b}^{ }$ and a ``right'' edge from $\graphGrv_{\!a}^{ }$, giving $\mu_{1\oneedge}^{ }(\graphGrv_{\!\!\:b}^{ })\mu_{1\oneedge}^{ }(\graphGrv_{\!a}^{ })$.  
Thus, 
\begin{align}
\mu_{2\twowedge}^{ }(\graphGrv_{\!a}^{ } \!\!\:\oplus \graphGrv_{\!\!\:b}^{ }) = \mu_{2\twowedge}^{ }(\graphGrv_{\!a}^{ }) + \mu_{2\twowedge}^{ }(\graphGrv_{\!\!\:b}^{ }) + 2 \mu_{1\oneedge}^{ }(\graphGrv_{\!a}^{ })\mu_{1\oneedge}^{ }(\graphGrv_{\!\!\:b}^{ }) \label{eq:muwedgesum}
\end{align}
Substituting \ref{eq:muedgesum} and \ref{eq:muwedgesum} into the expression for $\kappa_{2\twowedge}^{ }$ (equation~ \ref{Eq:WedgeCumulant}), we have
\begin{align}
\kappa_{2\twowedge}^{ }(\graphGrv_{\!a}^{ } \!\!\:\oplus \graphGrv_{\!\!\:b}^{ }) &= \mu_{2\twowedge}^{ }(\graphGrv_{\!a}^{ } \!\!\:\oplus \graphGrv_{\!\!\:b}^{ }) - \mu_{1\oneedge}(\graphGrv_{\!a}^{ } \!\!\:\oplus \graphGrv_{\!\!\:b}^{ })^{2} \nonumber \\
 &= \big( \mu_{2\twowedge}^{ }(\graphGrv_{\!a}^{ }) + \mu_{2\twowedge}^{ }(\graphGrv_{\!\!\:b}^{ }) + 2 \mu_{1\oneedge}^{ }(\graphGrv_{\!a}^{ })\mu_{1\oneedge}^{ }(\graphGrv_{\!\!\:b}^{ }) \big) - \big( \mu_{1\oneedge}^{ }(\graphGrv_{\!a}^{ }) + \mu_{1\oneedge}^{ }(\graphGrv_{\!\!\:b}^{ }) \big)^2 \nonumber\\
&= \mu_{2\twowedge}^{ }(\graphGrv_{\!a}^{ }) + \mu_{2\twowedge}^{ }(\graphGrv_{\!\!\:b}^{ }) - \mu_{1\oneedge}^{2}(\graphGrv_{\!a}^{ }) - \mu_{1\oneedge}^{2}(\graphGrv_{\!\!\:b}^{ }) \nonumber\\
&= \kappa_{2\twowedge}^{ }(\graphGrv_{\!a}^{ }) + \kappa_{2\twowedge}^{ }(\graphGrv_{\!\!\:b}^{ }), \nonumber
\end{align}
as desired. 

Likewise, for \mbox{$\mu_{3\threetriangle}^{ }(\graphGrv_{\!a}^{ } \!\!\:\oplus \graphGrv_{\!\!\:b}^{ })$}, one must again consider the $2^3_{ }$ ways in which a triangle could be formed: all edges from $\graphGrv_{\!a}^{ }$, giving $\mu_{3\threetriangle}^{ }(\graphGrv_{\!a}^{ })$; all edges from $\graphGrv_{\!\!\:b}^{ }$, giving $\mu_{3\threetriangle}^{ }(\graphGrv_{\!\!\:b}^{ })$; a wedge from $\graphGrv_{\!a}^{ }$ and an edge from $\graphGrv_{\!\!\:b}^{ }$ (occurring for three configurations), giving \mbox{$3\mu_{2\twowedge}^{ }(\graphGrv_{\!a}^{ })\mu_{1\oneedge}^{ }(\graphGrv_{\!\!\:b}^{ })$}; and a wedge from $\graphGrv_{\!\!\:b}^{ }$ and an edge from $\graphGrv_{\!a}^{ }$ (again occurring for three configurations), giving \mbox{$3\mu_{2\twowedge}^{ }(\graphGrv_{\!\!\:b}^{ })\mu_{1\oneedge}^{ }(\graphGrv_{\!a}^{ })$}. 
Thus,
\begin{align}
\mu_{3\threetriangle}^{ }(\graphGrv_{\!a}^{ } \!\!\:\oplus \graphGrv_{\!\!\:b}^{ }) = \mu_{3\threetriangle}^{ }(\graphGrv_{\!a}^{ }) + \mu_{3\threetriangle}^{ }(\graphGrv_{\!\!\:b}^{ }) + 3\mu_{2\twowedge}^{ }(\graphGrv_{\!a}^{ })\mu_{1\oneedge}^{ }(\graphGrv_{\!\!\:b}^{ }) + 3\mu_{2\twowedge}^{ }(\graphGrv_{\!\!\:b}^{ })\mu_{1\oneedge}^{ }(\graphGrv_{\!a}^{ }). \label{eq:mutrianglesum}
\end{align}
Substituting \ref{eq:muedgesum}, \ref{eq:muwedgesum} and \ref{eq:mutrianglesum} into the expression for $\kappa_{3\threetriangle}^{ }$ (equation~ \ref{Eq:TriangleCumulant}), we again find that  
\begin{align}
\kappa_{3\threetriangle}^{ }(\graphGrv_{\!a}^{ } \!\!\:\oplus \graphGrv_{\!\!\:b}^{ }) &= \kappa_{3\threetriangle}^{ }(\graphGrv_{\!a}^{ }) + \kappa_{3\threetriangle}^{ }(\graphGrv_{\!\!\:b}^{ }). \nonumber
\end{align} 
as desired.
\todoprivate{Add proof in general. \ maybe prove it for stuff with additional structure too. maybe explain dependency and independency for graph-valued random variables.}

\todoprivate{More generally,
\begin{align}
\mu_{rg}^{ }(\graphGrv_{\!a}^{ } \oplus \graphGrv_{\!\!\:b}^{ }) &= \sum_{x\subseteq E(g)}^{ } \mu_{x}^{ }(\graphGrv_{\!a}^{ }) \mu_{E(g)\setminus x}^{ }(\graphGrv_{\!\!\:b}^{ }) \label{eq:proofsumline1}\\
 &= \sum_{x\subseteq E(g)}^{ } \Big(\sum_{\pi \in P_{x}^{ }} \prod_{b \in \pi} \kappa_{b}^{ }(\graphGrv_{\!a}^{ })\Big) \Big(\sum_{\pi \in P_{E(g)\setminus x}^{ }} \prod_{b \in \pi} \kappa_{b}^{ }(\graphGrv_{\!\!\:b}^{ })\Big)\label{eq:proofsumline2}\\
 &= \sum_{\pi \in P_{E(g)}^{ }} \sum_{x \subseteq \pi} \Big( \prod_{b \in x} \kappa_{b}^{ }(\graphGrv_{\!a}^{ }) \prod_{b \in \pi\setminus x} \kappa_{b}^{ }(\graphGrv_{\!\!\:b}^{ }) \Big)\label{eq:proofsumline3}\\
 &= \sum_{\pi \in P_{E(g)}^{ }} \prod_{b \in \pi} \Big(\kappa_{b}^{ }(\graphGrv_{\!a}^{ }) + \kappa_{b}^{ }(\graphGrv_{\!\!\:b}^{ })\Big) \label{eq:proofsumline4}\\
 &= \sum_{\pi \in P_{E(g)}^{ }} \prod_{b \in \pi} \kappa_{b}^{ }(\graphGrv_{\!a}^{ } \oplus \graphGrv_{\!\!\:b}^{ }). \label{eq:proofsumline5} 
\end{align}
where, for notational convenience, we have used $E(g)$ (the edge set of $g$) instead of the pair $rg$ for the subscripts. }

More generally,
\begin{align}
\mu_{g}^{ }(\graphGrv_{\!a}^{ } \!\!\:\oplus \graphGrv_{\!\!\:b}^{ }) &= \sum_{x\subseteq g}^{ } \mu_{x}^{ }(\graphGrv_{\!a}^{ }) \mu_{x^c_{ }}^{ }(\graphGrv_{\!\!\:b}^{ }) \label{eq:proofsumline1}\\
 &= \sum_{x\subseteq g}^{ } \Bigg(\sum_{\pi \in P_{x}^{ }} \prod_{p \in \pi} \kappa_{p}^{ }(\graphGrv_{\!a}^{ })\Bigg) \Bigg(\sum_{\pi \in P_{x^c_{ }}^{ }} \prod_{p \in \pi} \kappa_{p}^{ }(\graphGrv_{\!\!\:b}^{ })\Bigg)\label{eq:proofsumline2}\\
 &= \sum_{\pi \in P_{g}^{ }} \sum_{y \subseteq \pi} \Bigg( \prod_{p \in y} \kappa_{p}^{ }(\graphGrv_{\!a}^{ }) \prod_{p \in y^c_{ }} \kappa_{p}^{ }(\graphGrv_{\!\!\:b}^{ }) \Bigg)\label{eq:proofsumline3}\\
 &= \sum_{\pi \in P_{g}^{ }} \prod_{p \in \pi} \Bigg(\kappa_{p}^{ }(\graphGrv_{\!a}^{ }) + \kappa_{p}^{ }(\graphGrv_{\!\!\:b}^{ })\Bigg) \label{eq:proofsumline4}\\
 &= \sum_{\pi \in P_{g}^{ }} \prod_{p \in \pi} \kappa_{p}^{ }\big(\graphGrv_{\!a}^{ } \!\!\:\oplus \graphGrv_{\!\!\:b}^{ }\big). \label{eq:proofsumline5} 
\end{align}
where, for notational convenience, $g$ now represents the edge set of the subgraph (replacing the pair $rg$).  
The superscript $c$ denotes the complement with respect to the set from which this subset was taken. 
Line~\ref{eq:proofsumline1} enumerates all possible $2^{r}_{ }$ ways in which subgraph $g$ could be formed.  
Line~\ref{eq:proofsumline2} expands the moments in terms of cumulants (equation~\ref{eq:PartitionSI}).  
Line~\ref{eq:proofsumline3} exchanges the order of the summands.  
Line~\ref{eq:proofsumline4} applies the distributive law.  
Line~\ref{eq:proofsumline5} demonstrates consistency with the additive property of graph cumulants (equation~\ref{Eq:AdditiveCumulants}).

\section{Spectral motivation}
\label{SI:spectralinterpretation}

In this section, we describe a spectral motivation for graph moments and their generalizations. 

\subsection{Simple graphs}
\label{SI:spectralinterpretationundirected}

Consider the problem of parameterizing a distribution over simple graphs with $n$ nodes.  
We will represent such networks by ordered binary vectors of length ${n \choose{2}}$, where each entry represents an unordered pair of nodes, with $1$ indicating that these two nodes are connected by an edge and $0$ that they are not.  
Let $\mathcal{X}$ be the space of all such vectors.  
In general, the same graph can be represented by multiple vectors, and distributions over these graphs must give the same probability to all vectors that represent the same graph.  

When parametrizing a distribution, it is often desirable that the distribution be ``smooth'', in the sense that similar graphs are assigned similar probabilities. 
As our notion of similarity, we consider a ``graph edit distance'' \cite{gao2010survey}, defined as the minimum number of edge changes (i.e., additions or deletions) needed to transform one graph into the other.\todoprivate{(over graphs with the same number of nodes). 
For example, consider graphs with $3$ nodes, the wedge is distance $1$ from the single edge and from the triangle, and distance $2$ from the empty graph. 
Despite the fact that they have the same number of edges, the \mbox{$3$-star} is distance $2$ from the triangle and one node, as one edge needs to be removed and another edge added (in a different position).}   
For example, the wedge is distance $2$ from the empty graph, and distance $1$ from both the single edge and the triangle. 

One common method for parameterizing smooth functions is via a Fourier representation, i.e., in terms of the eigenfunctions of some Laplacian operator.  
For example, a \mbox{low-pass} filter is equivalent to giving preference to the \mbox{low-frequency} (i.e., \mbox{long-wavelength}) terms, where the location of the cutoff determines the smoothness of the output.  
To obtain similarly smooth parameterizations over the space of networks, we consider a Laplacian operator based on this graph edit distance.  

\begin{figure}[H]
\begin{center}
\centerline{\includegraphics[width=\columnwidth]{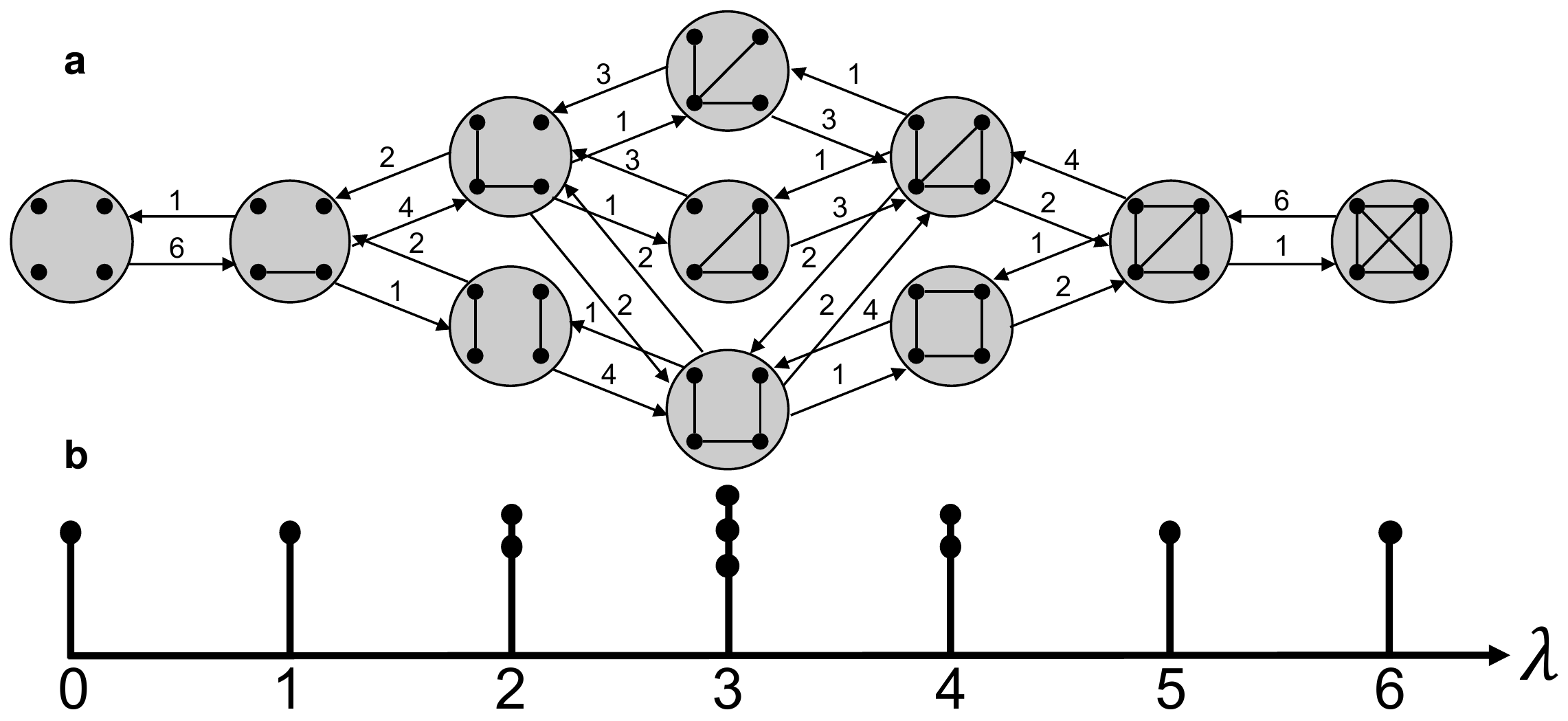}}
\caption{\textbf{Spectral motivation for graph moments.}
\textbf{a)} Directed weighted ``edit graph'' $H_4^{ }$ for simple networks with $4$ nodes. 
The large nodes represent each of the $11$ unique (\mbox{non-isomorphic}) networks.   
A directed edge from node $i$ to node $j$ indicates that network $i$ can be transformed into network $j$ by adding or removing an edge at a single location, where the weight of this directed edge is equal to the number of locations that could accomplish this transformation.  
Note that every node has an out-degree of $6$, corresponding to the ${4 \choose 2}$ possible locations.  
\textbf{b)} Schematic of the spectrum of the Laplacian of $H_4^{ }$. 
The eigenvalues are integers ranging from $0$ to ${4 \choose 2}$, and are degenerate with multiplicity equal to the number of distinct subgraphs (with at most $4$ nodes) with that number of edges.
This correspondence has a combinatorial interpretation: the span of left eigenvectors with eigenvalue \mbox{$\lambda \leq r$} is precisely the span of subgraph counts up to and including order $r$.  \todoscience{redo this picture...}
}
\label{fig:hamminggraph}
\end{center}
\end{figure}

To this end, we define a (weighted, directed) ``edit graph'' $H_n^{ }$, with nodes representing unique (i.e., \mbox{non-isomorphic}) networks with $n$ nodes. 
A directed edge from one node in $H_n^{ }$ to another appears whenever the network it represents can be transformed into the other by adding or removing an edge at a single location.  
The weight of an edge in $H_n^{ }$ is given by the number of locations that could be altered to effect this transformation (see Fig.~\ref{fig:hamminggraph} for the case of simple graphs with $4$ nodes).  
The Laplacian of $H_n^{ }$ is defined as \mbox{$L_{H}^{ } = D_{\textrm{out}}^{ } - A_H^{\top}$}, where $D_{\textrm{out}}^{ }$ is the diagonal matrix of \mbox{out-degrees}, and \mbox{$A_H^{ }=\{a_{i\!j}^{ }\}$} is the (asymmetric) adjacency matrix with entries $a_{i\!j}^{ }$ equal to the weight of the transition from $i$ to $j$. 

The lowest eigenvalue of $L_H^{ }$ is $0$, and is associated with a left eigenvector that is uniform over the \textit{unique} networks and a right eigenvector that is uniform over all representations of the networks (i.e., uniform over all binary vectors of length ${n \choose{2}}$), thus corresponding to the $\textrm{ER}_{n,\sfrac{1\!}{2}}$ distribution.  
The remainder of the spectrum contains additional structure.  
Its support is the set of positive integers up to and including ${n \choose{2}}$, and each integer has a predictable degeneracy: 
the multiplicity of an eigenvalue \mbox{$\lambda = r$} is equal to the number of distinct subgraphs with exactly $r$ edges (with at most $n$ nodes).  
This is not just a combinatorial coincidence; the span of left eigenvectors (with eigenvalue at most $r$) is precisely the span of subgraph counts (of order at most $r$) in each network.\todoprivate{in networks with $n$ nodes.}    

The structure of this spectrum gives rise to a hierarchical parameterization of distributions over networks that is equivalent to our proposed family of hierarchical ERGMs, namely
\begin{equation}
p(\graphG) \propto v_{d,0}^{ }(\graphG) \; \exp \!\left( \sum_{i=1}^{r_{\vphantom{i}}^{ }} \sum_{j=1}^{r_{\!i}^{ }} \beta_{i,j}^{ } v_{l,i,j}^{ }(\graphG) \right), \label{eq:ERGMwithEigeinvalue}
\end{equation}
where $v_{d,0}^{ }$ is the right eigenvector of $L_H^{ }$ with eigenvalue $0$ (i.e., the $\textrm{ER}_{n,\sfrac{1\!}{2}}$ base distribution); 
$v_{l,i}^{ }$ is the set of left eigenvectors of  $L_H^{ }$ with eigenvalue $i$, and $v_{l,i,j}^{ }$ is one such vector; 
$r_i^{ }$ is the number of eigenvectors with eigenvalue $i$; and $\beta_{i,j}^{ }$ are the parameters to be determined.  
In particular, as eigenvectors with the same eigenvalue are intrinsically intertwined, this degeneracy offers a principled motivation for the use of \textit{all} subgraphs up to some chosen order. 


\subsection{Generalizations}
\label{SI:spectralinterpretationgeneralization}

\todoprivate{the spectral motivation can be used to argue why all corresponding subgraphs should be taken.} 

We now generalize the concept of moments and cumulants for distributions over a set \mbox{$\mathcal{X} \equiv \mathcal{A}^\ell$}, i.e., the vectors of length \mbox{$\ell \in \mathbb{N}$} over the alphabet \mbox{$\mathcal{A} = \left\{ a_1, a_2, \ldots \right\}$}, invariant with respect to a group $\groupG$ acting on this set $\mathcal{X}$. 
The action of $\groupG$ induces an equivalence relation on \mbox{$\mathcal{X}$}: \mbox{$x\sim x' \Leftrightarrow \exists \groupg \in \groupG \: | \: x=\groupg\circ x'$}, partitioning it into orbits.
The distribution over $\mathcal{X}$ is then characterized by assigning a probability to each of these orbits.  
For example, for the case of simple graphs, 
the group $\groupG$ is $S_n$, acting by permuting the $n$ nodes of a graph.  
The set $\mathcal{X}$ consists of ordered binary vectors of length \mbox{$\ell = {n \choose 2}$}, where each entry represents an unordered pair of nodes, with $1$ indicating that these two nodes are connected by an edge and $0$ that they are not.  
Nonisomorphic graphs are in different orbits, and all the elements in a given orbit correspond to the same graph, with the probability associated to that orbit distributed uniformly over all of its elements.

With this framework, we can construct the weighted directed ``edit graph'' described in the previous section for an arbitrary set $\mathcal{X}$ and a group $\groupG$ acting upon it. 
We can then use the spectrum of the Laplacian of this edit graph to obtain the number of moments at each order.  
Again, the nodes of the edit graph are the orbits of $\mathcal{X}$ under the action of $\groupG$,  
and a directed edge from orbit $i$ to orbit $j$ indicates that an element in orbit $i$ can be transformed into an element in orbit $j$ by changing one of its $\ell$ entries. 
The weight of this directed edge is given by the number of elements in orbit $j$ that differ by a single entry from any single fixed element in orbit $i$.

This abstraction applies to a variety of situations, and naturally encompasses the generalizations previously presented in this paper. 
For example, for unweighted directed networks with no \mbox{self-loops}, $\mathcal{X}$ is the set of all ordered binary vectors of length \mbox{$\ell = 2{n \choose 2}$}, where each entry represents an \textit{ordered} pair of nodes, with $1$ indicating that there is an edge from the first node to the second and $0$ that there is not.  
The group $\groupG$ is again $S_n$.  
As another example, consider the case of undirected unweighted bipartite networks, i.e., every node has one of two possible ``flavors'' (``charm'' and ``strange''), and edges can occur only between nodes of different flavors. 
The set $\mathcal{X}$ consists of all ordered binary vectors of length \mbox{$\ell = n_{\textrm{ch}}^{ }n_{\textrm{str}}^{ }$}, where each entry represents a different unordered pair of nodes of different flavors, and $1$ indicates that these two nodes are connected by an edge and $0$ that they are not.  
The group $\groupG$ allows for permutations of nodes of the same flavor, namely \mbox{$S_{n_{\textrm{ch}}^{ }} \!\times S_{n_{\textrm{str}}^{ }}$}.\todoscience{figure out homogenization with dots..., (of/with???)}

We now illustrate the versatility of this formalism by describing an additional generalization, namely, \mbox{$k$-uniform} hypergraphs, i.e., where each 
(hyper)edge represents a connection between $k$ distinct nodes. 
As in the standard graph case (i.e., \mbox{$k=2$}), the group $\groupG$ is the symmetric group $S_n$ acting by permuting the $n$ nodes.  
The set $\mathcal{X}$ consists of all ordered binary vectors of length \mbox{$\ell = {n\choose k}$}, where each entry represents an unordered set of $k$ nodes, and a $1$ indicates the presence of a hyperedge between them and $0$ its absence. 
The orbits induced by this group action again partition the elements of $\mathcal{X}$ into equivalence classes, one for each of the unique hypergraphs.\todoprivate{where nonisomorphic hypergraphs are in different equivalence classes. or corresponding to each of the nonisomorphic hypergraphs. corresponding to each unique hypergraph. where each class corresponds to a unique hypergraph.}  
The eigenvalues of the Laplacian of the corresponding edit graph follow a similar pattern: 
associated to the eigenvalue of $0$ is a left eigenvector that is uniform over unique hypergraphs, and a right eigenvector that is uniform over all elements of $\mathcal{X}$. 
Likewise, for the remaining left eigenvectors, 
there is one eigenvector with associated eigenvalue of $1$ that is linear in the number of hyperedges.  
At second order (associated eigenvalue of $2$), there are now \textit{k} eigenvalues (for \mbox{$n\geq 2k$}), corresponding to the $k$ ways that two hyperedges can relate (sharing any number from $0$ to $k-1$ nodes). 

\section{Formulas for graph moments\\ and graph cumulants}
\label{SI:networkcumulantslist} 

Here, we provide the expressions for graph moments and graph cumulants used to obtain the results presented in this paper.  
These expressions are also included explicitly in our associated code, and we have automated their derivation to arbitrary order.  
%
We first give the normalizations for obtaining the graph moments as well as the expressions for efficiently computing the disconnected subgraph counts (see \SIreftextInSI~\ref{SI:ConnectedToDisconnected}).
We then give the expressions for their conversion to graph cumulants (by inverting equation~\ref{eq:PartitionSI}).

\todoprivate{Here, we present expressions for computing some relevant graph moments of a network $\graphG$ (with $n$ nodes) using only the counts $c_{g}^{ }$ of the \textit{connected} subgraphs $g$ in $\graphG$ (see \SIreftextInSI~\ref{SI:ConnectedToDisconnected}). 
We also present their conversion to the graph cumulants (by using equation~\ref{eq:PartitionSI}, see also \figreftext~\ref{SI:unbiasednetworkcumulants}). 
While we do not display all expressions here, we have automated their derivation to arbitrary order, and the results up to sixth order are included explicitly in the code associated with this paper. We first give the normalizations for obtaining the graph moments as well as the expressions for obtaining their counts in terms of the counts of the connected subgraphs only (see \SIreftextInSI~\ref{SI:ConnectedToDisconnected}). To encourage an efficient implementation, we express the graph moments in terms of the counts of the \textit{connected} subgraphs only (see \SIreftextInSI~\ref{SI:ConnectedToDisconnected}), and then present their conversion to graph cumulants (by inverting equation~\ref{eq:PartitionSI}, see also \figreftext~\ref{SI:unbiasednetworkcumulants}). } 

\subsection{Undirected, unweighted networks}
\label{SI:networkcumulantslist_SimpleNetworks}

For simple graphs, we now enumerate the expressions for all \mbox{third-order} graph moments and cumulants, as well as those that are necessary for computing the (\mbox{sixth-order}) cumulant associated with the complete graph with four nodes. 
The remaining expressions up to and including sixth order are explicitly included on our code. 

\subsubsection{Graph moments} 

\begingroup
\allowdisplaybreaks
\begin{align}
\mu_{1\oneedge}^{ } &= \frac{c_{\oneedge}^{ }}{{n\choose 2}}\label{eq:countedge}\\
\mu_{2\twowedge}^{ } &= \frac{c_{\twowedge}^{ }}{3{n\choose 3}}\label{eq:countwedge}\\
\mu_{2\twoparallel}^{ } &= \frac{c_{\twoparallel}^{ }}{3{n\choose 4}} = \frac{{c_{\oneedge}^{ } \choose 2} - c_{\twowedge}^{ }}{3 {n\choose 4}}\label{eq:counttwoparallels}\\
\mu_{3\threetriangle}^{ } &= \frac{c_{\threetriangle}^{ }}{{n\choose 3}}\label{eq:counttriangle} \\
\mu_{3\threeclaw}^{ } &= \frac{c_{\threeclaw}^{ }}{4{n\choose 4}}\label{eq:countclaw} \\
\mu_{3\threeline}^{ } &= \frac{c_{\threeline}^{ }}{12{n\choose 4}} \label{eq:counthreeline} \\
\mu_{3\threeedgewedge}^{ } &= \frac{c_{\threeedgewedge}^{ }}{30{n\choose 5}} = \frac{c_{\twowedge}^{ }(c_{\oneedge}^{ } - 2) - 3 c_{\threetriangle}^{ } - 3 c_{\threeclaw}^{ } - 2 c_{\threeline}^{ }}{30{n\choose 5}} \\
\mu_{3\threeparallel}^{ } &= \frac{c_{\threeparallel}^{ }}{15{n\choose 6}} = \frac{{c_{\oneedge}^{ } \choose 3} - c_{\threetriangle}^{ } - c_{\threeclaw}^{ } - c_{\threeline} - c_{\threeedgewedge}^{ }}{15{n\choose 6}}\\
\mu_{4\fourtriangleedge}^{ } &= \frac{c_{\fourtriangleedge}^{ }}{12{n\choose 4}}\\
\mu_{4\foursquare}^{ } &= \frac{c_{\foursquare}^{ }}{3{n\choose 4}}\\
\mu_{5\fivediagsquare}^{ } &= \frac{c_{\fivediagsquare}^{ }}{6{n\choose 4}}\\
\mu_{6\sixKfour}^{ } &= \frac{c_{\sixKfour}^{ }}{{n\choose 4}}\label{eq:countk4}
\end{align}
\endgroup

\subsubsection{Graph cumulants} 


\begingroup
\allowdisplaybreaks
\begin{align}
\kappa_{1\oneedge}^{ } &= \mu_{1\oneedge}^{ } \label{Eq:EdgeCumulant} \\
\kappa_{2\twowedge}^{ } &= \mu_{2\twowedge}^{ } - \mu_{1\oneedge}^{2} \label{Eq:WedgeCumulant}\\
\kappa_{2\twoparallel}^{ } &= \mu_{2\twoparallel}^{ } - \mu_{1\oneedge}^{2} \\
\kappa_{3\threetriangle}^{ } &= \mu_{3\threetriangle}^{ } - 3\mu_{2\twowedge}^{ }\mu_{1\oneedge}^{ } + 2\mu_{1\oneedge}^{3} \label{Eq:TriangleCumulant} \\
\kappa_{3\threeclaw}^{ } &= \mu_{3\threeclaw}^{ } - 3\mu_{2\twowedge}^{ }\mu_{1\oneedge}^{ } + 2\mu_{1\oneedge}^{3} \\
\kappa_{3\threeline}^{ } &= \mu_{3\threeline}^{ } - 2\mu_{2\twowedge}^{ }\mu_{1\oneedge}^{ } - \mu_{2\twoparallel}^{ }\mu_{1\oneedge}^{ } + 2\mu_{1\oneedge}^{3} \\
\kappa_{3\threeedgewedge}^{ } &= \mu_{3\threeedgewedge}^{ } - \mu_{2\twowedge}^{ }\mu_{1\oneedge}^{ } - 2\mu_{2\twoparallel}^{ }\mu_{1\oneedge}^{ } + 2\mu_{1\oneedge}^{3} \\
\kappa_{3\threeparallel}^{ } &= \mu_{3\threeparallel}^{ } - 3\mu_{2\twoparallel}^{ }\mu_{1\oneedge}^{ } + 2\mu_{1\oneedge}^{3} \\
\kappa_{4\fourtriangleedge}^{ } &= \mu_{4\fourtriangleedge}^{ } - \mu_{3\threetriangle}^{ } \mu_{1\oneedge}^{ } - \mu_{3\threeclaw}^{ } \mu_{1\oneedge}^{ } - 2\mu_{3\threeline}^{ } \mu_{1\oneedge}^{ } \nonumber\\
&\hphantom{=} - 2\mu_{2\twowedge}^{2} - \mu_{2\twowedge}^{ }\mu_{2\twoparallel}^{ } + 10\mu_{2\twowedge}^{ }\mu_{1\oneedge}^{2} + 2\mu_{2\twoparallel}^{ }\mu_{1\oneedge}^{2} - 6\mu_{1\oneedge}^{4} \\
\kappa_{4\foursquare}^{ } &= \mu_{4\foursquare}^{ } - 4\mu_{3\threeline}^{ } \mu_{1\oneedge}^{ } - 2\mu_{2\twowedge}^{2} - \mu_{2\twoparallel}^{2} + 8\mu_{2\twowedge}^{ }\mu_{1\oneedge}^{2} + 4\mu_{2\twoparallel}^{ }\mu_{1\oneedge}^{2} - 6\mu_{1\oneedge}^{4} \\
\kappa_{5\fivediagsquare}^{ } &= \mu_{5\fivediagsquare}^{ } - 4\mu_{4\fourtriangleedge}^{ }\mu_{1\oneedge}^{ } - \mu_{4\foursquare}^{ }\mu_{1\oneedge}^{ } - 2\mu_{3\threetriangle}^{ }\mu_{2\twowedge}^{ } - 2\mu_{3\threeclaw}^{ }\mu_{2\twowedge}^{ } - 4\mu_{3\threeline}^{ }\mu_{2\twowedge}^{ } - 2\mu_{3\threeline}^{ }\mu_{2\twoparallel}^{ } \nonumber\\
&\hphantom{=} + 4\mu_{3\threetriangle}^{ }\mu_{1\oneedge}^{2} + 4\mu_{3\threeclaw}^{ }\mu_{1\oneedge}^{2} + 8\mu_{3\threeline}^{ }\mu_{1\oneedge}^{2} + 4\mu_{3\threeline}^{ }\mu_{1\oneedge}^{2} \nonumber\\
&\hphantom{=} + 20\mu_{2\twowedge}^{2}\mu_{1\oneedge}^{ } + 8\mu_{2\twowedge}^{ }\mu_{2\twoparallel}^{ }\mu_{1\oneedge}^{ } + 2\mu_{2\twoparallel}^{2}\mu_{1\oneedge}^{ } \nonumber\\
&\hphantom{=} - 48\mu_{2\twowedge}^{ }\mu_{1\oneedge}^{3} - 12\mu_{2\twoparallel}^{ }\mu_{1\oneedge}^{3} + 24\mu_{1\oneedge}^{5} \\
\kappa_{6\sixKfour}^{ } &= \mu_{6\sixKfour}^{ } - 6\mu_{5\fivediagsquare}^{ }\mu_{1\oneedge}^{ } - 12\mu_{4\fourtriangleedge}^{ }\mu_{2\twowedge}^{ } - 3\mu_{4\foursquare}^{ }\mu_{2\twoparallel}^{ } + 24\mu_{4\fourtriangleedge}^{ }\mu_{1\oneedge}^{2} + 6\mu_{4\foursquare}^{ }\mu_{1\oneedge}^{2} \nonumber\\
&\hphantom{=} - 4\mu_{3\threetriangle}^{ }\mu_{3\threeclaw}^{ } - 6\mu_{3\threeline}^{2} \nonumber\\
&\hphantom{=} + 24\mu_{3\threetriangle}^{ }\mu_{2\twowedge}^{ }\mu_{1\oneedge}^{ } + 24\mu_{3\threeclaw}^{ }\mu_{2\twowedge}^{ }\mu_{1\oneedge}^{ } + 48\mu_{3\threeline}^{ }\mu_{2\twowedge}^{ }\mu_{1\oneedge}^{ } + 24\mu_{3\threeline}^{ }\mu_{2\twoparallel}^{ }\mu_{1\oneedge}^{ } \nonumber\\
&\hphantom{=} - 24\mu_{3\threetriangle}^{ }\mu_{1\oneedge}^{3} - 24\mu_{3\threeclaw}^{ }\mu_{1\oneedge}^{3} - 72\mu_{3\threeline}^{ }\mu_{1\oneedge}^{3} \nonumber\\
&\hphantom{=} + 12\mu_{2\twowedge}^{3} + 15\mu_{2\twowedge}^{2}\mu_{2\twoparallel}^{ } + 3\mu_{2\twoparallel}^{3} \nonumber\\
&\hphantom{=} - 153\mu_{2\twowedge}^{2}\mu_{1\oneedge}^{2} - 90\mu_{2\twowedge}^{ }\mu_{2\twoparallel}^{ }\mu_{1\oneedge}^{2} - 27\mu_{2\twoparallel}^{2}\mu_{1\oneedge}^{2} \nonumber\\
&\hphantom{=} + 288\mu_{2\twowedge}^{ }\mu_{1\oneedge}^{4} + 72\mu_{2\twoparallel}^{ }\mu_{1\oneedge}^{4} - 120\mu_{1\oneedge}^{6}
\end{align}
\endgroup

\subsection{Directed networks}
\label{SI:networkcumulantslist_DirectedNetworks}

We now enumerate the expressions necessary for computing the graph moments and cumulants of all directed subgraphs with three nodes, including the \mbox{sixth-order} graph cumulant associated with the complete directed triad. 
The remaining expressions up to and including fifth order are explicitly included on our code. 

\subsubsection{Graph moments} 

\begingroup
\allowdisplaybreaks
\begin{align}
\mu_{1\oneedgedir}^{ } &= \frac{c_{1\oneedgedir}^{ }}{2{n\choose2}} \\
\mu_{2\twowedgeii}^{ } &= \frac{c_{2\twowedgeii}^{ }}{3{n\choose3}} \\
\mu_{2\twowedgeoo}^{ } &= \frac{c_{2\twowedgeoo}^{ }}{3{n\choose3}} \\
\mu_{2\twowedgeio}^{ } &= \frac{c_{2\twowedgeio}^{ }}{6{n\choose3}} \\
\mu_{2\tworeciprocaldir}^{ } &= \frac{c_{2\tworeciprocaldir}^{ }}{{n\choose2}} \\
\mu_{3\threereciprocalwedgei}^{ } &= \frac{c_{3\threereciprocalwedgei}^{ }}{6{n\choose3}} \\
\mu_{3\threereciprocalwedgeo}^{ } &= \frac{c_{3\threereciprocalwedgeo}^{ }}{6{n\choose3}} \\
\mu_{3\threetrianglelinear}^{ } &= \frac{c_{3\threetrianglelinear}^{ }}{6{n\choose3}} \\
\mu_{3\threetrianglecyclic}^{ } &= \frac{c_{3\threetrianglecyclic}^{ }}{2{n\choose3}} \\
\mu_{4\fourdirectedtriangleii}^{ } &= \frac{c_{4\fourdirectedtriangleii}^{ }}{3{n\choose3}} \\
\mu_{4\fourdirectedtriangleoo}^{ } &= \frac{c_{4\fourdirectedtriangleoo}^{ }}{3{n\choose3}} \\
\mu_{4\fourdirectedtriangleio}^{ } &= \frac{c_{4\fourdirectedtriangleio}^{ }}{6{n\choose3}} \\
\mu_{4\fourreciprocalreciprocaldir}^{ } &= \frac{c_{4\fourreciprocalreciprocaldir}^{ }}{3{n\choose3}} \\
\mu_{5\fivedirectedtriangledir}^{ } &= \frac{c_{5\fivedirectedtriangledir}^{ }}{6{n\choose3}} \\
\mu_{6\sixdirectedtriangledir}^{ } &= \frac{c_{6\sixdirectedtriangledir}^{ }}{{n\choose3}}
\end{align}
\endgroup

%

\subsubsection{Graph cumulants} 

\begingroup
\allowdisplaybreaks
\begin{align}
\kappa_{1\oneedgedir}^{ } &= \mu_{1\oneedgedir}^{ }\\
\kappa_{2\twowedgeii}^{ } &= \mu_{2\twowedgeii}^{ } - \mu_{1\oneedgedir}^{2}\\
\kappa_{2\twowedgeoo}^{ } &= \mu_{2\twowedgeoo}^{ } - \mu_{1\oneedgedir}^{2}\\
\kappa_{2\twowedgeio}^{ } &= \mu_{2\twowedgeio}^{ } - \mu_{1\oneedgedir}^{2}\\
\kappa_{2\tworeciprocaldir}^{ } &= \mu_{2\tworeciprocaldir}^{ } - \mu_{1\oneedgedir}^{2}\\
\kappa_{3\threereciprocalwedgei}^{ } &= \mu_{3\threereciprocalwedgei}^{ } - \mu_{2\twowedgeii}^{ } \mu_{1\oneedgedir}^{ } - \mu_{2\twowedgeio}^{ } \mu_{1\oneedgedir}^{ } - \mu_{2\tworeciprocaldir}^{ } \mu_{1\oneedgedir}^{ } + 2 \mu_{1\oneedgedir}^{3} \\
\kappa_{3\threereciprocalwedgeo}^{ } &= \mu_{3\threereciprocalwedgeo}^{ } - \mu_{2\twowedgeoo}^{ } \mu_{1\oneedgedir}^{ } - \mu_{2\twowedgeio}^{ } \mu_{1\oneedgedir}^{ } - \mu_{2\tworeciprocaldir}^{ } \mu_{1\oneedgedir}^{ } + 2 \mu_{1\oneedgedir}^{3} \\
\kappa_{3\threetrianglelinear}^{ } &= \mu_{3\threetrianglelinear}^{ } - \mu_{2\twowedgeii}^{ } \mu_{1\oneedgedir}^{ } - \mu_{2\twowedgeoo}^{ } \mu_{1\oneedgedir}^{ } - \mu_{2\twowedgeio}^{ } \mu_{1\oneedgedir}^{ } + 2 \mu_{1\oneedgedir}^{3} \\
\kappa_{3\threetrianglecyclic}^{ } &= \mu_{3\threetrianglecyclic}^{ } - 3 \mu_{2\twowedgeio}^{ } \mu_{1\oneedgedir}^{ } + 2 \mu_{1\oneedgedir}^{3} \\
\kappa_{4\fourdirectedtriangleii}^{ } &= \mu_{4\fourdirectedtriangleii}^{ } - 2 \mu_{3\threereciprocalwedgeo}^{ } \mu_{1\oneedgedir}^{ } - 2 \mu_{3\threetrianglelinear}^{ } \mu_{1\oneedgedir}^{ } - \mu_{2\twowedgeoo}^{2} - \mu_{2\twowedgeio}^{2} - \mu_{2\twowedgeii}^{ } \mu_{2\tworeciprocaldir}^{ } \nonumber\\
&\quad+ 2 \mu_{2\twowedgeii}^{ }\mu_{1\oneedgedir}^{2} + 4 \mu_{2\twowedgeoo}^{ }\mu_{1\oneedgedir}^{2} + 4 \mu_{2\twowedgeio}^{ }\mu_{1\oneedgedir}^{2} + 2 \mu_{2\tworeciprocaldir}^{ }\mu_{1\oneedgedir}^{2} - 6 \mu_{1\oneedgedir}^{4}\\
\kappa_{4\fourdirectedtriangleoo}^{ } &= \mu_{4\fourdirectedtriangleoo}^{ } - 2 \mu_{3\threereciprocalwedgei}^{ } \mu_{1\oneedgedir}^{ } - 2 \mu_{3\threetrianglelinear}^{ } \mu_{1\oneedgedir}^{ } - \mu_{2\twowedgeii}^{2} - \mu_{2\twowedgeio}^{2} - \mu_{2\twowedgeoo}^{ } \mu_{2\tworeciprocaldir}^{ } \nonumber\\
&\quad+ 4 \mu_{2\twowedgeii}^{ }\mu_{1\oneedgedir}^{2} + 2 \mu_{2\twowedgeoo}^{ }\mu_{1\oneedgedir}^{2} + 4 \mu_{2\twowedgeio}^{ }\mu_{1\oneedgedir}^{2} + 2 \mu_{2\tworeciprocaldir}^{ }\mu_{1\oneedgedir}^{2} - 6 \mu_{1\oneedgedir}^{4}\\
\kappa_{4\fourdirectedtriangleio}^{ } &= \mu_{4\fourdirectedtriangleio}^{ } - \mu_{3\threereciprocalwedgei}^{ } \mu_{1\oneedgedir}^{ } - \mu_{3\threereciprocalwedgeo}^{ } \mu_{1\oneedgedir}^{ } - \mu_{3\threetrianglelinear}^{ } \mu_{1\oneedgedir}^{ } - \mu_{3\threetrianglecyclic}^{ } \mu_{1\oneedgedir}^{ } - \mu_{2\twowedgeii}^{ }\mu_{2\twowedgeio}^{ } - \mu_{2\twowedgeoo}^{ }\mu_{2\twowedgeio}^{ } - \mu_{2\twowedgeio}^{ } \mu_{2\tworeciprocaldir}^{ } \nonumber\\
&\quad+ 2 \mu_{2\twowedgeii}^{ }\mu_{1\oneedgedir}^{2} + 2 \mu_{2\twowedgeoo}^{ }\mu_{1\oneedgedir}^{2} + 6 \mu_{2\twowedgeio}^{ }\mu_{1\oneedgedir}^{2} + 2 \mu_{2\tworeciprocaldir}^{ }\mu_{1\oneedgedir}^{2} - 6 \mu_{1\oneedgedir}^{4}\\
\kappa_{4\fourreciprocalreciprocaldir}^{ } &= \mu_{4\fourreciprocalreciprocaldir}^{ } - 2 \mu_{3\threereciprocalwedgei}^{ } \mu_{1\oneedgedir}^{ } - 2 \mu_{3\threereciprocalwedgeo}^{ } \mu_{1\oneedgedir}^{ } - \mu_{2\twowedgeii}^{ }\mu_{2\twowedgeoo}^{ } - \mu_{2\twowedgeio}^{2} - \mu_{2\tworeciprocaldir}^{2} \nonumber\\
&\quad+ 2 \mu_{2\twowedgeii}^{ }\mu_{1\oneedgedir}^{2} + 2 \mu_{2\twowedgeoo}^{ }\mu_{1\oneedgedir}^{2} + 4 \mu_{2\twowedgeio}^{ }\mu_{1\oneedgedir}^{2} + 4 \mu_{2\tworeciprocaldir}^{ }\mu_{1\oneedgedir}^{2} - 6 \mu_{1\oneedgedir}^{4}\\
\kappa_{5\fivedirectedtriangledir}^{ } &= \mu_{5\fivedirectedtriangledir}^{ } - \mu_{4\fourdirectedtriangleii}^{ } \mu_{1\oneedgedir}^{ } - \mu_{4\fourdirectedtriangleoo}^{ } \mu_{1\oneedgedir}^{ } - 2 \mu_{4\fourdirectedtriangleio}^{ } \mu_{1\oneedgedir}^{ } - \mu_{4\fourreciprocalreciprocaldir}^{ } \mu_{1\oneedgedir}^{ } - \nonumber\\
&\quad- \mu_{3\threereciprocalwedgei}^{ }\mu_{2\twowedgeii}^{ } - \mu_{3\threereciprocalwedgei}^{ }\mu_{2\twowedgeio}^{ } - \mu_{3\threereciprocalwedgei}^{ }\mu_{2\tworeciprocaldir}^{ } - \mu_{3\threereciprocalwedgeo}^{ }\mu_{2\twowedgeoo}^{ } - \mu_{3\threereciprocalwedgeo}^{ }\mu_{2\twowedgeio}^{ } \nonumber\\
&\quad- \mu_{3\threereciprocalwedgeo}^{ }\mu_{2\tworeciprocaldir}^{ } - \mu_{3\threetrianglelinear}^{ }\mu_{2\twowedgeii}^{ } - \mu_{3\threetrianglelinear}^{ }\mu_{2\twowedgeoo}^{ } - \mu_{3\threetrianglelinear}^{ }\mu_{2\twowedgeio}^{ } - \mu_{3\threetrianglecyclic}^{ }\mu_{2\twowedgeio}^{ } \nonumber\\
&\quad + 6 \mu_{3\threereciprocalwedgei}^{ }\mu_{1\oneedgedir}^{2} + 6 \mu_{3\threereciprocalwedgeo}^{ }\mu_{1\oneedgedir}^{2} + 9 \mu_{3\threetrianglelinear}^{ }\mu_{1\oneedgedir}^{2} + 2 \mu_{3\threetrianglecyclic}^{ }\mu_{1\oneedgedir}^{2} \nonumber\\
&\quad + 2 \mu_{2\twowedgeii}^{2} + 2 \mu_{2\twowedgeoo}^{2} + 6 \mu_{2\twowedgeio}^{2} + 2 \mu_{2\tworeciprocaldir}^{2} + 2 \mu_{2\twowedgeii}^{ }\mu_{2\twowedgeoo}^{ } \nonumber\\
&\quad + 4 \mu_{2\twowedgeii}^{ }\mu_{2\twowedgeio}^{ } + 4 \mu_{2\twowedgeoo}^{ }\mu_{2\twowedgeio}^{ } + 2 \mu_{2\twowedgeii}^{ }\mu_{2\tworeciprocaldir}^{ } + 2 \mu_{2\twowedgeoo}^{ }\mu_{2\tworeciprocaldir}^{ } + 4 \mu_{2\twowedgeio}^{ }\mu_{2\tworeciprocaldir}^{ } \nonumber\\
&\quad - 12 \mu_{2\twowedgeii}^{ }\mu_{1\oneedgedir}^{3} - 12 \mu_{2\twowedgeoo}^{ }\mu_{1\oneedgedir}^{3} - 24 \mu_{2\twowedgeio}^{ }\mu_{1\oneedgedir}^{3} - 12 \mu_{2\tworeciprocaldir}^{ }\mu_{1\oneedgedir}^{3} 
+ 24 \mu_{1\oneedgedir}^{5} \\
\kappa_{6\sixdirectedtriangledir}^{ } &= \mu_{6\sixdirectedtriangledir}^{ } - 6 \mu_{5\fivedirectedtriangledir}^{ }\mu_{1\oneedgedir}^{ } - 3 \mu_{4\fourdirectedtriangleii}^{ }\mu_{2\twowedgeoo}^{ } - 3 \mu_{4\fourdirectedtriangleoo}^{ }\mu_{2\twowedgeii}^{ } - 6 \mu_{4\fourdirectedtriangleio}^{ }\mu_{2\twowedgeio}^{ } - 3 \mu_{4\fourreciprocalreciprocaldir}^{ }\mu_{2\tworeciprocaldir}^{ } \nonumber\\
&\quad + 6 \mu_{4\fourdirectedtriangleii}^{ }\mu_{1\oneedgedir}^{2} + 6 \mu_{4\fourdirectedtriangleoo}^{ }\mu_{1\oneedgedir}^{2} + 12 \mu_{4\fourdirectedtriangleio}^{ }\mu_{1\oneedgedir}^{2} + 6 \mu_{4\fourreciprocalreciprocaldir}^{ }\mu_{1\oneedgedir}^{2} 
- 3 \mu_{3\threereciprocalwedgei}^{2} - 3 \mu_{3\threereciprocalwedgeo}^{2} - 3 \mu_{3\threetrianglelinear}^{2} - \mu_{3\threetrianglecyclic}^{2} \nonumber\\
&\quad + 12 \mu_{3\threereciprocalwedgei}^{ }\mu_{2\twowedgeii}^{ }\mu_{1\oneedgedir}^{ } + 12 \mu_{3\threereciprocalwedgei}^{ }\mu_{2\twowedgeio}^{ }\mu_{1\oneedgedir}^{ } + 12 \mu_{3\threereciprocalwedgei}^{ }\mu_{2\tworeciprocaldir}^{ }\mu_{1\oneedgedir}^{ } + 12 \mu_{3\threereciprocalwedgeo}^{ }\mu_{2\twowedgeoo}^{ }\mu_{1\oneedgedir}^{ } + 12 \mu_{3\threereciprocalwedgeo}^{ }\mu_{2\twowedgeio}^{ }\mu_{1\oneedgedir}^{ } \nonumber\\
&\quad + 12 \mu_{3\threereciprocalwedgeo}^{ }\mu_{2\tworeciprocaldir}^{ }\mu_{1\oneedgedir}^{ } + 12 \mu_{3\threetrianglelinear}^{ }\mu_{2\twowedgeii}^{ }\mu_{1\oneedgedir}^{ } + 12 \mu_{3\threetrianglelinear}^{ }\mu_{2\twowedgeoo}^{ }\mu_{1\oneedgedir}^{ } + 12 \mu_{3\threetrianglelinear}^{ }\mu_{2\twowedgeio}^{ }\mu_{1\oneedgedir}^{ } + 12 \mu_{3\threetrianglecyclic}^{ }\mu_{2\twowedgeio}^{ }\mu_{1\oneedgedir}^{ } \nonumber\\
&\quad - 36 \mu_{3\threereciprocalwedgei}^{ }\mu_{1\oneedgedir}^{3} - 36 \mu_{3\threereciprocalwedgeo}^{ }\mu_{1\oneedgedir}^{3} - 36 \mu_{3\threetrianglelinear}^{ }\mu_{1\oneedgedir}^{3} - 12 \mu_{3\threetrianglecyclic}^{ }\mu_{1\oneedgedir}^{3} \nonumber\\
&\quad + 2 \mu_{2\twowedgeii}^{3} + 6 \mu_{2\twowedgeii}^{ }\mu_{2\twowedgeoo}^{ }\mu_{2\tworeciprocaldir}^{ } + 6 \mu_{2\twowedgeii}^{ }\mu_{2\twowedgeio}^{2} + 2 \mu_{2\twowedgeoo}^{3} + 6 \mu_{2\twowedgeoo}^{ }\mu_{2\twowedgeio}^{2} + 6 \mu_{2\twowedgeio}^{2}\mu_{2\tworeciprocaldir}^{ } + 2 \mu_{2\tworeciprocaldir}^{3} \nonumber\\
&\quad - 18 \mu_{2\twowedgeii}^{2}\mu_{1\oneedgedir}^{2} - 18 \mu_{2\twowedgeii}^{ }\mu_{2\twowedgeoo}^{ }\mu_{1\oneedgedir}^{2} - 36 \mu_{2\twowedgeii}^{ }\mu_{2\twowedgeio}^{ }\mu_{1\oneedgedir}^{2} - 18 \mu_{2\twowedgeii}^{ }\mu_{2\tworeciprocaldir}^{ }\mu_{1\oneedgedir}^{2} - 18 \mu_{2\twowedgeoo}^{2}\mu_{1\oneedgedir}^{2} \nonumber\\
&\quad - 36 \mu_{2\twowedgeoo}^{ }\mu_{2\twowedgeio}^{ }\mu_{1\oneedgedir}^{2} - 18 \mu_{2\twowedgeoo}^{ }\mu_{2\tworeciprocaldir}^{ }\mu_{1\oneedgedir}^{2} - 54 \mu_{2\twowedgeio}^{2}\mu_{1\oneedgedir}^{2} - 36 \mu_{2\twowedgeio}^{ }\mu_{2\tworeciprocaldir}^{ }\mu_{1\oneedgedir}^{2} - 18 \mu_{2\tworeciprocaldir}^{2}\mu_{1\oneedgedir}^{2} \nonumber\\
&\quad + 60 \mu_{2\twowedgeii}^{ }\mu_{1\oneedgedir}^{4} + 60 \mu_{2\twowedgeoo}^{ }\mu_{1\oneedgedir}^{4} + 120 \mu_{2\twowedgeio}^{ }\mu_{1\oneedgedir}^{4} + 60 \mu_{2\tworeciprocaldir}^{ }\mu_{1\oneedgedir}^{4} - 120 \mu_{1\oneedgedir}^{6}
\end{align}
\endgroup

\subsection{Networks with node attributes}
\label{SI:networkcumulantslist_NodeAttributeNetworks}


Often, networks have additional attributes associated with the nodes.  
By incorporating this information into the graph cumulant formalism, one can reveal structure that is correlated with these attributes.  
Our example in Figure~\ref{fig:NodeLabel_GenderedPrimaryClass} considers cumulants associated with subgraphs containing up to three nodes for a network with a binary node attribute.  
%
We now enumerate the expressions for computing the graph moments and cumulants required for this analysis (as well as for the \mbox{$3$-star} subgraphs).  
Note that the mapping from attributes to colors is arbitrary; the colors may be reversed in any expression (e.g., the expression for $\kappa_{2\TwoPurpleGreenPurple}^{ }$ can be obtained from that for $\kappa_{2\TwoGreenPurpleGreen}^{ }$ by exchanging all instances of purple and green with each other).  


\subsubsection{Graph moments}

\begingroup
\allowdisplaybreaks
\begin{align}
\mu_{1\OnePurpleToPurpleEdge}^{ } &= \frac{c_{1\OnePurpleToPurpleEdge}^{ }}{{n_{\PurpleDot}^{ } \choose 2}}\\
\mu_{1\OnePurpleToGreenEdge}^{ } &= \frac{c_{1\OnePurpleToGreenEdge}^{ }}{n_{\PurpleDot}^{ } n_{\GreenDot}^{ }}\\
\mu_{2\TwoPurplePurplePurple}^{ } &= \frac{c_{2\TwoPurplePurplePurple}^{ }}{3{n_{\PurpleDot}^{ } \choose 3}}\\
\mu_{2\TwoPurplePurpleGreen}^{ } &= \frac{c_{2\TwoPurplePurpleGreen}^{ }}{2{n_{\PurpleDot}^{ } \choose 2}n_{\GreenDot}^{ }}\\
\mu_{2\TwoGreenPurpleGreen}^{ } &= \frac{c_{2\TwoGreenPurpleGreen}^{ }}{n_{\PurpleDot}^{ }{n_{\GreenDot}^{ } \choose 2}}\\
\mu_{3\ThreePurplePurplePurple}^{ } &= \frac{c_{3\ThreePurplePurplePurple}^{ }}{{n_{\PurpleDot}^{ } \choose 3}}\\
\mu_{3\ThreePurplePurpleGreen}^{ } &= \frac{c_{3\ThreePurplePurpleGreen}^{ }}{{n_{\PurpleDot}^{ } \choose 2}n_{\GreenDot}^{ }}\\
\mu_{3\ThreeClawPurplePurplePurple}^{ } &= \frac{c_{3\ThreeClawPurplePurplePurple}^{ }}{4{n_{\PurpleDot}^{ } \choose 4}}\\
\mu_{3\ThreeClawPurplePurpleGreen}^{ } &=  \frac{c_{3\ThreeClawPurplePurpleGreen}^{ }}{3{n_{\PurpleDot}^{ } \choose 3}n_{\GreenDot}^{ }}\\
\mu_{3\ThreeClawPurpleGreenGreen}^{ } &= \frac{c_{3\ThreeClawPurpleGreenGreen}^{ }}{{2{n_{\PurpleDot}^{ } \choose 2}{n_{\GreenDot}^{ } \choose 2}}}\\
\mu_{3\ThreeClawGreenGreenGreen}^{ } &= \frac{c_{3\ThreeClawGreenGreenGreen}^{ }}{n_{\PurpleDot}^{ }{n_{\GreenDot}^{ } \choose 3}}
\end{align}
\endgroup

\subsubsection{Graph cumulants}

\begingroup
\allowdisplaybreaks
\begin{align}
\kappa_{1\OnePurpleToPurpleEdge}^{ } &= \mu_{1\OnePurpleToPurpleEdge}^{ } \\
\kappa_{1\OnePurpleToGreenEdge}^{ } &= \mu_{1\OnePurpleToGreenEdge}^{ } \\
\kappa_{2\TwoPurplePurplePurple}^{ } &= \mu_{2\TwoPurplePurplePurple}^{ } - \mu_{1\OnePurpleToPurpleEdge}^2 \\
\kappa_{2\TwoPurplePurpleGreen}^{ } &= \mu_{2\TwoPurplePurpleGreen}^{ } - \mu_{1\OnePurpleToPurpleEdge}^{ }\mu_{1\OnePurpleToGreenEdge}^{ }\\
\kappa_{2\TwoGreenPurpleGreen}^{ } &= \mu_{2\TwoGreenPurpleGreen}^{ } - \mu_{1\OnePurpleToGreenEdge}^2 \\
\kappa_{3\ThreePurplePurplePurple}^{ } &= \mu_{3\ThreePurplePurplePurple}^{ } - 3\mu_{2\TwoPurplePurplePurple}^{ }\mu_{1\OnePurpleToPurpleEdge}^{ } + 2\mu_{1\OnePurpleToPurpleEdge}^3 \\
\kappa_{3\ThreePurplePurpleGreen}^{ } &= \mu_{3\ThreePurplePurpleGreen}^{ } - 2\mu_{2\TwoPurplePurpleGreen}^{ }\mu_{1\OnePurpleToGreenEdge}^{ }  - \mu_{2\TwoPurpleGreenPurple}^{ }\mu_{1\OnePurpleToPurpleEdge}^{ } + 2\mu_{1\OnePurpleToPurpleEdge}^{ }\mu_{1\OnePurpleToGreenEdge}^2 \\
\kappa_{3\ThreeClawPurplePurplePurple}^{ } &= \mu_{3\ThreeClawPurplePurplePurple}^{ } - 3\mu_{2\TwoPurplePurplePurple}^{ }\mu_{1\OnePurpleToPurpleEdge}^{ } + 2\mu_{1\OnePurpleToPurpleEdge}^{3} \\
\kappa_{3\ThreeClawPurplePurpleGreen}^{ } &= \mu_{3\ThreeClawPurplePurpleGreen}^{ } - 2\mu_{2\TwoPurplePurpleGreen}^{ }\mu_{1\OnePurpleToPurpleEdge}^{ } - \mu_{2\TwoPurplePurplePurple}^{ }\mu_{1\OnePurpleToGreenEdge}^{ } + 2\mu_{1\OnePurpleToPurpleEdge}^{2}\mu_{1\OnePurpleToGreenEdge}^{ } \\
\kappa_{3\ThreeClawPurpleGreenGreen}^{ } &= \mu_{3\ThreeClawPurpleGreenGreen}^{ } - 2\mu_{2\TwoPurplePurpleGreen}^{ }\mu_{1\OnePurpleToGreenEdge}^{ } - \mu_{2\TwoGreenPurpleGreen}^{ }\mu_{1\OnePurpleToPurpleEdge}^{ } + 2\mu_{1\OnePurpleToPurpleEdge}^{ }\mu_{1\OnePurpleToGreenEdge}^{2} \\
\kappa_{3\ThreeClawGreenGreenGreen}^{ } &= \mu_{3\ThreeClawGreenGreenGreen}^{ } - 3\mu_{2\TwoGreenPurpleGreen}^{ }\mu_{1\OnePurpleToGreenEdge}^{ } + 2\mu_{1\OnePurpleToGreenEdge}^{3} 
%
%
\end{align}
\endgroup

\todoscience{expression for bipartite clustering analysis}

\subsection{Weighted networks}
\label{SI:networkcumulantslist_WeightedNetworks}

For weighted networks, each instance of a subgraph is counted with weight equal to the product of its edge weights (see \figreftextSI~\ref{fig:SchemeWeightsWedge} and \SIreftextInSI~\ref{SI:graphcumulantsadd})), but the normalization to moments and conversion to cumulants remain the same as in the unweighted case.  
However, the expressions for computing the disconnected counts from the connected counts requires a slight modification.  
For example, the counts of two weighted edges that do not share any node now includes a \mbox{second-order} term related to the square of the edge weights.

\subsubsection{Graph moments}

\begingroup
\allowdisplaybreaks
\begin{align}
\mu_{1\oneedge}^{ } &= \frac{1}{{n\choose 2}} \sum_{0\leq i<j \leq n}^{ } w_{i\!j}^{ }\\
\mu_{2\tworeciprocal}^{ } &= \frac{1}{{n\choose 2}} \sum_{0\leq i<j \leq n}^{ } w_{i\!j}^{2}\\
\mu_{2\twowedge}^{ } &= \frac{1}{3{n\choose 3}} \sum_{0\leq i<j<k \leq n}^{ } \!\!\Big( w_{i\!j}^{ }w_{j\!k}^{ } + w_{j\!k}^{ }w_{k\!i}^{ } + w_{k\!i}^{ }w_{i\!j}^{ } \Big)\\
\mu_{2\twoparallel}^{ } &= \frac{1}{3{n\choose 4}} \Bigg( \frac{1}{2}\bigg( \mu_{1\oneedge}^{ } {n \choose 2} \bigg)^2 - \frac{1}{2} \bigg( \mu_{2\tworeciprocal}^{ } {n \choose 2}\bigg) - \bigg( \mu_{2\twowedge}^{ } 3 {n \choose 3}\bigg) \Bigg)\\
\mu_{3\threereciprocal}^{ } &= \frac{1}{{n\choose 2}} \sum_{0\leq i<j \leq n}^{ } w_{i\!j}^{3}\\
\mu_{3\threereciprocalwedge}^{ } &= \frac{1}{6{n\choose 3}} \sum_{0\leq i<j<k \leq n}^{ } \!\!\Big( w_{i\!j}^{2}w_{j\!k}^{ } + w_{j\!k}^{2}w_{k\!i}^{ } + w_{k\!i}^{2}w_{i\!j}^{ } + w_{i\!j}^{ }w_{j\!k}^{2} + w_{j\!k}^{ }w_{k\!i}^{2} + w_{k\!i}^{ }w_{i\!j}^{2} \Big)
\end{align}
\endgroup

\subsubsection{Graph cumulants}

\begingroup
\allowdisplaybreaks
\begin{align}
\kappa_{1\oneedge}^{ } &= \mu_{1\oneedge}^{ } \\
\kappa_{2\tworeciprocal}^{ } &= \mu_{2\tworeciprocal}^{ } - \mu_{1\oneedge}^{2} \\
\kappa_{2\twowedge}^{ } &= \mu_{2\twowedge}^{ } - \mu_{1\oneedge}^{2} \\
\kappa_{2\twoparallel}^{ } &= \mu_{2\twoparallel}^{ } - \mu_{1\oneedge}^{2} \\
\kappa_{3\threereciprocal}^{ } &= \mu_{3\threereciprocal}^{ } - 3\mu_{2\tworeciprocal}^{ }\mu_{1\oneedge}^{ } + 2\mu_{1\oneedge}^{3} \\
\kappa_{3\threereciprocalwedge}^{ } &= \mu_{3\threereciprocal}^{ } - 2\mu_{2\twowedge}^{ }\mu_{1\oneedge}^{ } - \mu_{2\tworeciprocal}^{ }\mu_{1\oneedge}^{ } + 2\mu_{1\oneedge}^{3} 
\end{align}
\endgroup


%
%

\end{document}